\documentclass[11pt]{article}
\input amssymb.sty
\input amssym.def
\input amssym
\input epsf

\newtheorem{theorem}{Theorem}[section]
\newtheorem{lemma}[theorem]{Lemma}
\newtheorem{proposition}[theorem]{Proposition}
\newtheorem{corollary}[theorem]{Corollary}
\newtheorem{conjecture}[theorem]{Conjecture}
\newtheorem{definition}[theorem]{Definition}

\newtheorem{theorem-construction}[theorem]{Theorem--Construction}
\newtheorem{lemma-construction}[theorem]{Lemma--Construction}

\begin{document}

\newcommand{\Z}{{\Bbb Z}}
\newcommand{\R}{{\Bbb R}}
\newcommand{\Q}{{\Bbb Q}}
\newcommand{\C}{{\Bbb C}}
\newcommand{\lra}{\longrightarrow}
\newcommand{\lms}{\longmapsto}
\newcommand{\AAA}{{\Bbb A}}
\newcommand{\Alt}{{\rm Alt}}
\newcommand{\wg}{\wedge}
\newcommand{\ol}{\overline}
\newcommand{\CP}{{\Bbb C}P}
\newcommand{\bwg}{\bigwedge}
\newcommand{\caL}{{\cal L}}
\newcommand{\PP}{{\Bbb P}}
\newcommand{\HH}{{\Bbb H}}
\newcommand{\LL}{{\Bbb L}}
\newcommand{\hra}{{\hookrightarrow}}

\begin{titlepage}
\title{Galois symmetries of fundamental groupoids
 and noncommutative geometry}
\author{A. B.  Goncharov}
\date{}
\end{titlepage}
%\stepcounter{page}
\maketitle
\tableofcontents

\section{Introduction}

{\bf Abstract}. We define a Hopf algebra of motivic iterated 
integrals on the line and prove an explicit  formula for  the 
coproduct $\Delta$ in this Hopf algebra. 
We show that this formula encodes the group law 
of the automorphism group  of a certain noncommutative variety. 
We relate the 
coproduct $\Delta$ with the coproduct in the 
Hopf algebra of decorated rooted plane trivalent trees,  
which is a plane decorated  version of the one 
defined by Connes and  Kreimer [CK]. 
As an application we derive explicit formulas for the coproduct 
in the motivic multiple polylogarithm Hopf algebra. 
These formulas play a key role in the mysterious correspondence between 
the structure of the motivic fundamental group of ${\Bbb P}^1 - (\{0,
\infty\}\cup  \mu_N)$, 
where  $\mu_N$ is the group of all 
$N$-th roots of unity, 
and modular varieties for ${GL_m}$ ([G1--2]). 

In Chapter 7 we discuss some general principles 
relating Feynman integrals and  mixed motives. They are  suggested by 
Chapter 4  and the Feynman integral 
approach for multiple polylogarithms on curves  given in [G2]. 
Chapter 8 contains background material. 

{\bf 1. The Hopf algebra of motivic iterated integrals}. 
Consider the iterated integral 
\begin{equation} \label{7.23.02.5}
{\rm I}_{\gamma}(a_0; a_1, ..., a_n; a_{n+1}):= (2\pi i)^{-n} \cdot 
\int_{\Delta_{n, \gamma}} \frac{dt_1}{t_1-a_1} \wedge 
\frac{dt_2}{t_2-a_2} \wedge ... \wedge \frac{dt_n}{t_n-a_n}
\end{equation}
Here $\gamma$ is a path from $a_0$ to $a_{n+1}$ in 
$\C - \{a_1 \cup ... \cup a_n\}$, and integration is over a 
simplex $\Delta_{n, \gamma}$ consisting of all ordered $n$--tuples of points 
$(t_1, ..., t_n)$ on  $\gamma$.  
We assume  no restrictions 
on the points $a_1, ..., a_{n} \in \C$. So
 the iterated integral (\ref{7.23.02.5}) can be divergent, 
and in this case it has to be regularized. The reguralization 
depends on the choice of tangent vectors in the points $a_0$ and $a_{n+1}$. 
Below we will assume that these tangent vectors are $\partial/\partial t$, 
where $t$ is the standard coordinate on ${\Bbb A}^1$.

The numbers  (\ref{7.23.02.5}) are periods of $\Q$-rational framed 
Hodge-Tate structures. A Hodge-Tate structure is a 
 mixed Hodge structure such that the Hodge numbers $h^{p,q}=0$ if $p \not = q$. 
The framing is an additional data which, together with a splitting of
the weight filtration, allows   to consider  
a specific period of a mixed Hodge structure, which is simply 
a complex number. If we want to keep the information 
about this period only we are led to an  equivalence relation on the set of all framed 
Hodge-Tate structures. 
The equivalence classes form a commutative graded Hopf algebra ${\cal H}_{\bullet}$ 
over $\Q$. 
 We review these definitions in the Appendix. 
The product structure on ${\cal H}_{\bullet}$ 
is compatible with the product of the periods. The coproduct 
is something really new: it is invisible on the level of numbers.

Let us choose an embedding $\overline \Q \hookrightarrow \C$. 
Let us assume that the parameters $a_i$ of 
the iterated integral (\ref{7.23.02.5}) are algebraic numbers. Then 
we can do even better, and upgrade (\ref{7.23.02.5}) 
to a  framed mixed Tate motive  over $\overline \Q$,
called {\it motivic iterated integral}:
\begin{equation} \label{7.26.02.11}
{\rm I}^{\cal M}(a_0; a_1, ..., a_n; a_{n+1})\in {\cal A}_n(\overline \Q)
\end{equation}
 By its very definition it lies in a  commutative, graded Hopf algebra 
${\cal A}_{\bullet}(\overline \Q)$ with a
coproduct $\Delta$, defined in the Appendix.  

More precisely, let $F$ be a number field. Then there is a graded, commutative 
Hopf algebra ${\cal A}_{\bullet}(F)$, the fundamental Hopf algebra 
of the abelian category ${\cal M}_T(F)$ 
of mixed Tate motives over $F$, see the Appendix for the background. It is 
isomorphic to the 
Hopf algebra of 
regular functions on the unipotent part of the 
motivic Tate Galois group of $\overline \Q$ (loc. cit.). 

The Hopf algebra ${\cal A}_{\bullet}(F)$ depends functorially on $F$, 
and have the following structure. Let $T(V_{\bullet})$ denotes 
the tensor algebra of a graded $\Q$-vector space $V_{\bullet}$. 
It has a natural commutative graded Hopf algebra structure. Then there is an isomorphism
$$
{\cal A}_{\bullet}(F) \stackrel{\sim}{=} T\Bigl(\oplus_{n\geq 1}K_{2n-1}(F)\otimes \Q\Bigr)
$$
where $K_{2n-1}$ on the right sits in the degree $n$. 

 Given an embedding $\sigma: F \hra \C$ of 
$F$ into $\C$, we define  in Chapter 8.5 a filtered algebra ${\cal P}^{\sigma}(F)$ over $\Q$ 
of periods 
of mixed Tate motives over $F$ in the Hodge realization provided by $\sigma$:
$$
\Q = {\cal P}^{\sigma}_{\leq 0}(F) \subset {\cal P}^{\sigma}_{\leq 1}(F) \
\subset {\cal P}^{\sigma}_{\leq 2}(F)  \subset ...; \qquad {\cal P}^{\sigma}(F)= \cup 
{\cal P}^{\sigma}_{\leq n}(F)
$$
Let ${\cal P}^{\sigma}_{\bullet}(F)$ be its associate graded. We prove in Theorem 
\ref{6.14.04.5} that 
there is canonical surjective homomorphism of graded commutative algebras
\begin{equation} \label{8.04.03.3}
p_{\sigma}: {\cal A}_{\bullet}(F) \lra {\cal P}^{\sigma}_{\bullet}(F)
\end{equation}

If the parameters of the iterated integral (\ref{7.23.02.5}) are in the subfield 
$\sigma(F) \subset \C$, one can prove ([G3]) that the projection of the value of integral 
(\ref{7.23.02.5}) to ${\cal P}^{\sigma}_{n}(F)$ does not depend on the
 path $\gamma$, i.e. on the monodromy. Thus the values of 
the iterated integral (\ref{7.23.02.5}) provide a well 
defined element 
$$
\overline {\rm I}(a_0; a_1, ..., a_n; a_{n+1}) \in {\cal P}^{\sigma}_{\bullet}(F), 
\quad a_i \in \sigma(F) \subset \C
$$ 

\begin{theorem} \label{6.17.04.1} Suppose that $a_i$ are elements of a number field $F$. 
Then there exist motivic iterated integrals
\begin{equation} \label{6.17.04.2}
{\rm I}^{\cal M}(a_0; a_1, ..., a_n; a_{n+1})\in {\cal A}_n(F)
\end{equation}
such that 
$$
p_{\sigma}({\rm I}^{\cal M}(a_0; a_1, ..., a_n; a_{n+1})) = \overline {\rm I}(\sigma(a_0); 
\sigma(a_1), ..., \sigma(a_n); \sigma(a_{n+1})) 
$$ 
\end{theorem}

{\bf Example}. Let $n=1$. Then 
$$
{\rm I}(a,b,c) = (2\pi i)^{-1}\int_a^c\frac{dt}{t-b} = 
(2\pi i)^{-1}\log\frac{c-b}{a-b}
$$
This integral diverges if $a=b$ or $b=c$. To regularize it we calculate 
the asymptotic expansion when $\varepsilon \to 0$ of the integral 
 $$
(2\pi i)^{-1}\int_{a+\varepsilon}^{c+\varepsilon}\frac{dt}{t-b} = (2\pi i)^{-1}
\log\frac{c-b +\varepsilon}{a-b +\varepsilon}, 
$$
 getting a polynomial in $\log 
\varepsilon $, 
and take its free term as the regularized value of the integral. We get 
the following answer. 
Let
\begin{equation} \label{6.17.04.10}
\widetilde r(a,b,c):= \left\{ \begin{array}{ll} 
\frac{(c-b)}{(a-b)} & 
\mbox{if $a \not = b, b \not = c$}
\\ c-b & \mbox{if $a=b$, but $b \not = c$}\\
(a-b)^{-1} & \mbox{if $b = c$, but $a\not = b$}\\
1 & \mbox{if $a = b = c$}
\end{array}\right.
\end{equation}
Then 
$$
{\rm I}(a,b,c) = (2\pi i)^{-1}\log \widetilde r(a,b,c)
$$
The right hand side provides a well defined element $\overline {\rm I}(a,b,c) \in \C/\Q$. 
There is a canonical isomorphism ${\cal A}_1(F) = F^* \otimes \Q$, and 
$$
{\rm I}^{\cal M}(a,b,c) =  \widetilde r(a,b,c)
$$

One of our  results is  the following explicit formula 
for the coproduct of the elements (\ref{6.17.04.2}).

\begin{theorem} \label{8.3.02.1}
The coproduct $\Delta$ is computed by the formula 
\begin{equation} \label{CP2*q}
\Delta {\rm I}^{\cal M}(a_0; a_1, a_2, ..., a_n; a_{n+1})= 
\end{equation}
$$
\sum_{0 = i_0 < i_1 < ... < i_k < i_{k+1} = n+1} 
{\rm  I}^{\cal M}(a_0; a_{i_1}, ..., a_{i_k}; a_{n+1}) \otimes \prod_{p =0}^k
{\rm  I}^{\cal M}(a_{i_{p}}; a_{i_{p}+1}, ..., a_{i_{p+1}-1}; a_{i_{p+1}})
$$ Here $0 \leq k \leq n$ and $a_i \in \overline \Q$. 
\end{theorem}
We prove a similar result for the Hodge and $l$-adic analogs of the
motivic 
iterated integrals, see Theorem \ref{7.24.02.12} for a precise
statement. Another proof in the Hodge set-up was given in
Sections 5-6 of [G3]. 

The terms in the formula (\ref{CP2*q}) 
are in one--to--one correspondence with the  subsequences 
\begin{equation} \label{7.23.2.6}
\{a_{i_1}, ..., a_{i_k}\} \subset \{a_{1}, ..., a_{n}\}
\end{equation}
If we locate the ordered sequence $\{a_{0}, ..., a_{n+1}\}$ 
on a semicircle 
then the terms in 
(\ref{CP2*q}) correspond to 
the polygons with vertices at the points $a_i$, 
containing $a_0$ and $a_{n+1}$, inscribed into the semicircle. 
\begin{center}
\hspace{4.0cm}
\epsffile{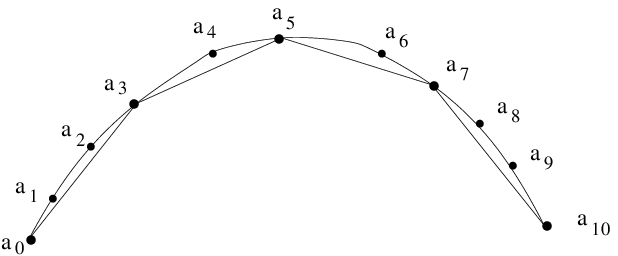}
\end{center}
The  picture illustrates the term 
$$
{\rm  I}^{\cal M}(a_0; a_{3},  a_{5},   a_{7}; a_{10}) \otimes 
{\rm  I}^{\cal M}(a_{0}; a_{1}, a_{2}; a_{3}){\rm  I}^{\cal M}(a_{3}; a_{4}; a_{5})
{\rm  I}^{\cal M}(a_{5}; a_{6}; a_{7}){\rm  I}^{\cal M}(a_{7}; a_{8}, a_{9}; a_{10})
$$

{\bf Example}. For $m=2$ formula  (\ref{CP2*q}) gives
$$
\Delta {\rm I}^{\cal M}(a_0; a_1, a_2; a_3) \quad = 
\quad 1 \otimes {\rm I}^{\cal M}(a_0; a_1, a_2; a_3) +
$$
$$
{\rm I}^{\cal M}(a_0; a_1; a_3) \otimes {\rm I}^{\cal M}(a_1; a_2; a_3) +  
{\rm I}^{\cal M}(a_0; a_2; a_3) \otimes {\rm I}^{\cal M}(a_0; a_1; a_2) +
{\rm I}^{\cal M}(a_0; a_1, a_2; a_3) \otimes 1 
$$
The $m=3$ case is worked out at the end of Chapter 3.

Our proof penetrates some algebraic structures 
staying behind this formula. 
They include automorphism groups of 
certain noncommutative varieties discussed in Chapters 2 and 3, and  
a Hopf algebra of decorated rooted plane trees in Chapter 4.

{\bf 2. Why do we care about motivic iterated integrals?} In other words, 
what do we gain by  upgrading the iterated integrals (\ref{7.23.02.5}) 
to the elements (\ref{6.17.04.2}) of the motivic Hopf algebra 
${\cal A}_{\bullet}(F)$? 

First, by doing this we {\it conjecturally} do not lose any information. 
\begin{conjecture} \label{8.04.03.2} The map (\ref{8.04.03.3}) 
is an isomorphism of $\Q$--vector spaces. 
\end{conjecture}

Our point is that 
it is much easier to deal with  motivic iterated integrals than with  numbers! 

Here is what  we gain:

A) Unlike the numbers, the motivic iterated integrals form a Hopf algebra.

B) Working with the more  sophisticated   motivic iterated integrals we 
eliminated the transcendental 
aspect of the problem.

C) Motivic iterated integrals serve as a bridge between their Hodge (i.e. analytic) and $l$-adic 
(i.e. arithmetic)  realizations. This 
allows to use arithmetic insights in analytic problems and vice versa.

Here is a simplest example. 
The 
Leibniz formula expresses special values of the Riemann $\zeta$--function via 
iterated integrals:
$$
\zeta(m) 
%= (2\pi i)^{m}{\rm I}(0; 1, \underbrace{0,  ...,0}_{m-1}; 1) 
= \int_{0 \leq t_1 \leq ... 
\leq t_m \leq 1} \frac{dt_1}{1-t_1}\wedge 
\frac{dt_2}{t_2}\wedge ... \wedge 
\frac{dt_m}{t_m}
$$
According to Euler $(2\pi i)^{-2k}\zeta(2k) \in \Q$. 
In accordance with this one has $\zeta^{\cal M}(2k) =0$. 
Nobody can  prove so far that the numbers  $\zeta(2k+1)$ 
are linearly independent over $\Q$, e.g.   
that ${ \zeta}(5) \not \in \Q$. 
However it is easy to show that the motivic elements
$$
\zeta^{\cal M}(2k+1) \in {\cal A}_{2k+1}(\Q)
$$ 
 are non-zero. Moreover they are  linearly independent over $\Q$. 
Indeed, they belong to  components of different degrees in 
${\cal A}_{\bullet}(\Q)$.

 Another example is worked out in 
Chapter 6.7, where  the 
motivic double  zeta's are investigated.  
Recall that the double zeta values $\zeta(m,n)$ can be presented 
both as power series and iterated integrals: 
$$
\zeta(m, n)= \sum_{0< k_1 < k_2}\frac{1}{k_1^mk_2^n} = {\rm I}(0; \underbrace{1,
0, ..., 0}_m, \underbrace{1, 0, ..., 0}_n; 1)
$$ 
They satisfy the (regularized) double shuffle
relations, obtianed by expressing the product of two classical 
$\zeta$-values by using either the power series or the iterated integral
presentations. The iterated integral presentation leads to the definition of 
motivic double zeta's $\zeta^{\cal M}(m, n)$. 
It follows from  [G8] 
that the motivic double $\zeta$'s satisfy the same double shuffle
relations. We will reprove in Chapter 6.7 a version of this result, and prove that 
{\it any relation between the motivic double $\zeta$'s is a
  consequence of the motivic double shuffle relations}. 

I  think that an understanding 
of the transcendental aspects of the iterated integrals 
(\ref{7.23.02.5}) is impossible without 
  investigation of the corresponding motivic objects. 

{\bf 3. Motivic iterated integrals unramified outside of a finite set $S$ of prime ideals}. 
Given the parameters $a_i \in F$, one can describe rather precisely the 
subspace in ${\cal A}_{n}(F)$ where the motivic iterated integral (\ref{6.17.04.2}) land. 
Namely, let $S$ be a collection of primes in a number field $F$, and ${\cal O}_{F, S}$ 
the ring of $S$-integers in $F$. Then there exists the 
category of mixed Tate motives over ${\cal O}_{F, S}$, 
([DG]). Let 
${\cal A}_{\bullet}({\cal O}_{F, S})$ be its fundamental Hopf algebra (Section 8.4). 
We prove the following result (Theorem \ref{6.12.04.1}).  
Recall $\widetilde r(a, b, c)$ from (\ref{6.17.04.10}).  

\begin{theorem} \label{6.12.04.1a}
Let $a_0, ..., a_{n+1}$ be elements of a number field $F$. 
If 
for any $0 \leq i < j < k \leq n+1$ one has 
$
\widetilde r(a_i, a_j, a_k) \in {\cal O}_{F,S}^*
$, then
\begin{equation}\label{6.12.04.4}
{\rm I}^{\cal M}(a_0; a_1, ..., a_n; 
a_{n+1}) \in {\cal A}_n({\cal O}_{F, S})
\end{equation}
\end{theorem}

\begin{corollary} \label{6.12.04.5}
Suppose that $a_0, ..., a_{n+1}\in \{\zeta^{a}_N\} \cup \{0\}$, where $\zeta_N$ is a primitive root of unity. 
Then 
\begin{equation}\label{6.17.04.11}
{\rm I}^{\cal M}(a_0; a_1, ..., a_n; 
a_{n+1}) \in {\cal A}_n(\Z[\zeta_N][\frac{1}{N}])
\end{equation}
\end{corollary}
In particular $\zeta^{\cal M}(n_1, ..., n_m) \in {\cal A}_w(\Z)$, where $w = n_1 + ... + n_m$.

Indeed, $\zeta_N^a - \zeta_N^b = \zeta_N^a(1 - \zeta_N^{b-a}) \in  
(\Z[\zeta_N][\frac{1}{N}])^*$. So the corollary follows from Theorem \ref{6.12.04.1a}. 

{\bf Example}. $\zeta^{\cal M}(n_1, ..., n_m) \in {\cal A}_w(\Z)$, where $w = n_1 + ... + n_m$.

The Hopf algebra ${\cal A}_{\bullet}({\cal O}_{F, S})$ has the following structure. 
There exists an isomorphism 
of commutative graded Hopf algebras, see Section 8.5:
$$
{\cal A}_{\bullet}({\cal O}_{F, S}) \stackrel{\sim}{=} 
T\Bigl(\oplus_{n>1}K_{2n-1}({\cal O}_{F, S})\otimes {\Q}\Bigr)
$$
where $K_{2n-1}$ sits in the degree $n$. In particular, if $S$ is a finite set, then 
${\cal A}_{n}({\cal O}_{F, S})$ is a finite-dimensional $\Q$-vector space! 
Indeed, 
if $n>1$ then for any $S$ 
$K_{2n-1}({\cal O}_{F, S})\otimes \Q  = K_{2n-1}(F)\otimes \Q $ is finite dimensional 
by Borel's theorem. But $K_{2n-1}({\cal O}_{F, S}) = {\cal O}_{F, S}^*$ is of finite rank only if $S$ is finite. 
This combined with Theorem \ref{6.12.04.1a} give a strong estimate from above on the 
dimension of the subspaces in $\C$  generated by 
the values of iterated integras (\ref{7.23.02.5}),
see Section 6.8.

Here is a typical example. 
Combining Corollary \ref{6.12.04.5} with (\ref{6.14.04.4df}), we get an estimate from 
above on the $\Q$-vector space ${\cal Z}_n(\mu_N \cup \{0\})$ spanned by the iterated integrals 
${\rm I}(a_0; a_1, ..., a_n; 
a_{n+1})$, where  $a_0, ..., a_{n+1}\in \{\zeta^{a}_N\} \cup \{0\}$, modulo 
similar integrals of weight $< n$: 
\begin{equation}\label{6.12.04.xf}
{\rm dim}{\cal Z}_n(\mu_N \cup \{0\}) \leq {\rm dim}
T_{(n)}\Bigl(\oplus_{n>0}K_{2n-1}(\Z[\zeta_N][\frac{1}{N}])\Bigr)
\end{equation}
Here on the right stays the dimension of the degree $n$ part of the tensor algebra. 
A bit different way to get the same estimate see in [DG].

To clarify the relationship between the numbers and their more
sofisticated 
motivic counterparts we discuss below 
the Hodge and $l$-adic versions of iterated integrals. 

{\bf 4. The Hodge iterated integrals}. By definition 
(see the Appendix) an $n$-framing on a Hodge-Tate structure $H$ is a 
choice of nonzero morphisms of Hodge-Tate structures
$$
 v: \Q(0) \lra {\rm gr}_{0}^WH; \qquad f: {\rm gr}_{-2m}^WH \lra \Q(m)
$$
Consider the coarsest equivalence relation $\sim$ on the set of all 
$n$-framed Hodge-Tate structures such that $(H_1, v_1, f_1) \sim (H_2, v_2, f_2)$ if there 
is a morphism of Hodge-Tate structures $H_1 \to H_2$ respecting the frames. 
The equivalence classes 
of $n$-framed Hodge-Tate structures form an abelian group ${\cal H}_{n}$, and 
${\cal H}_{\bullet}:= \oplus_{n \geq 0}{\cal H}_{n}$ has a commutative 
graded Hopf algebra structure.

The framed Hodge-Tate structure corresponding to the iterated integrals 
(\ref{7.23.02.5}) arises as follows. Let $S:= \{z_1, ..., z_m\}$. Denote by 
$v_a, v_b$ the tangent vectors at the points $a,b$ dual to the canonical 
differential $dt$ on $\C$. Then the pronilpotent completion 
\begin{equation} \label{3.6.01.7}
{\cal P}^{\cal H}(\C  - S; v_a, v_b)
\end{equation}
of the topological torsor of paths   from 
$v_a$ to $v_b$ (see s. 3.1) is equipped with the structure of a projective limit of 
 Hodge-Tate structures ([D]). In particular 
 it carries a weight filtration $W_{\bullet}$ such that 
\begin{equation} \label{3.6.01.11}
{\rm gr}_{-2m}^W{\cal P}^{\cal H}(\C - S; v_a, v_b)  \stackrel{\sim}{=} \otimes^m 
H_1(\C  - S)
\end{equation}
The forms $(2\pi i)^{-1}d\log (t-z_i)$ provide  a basis of 
$H^1(\C - S)$. We define 
\begin{equation} \label{3.6.01.12}
{\rm I}^{\cal H}(a; z_1, ..., z_m; b) \in {\cal H}_m
\end{equation}
as  the 
Hodge-Tate structure (\ref{3.6.01.7}) framed by 
$1 \in \Q\stackrel{\sim}{=} {\rm gr}_{0}^W{\cal P}^{\cal H}(\C - S; v_a, v_b)$ and 
\begin{equation} \label{3.6.01.5}
(2\pi i)^{-1} d\log(t-z_1)  \otimes
... \otimes (2\pi i)^{-1} d\log(t-z_m) 
 \in {\rm gr}_{-2m}^W{\cal P}^{\cal H}(\C - S; v_a, v_b)^{\vee}
\end{equation}
If we choose in addition a splitting of the weight filtration in
(\ref{3.6.01.7}), one  gets a complex number,  period, out of (\ref{3.6.01.12}),
given by the usual iterated integral. 

To see why the coproduct appears,  
let us discuss  the $l$-adic case. 

{\bf 5. $l$-adic iterated integrals as functions on the Galois
  group}. Let $F$ be a field that does not contain all $l^{\infty}$-th
  roots of unity. For example, take a number field $F$. Let us choose a finite subset
  $S =\{z_1, ..., z_m\} \subset F$. The pro-$l$ completion of
  the topological torsor of path from $v_a$ to $v_b$ (s. 3.1) 
\begin{equation} \label{cv}
{\cal P}^{(l)}({\Bbb A}^1 -S; v_a, v_b)
\end{equation}
 is a Galois module.  We restrict it to the subgroup 
${\rm Gal}(\overline \Q/F(\zeta_{l^{\infty}}))$. Let us
make the following choices:

i) Choose a prime $p$ such that (\ref{cv})  is unramified at $p$, and choose a
  Frobenious element $F_p$ at $p$. 

ii) Choose a generator
  $\{\zeta_{l^n}\}$ of the Tate module.  

Then,  splitting  the weight filtration on
  (\ref{cv}) by the  eigenspaces of $F_p$, we get an isomorphism
$$
\prod_{n \geq 0} 
{\rm Hom}(\Q_l(n), H_1({\Bbb
  A}^1 -S, \Q_l)^{\otimes n}) \stackrel{\sim}{\lra} {\cal P}^l({\Bbb A}^1 -S; v_a, v_b)  
$$
$H_1({\Bbb
  A}^1 -S, \Q_l)$ is identified with a direct sum of the Tate modules
$\Q_l(1)$ indexed by the points $z_i$. Thanks to ii), it  has a natural basis, denoted $d\log(t-z_i)$. Take the
dual basis.  It provides, for each $m \geq 0$,  a basis in 
the weight $[0, -m]$ quotient of (\ref{cv}). In particular there is
canonical isomorphism $\Q_l = {\rm Gr}^W_0{\cal P}^{(l)}({\Bbb A}^1 -S;
v_a, v_b)$. We denote by ${ 1}$ the image of $1\in \Q_l$. 
\begin{definition} \label{2.9.04.1} 
  ${\rm I}_{F_p}^{(l)}(a; z_1, ..., z_m; b)$ is the function on ${\rm Gal}(\overline
  \Q/F(\zeta_{l^{\infty}}))$ given by the matrix element
\begin{equation} \label{2.10.04.1}
g \lms <d\log(t-z_1) \otimes ... \otimes d\log(t-z_m) | g | 1>
\end{equation} 
\end{definition} 

The $\Q_l$-valued functions on the Galois group form a
Hopf algebra. Its coproduct $\Delta^{(l)}$ is provided by the group
multiplication. 
We prove that the coproduct  $\Delta^{(l)}{\rm I}_{F_p}^{(l)}(a; z_1,
..., z_m; b)$ is calculated by the same formula (\ref{CP2*q}). 
It follows that when $a, b, z_1, ..., z_m$ run
through the elements of a finite set $S$, the functions
(\ref{2.10.04.1}) span over $\Q_l$ a commutative graded Hopf algebra
${\cal H}^{(l)}_{\bullet}(S)$.  

{\bf Example}. By Kummer's theory an element $a \in F^*$
gives rise to a $\Z_l$-valued function $\chi_a(g)$ on ${\rm Gal}(\overline \Q/F(\zeta_{l^{\infty}}))$ given
by 
$
\zeta^{\chi_a(g)}_{l_n} := g(a^{1/l^n})/a^{1/l^n}
$. One has $\Delta^{(l)} \chi_a = \chi_a \otimes 1 + 1 \otimes
\chi_a$. We have 
$$
\chi_a = {\rm I}^{(l)}(1; 0, a) = {\rm I}^{(l)}(0; 0, a) 
$$
It turns out to be independent of a choice of Frobenious or  splitting. 

On the other hand, the Galois module (\ref{cv}), framed by $1$ and
$d\log(t-z_1) \otimes ... \otimes d\log(t-z_m)$, provides 
an element ${\rm I}^{(l)}(a; z_1, ..., z_m; b)$ of the commutative graded algebra of framed mixed Tate
$l$-adic representations of the Galois group, see the Appendix. 
They are the $l$-adic counterparts of the elements ${\rm I}^{\cal H}(a; z_1, ...,
z_m; b)$. The functions ${\rm
  I}^{(l)}_{F_p}(a; z_1, ..., z_m; b)$ 
on the Galois group are the $l$-adic counterparts of periods.  One can
avoid making the choise ii), getting $\Q_l(n)$-valued functions on the
Galois group. However, just as in the analyitc case, without the  choice i) we
can not get functions on the Galois group. The
elements ${\rm
  I}^{(l)}$ 
are 
canonical: their definition does not depend on the choices i) and ii). They span a Hopf algebra, which is canonically 
isomorphic to ${\cal H}^{(l)}_{\bullet}(S)$, see the Appendix.

The definition of the motivic iterated integrals as
framed objects is similar.  
It uses the motivic version of 
the pronilpotent torsor of paths
(\ref{3.6.01.7}). 

A version of Definition \ref{2.9.04.1} is obtained if one does not use the choice ii).
Then formula (\ref{2.10.04.1}) yelds 
a $\Q_l(n)$-valued function on the Galois group ${\rm Gal}(\overline \Q/F)$. 

{\bf 6. The algebraic structures underlying  formula (\ref{CP2*q})}. 
The integrals (\ref{7.23.02.5}) 
satisfy the following basic properties ([Chen]):

{\it The shuffle product formula}. Let $\Sigma_{m,n}$ be the set of all shuffles 
of the ordered sets $\{1, ..., m\}$ and $\{m+1, ..., m+n\}$. Then 
$$
{\rm I}_{\gamma}(a; z_1, ..., z_m; b)\cdot {\rm I}_{\gamma}
(a; z_{m+1}, ..., z_{m+n}; b) = 
\sum_{\sigma \in \Sigma_{m,n}}{\rm I}_{\gamma}(a; z_{\sigma(1)}, ..., z_{\sigma(m+n)}; b)
$$

{\it The path composition formula}. If  $\gamma = \gamma_1 \gamma_2 $, where $\gamma_1$ is a path from $a$ to $x$, 
and $\gamma_2$ is a path from $x$ to $b$ in 
$\C - \{z_1 \cup ... \cup z_n \cup a \cup b\}$, then 
$$
{\rm I}_{\gamma}(a; z_1, ..., z_m; b) = 
\sum_{k=0}^{m}{\rm I}_{\gamma_1}(a; z_1, ..., z_k; x) \cdot 
{\rm I}_{\gamma_2}(x; z_{k+1}, ..., z_m; b) 
$$
 
Let $S$ be an arbitrary set. In Chapter 2 we  define a commutative, graded 
 Hopf algebra ${\cal I}_{\bullet}(S)$.  As a commutative algebra it is 
generated  
by symbols ${\Bbb I}(s_0; s_1, ..., s_n; s_{n+1})$ with $s_i \in S$ 
 subject to the relations as for $I_{\gamma}$ above. 
We use  formula (\ref{CP2*q}) to {\it define} 
the coproduct in ${\cal I}_{\bullet}(S)$. 
If $S$ is a subset of points of the affine line ${\Bbb A}^1$ the Hopf algebra 
${\cal I}_{\bullet}(S)$ 
reflects the   basic properties 
of the iterated integrals on ${\Bbb A}^1_S:= {\Bbb A}^1 -S$ between the  
tangential base points at $S$. 
To show that ${\cal I}_{\bullet}(S)$ is indeed a Hopf algebra
we interpret it as the algebra of regular functions on a certain 
pro--unipotent group scheme. 
Namely, let $\Gamma_S$ be the graph whose vertices are elements of $S$, and 
every two vertices are connected by a unique edge. 
Let  ${P}(S)$ be the 
path algebra 
 of this graph. Its basis over $\Q$ is formed by the paths in the graph.
The
 product of two basis elements is given by the composition of paths, 
provided that the end of the first
 coincides with the beginning of the second. Otherwise the product is zero. 

\begin{center}
\hspace{4.0cm}
\epsffile{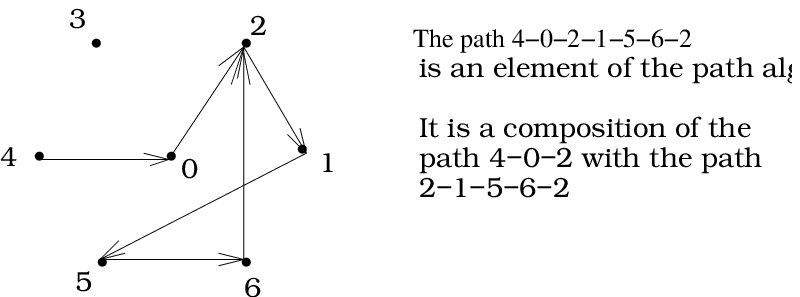}
\end{center}
It  has some additional structures: two algebra structures $\circ$ and $\ast$,  
plus a coproduct  $\Delta$. 
Let  
$G(S)$ be  
the  group  of automorphisms  
 of  ${P}(S)$ 
preserving these structures, and acting as the identity on 
${P}_+(S)/{P}^2_+(S)$.  Here ${P}_+(S)$ 
is given by paths of positive length. 
In other words it 
is the group of automorphisms of 
the corresponding noncommutative variety preserving 
some natural structures on it. 
We prove in Theorem \ref{7.18.02.2} that 
$$
G(S) = {\rm Spec}({\cal I}_{\bullet}(S))
$$

In Chapter 3 we show that the  algebra ${P}(S)$
is provided by the motivic fundamental groupoid 
\begin{equation} \label{8.3.02.2}
{\cal P}^{\cal M}({\Bbb A}^1_S, S)
\end{equation}
of paths on the affine line punctured at $S$,  
between the chosen tangential base points at  $S$. 
We consider (\ref{8.3.02.2}) 
as a pro--object in the abelian category 
of mixed Tate motives over a number field or in one of the realization
categories, 
see the 
Appendix for the background. 
Such a category is equipped with  a canonical fiber 
functor $\omega$. 
The motivic Galois group acts on 
$\omega({\cal P}^{\cal M}({\Bbb A}^1_S, S))$. Let  
$G_{\cal M}(S)$ be its image. 
The motivic 
fundamental groupoid
 carries some additional structures provided by the composition of paths and 
``canonical loops'' near the punctures. Using these structures 
we identify 
$\omega({\cal P}^{\cal M}({\Bbb A}^1_S, S))$ 
 with the path algebra ${P}(S)$.   
Therefore $G_{\cal M}(S)$
  is realized as a subgroup of $G(S)$:
\begin{equation} \label{8.3.02.q12}
G_{\cal M}(S) \hookrightarrow G(S)  
\end{equation}
Theorem \ref{8.3.02.1} follows immediately from this. 

In Chapter 4 we show how the collection of all 
  plane trivalent rooted trees decorated 
by the ordered set $\{a_0, a_1, ..., a_m, a_{m+1}\}$ governs the 
fine structure of the motivic object 
(\ref{7.26.02.11}). Namely, in Sections 4.1-4.2 we define
a commutative 
Hopf algebra  ${\cal T}_{\bullet}(S)$   of 
$S$--decorated planar rooted trivalent trees, and  
relate the groups $G(S)$ and ${\rm Spec}({\cal T}_{\bullet}(S))$. 
Precisely, let $\widetilde {\cal I}_{\bullet}(S)$ be the commutative graded algebra 
freely generated by the elements 
${\Bbb I}(s_0; s_1, ..., s_m; s_{m+1})$. 
We use  formula (\ref{CP2*q}) to define the coproduct in 
$\widetilde {\cal I}_{\bullet}(S)$, and show that it is indeed a Hopf algebra.  
Consider the map 
\begin{equation} \label{8.22.01.1}
t: {\Bbb I}(s_0; s_1, ..., s_m; s_{m+1})  \longmapsto  \mbox{\rm formal sum of all 
plane rooted trivalent trees}
\end{equation}
$$
\mbox{\rm decorated by the ordered set 
$\{s_0, s_1, ..., s_m, s_{m+1}\}$}
$$
{\bf Example}. Here is how the map $t$ looks in the two simplest cases
\begin{center}
\hspace{4.0cm}
\epsffile{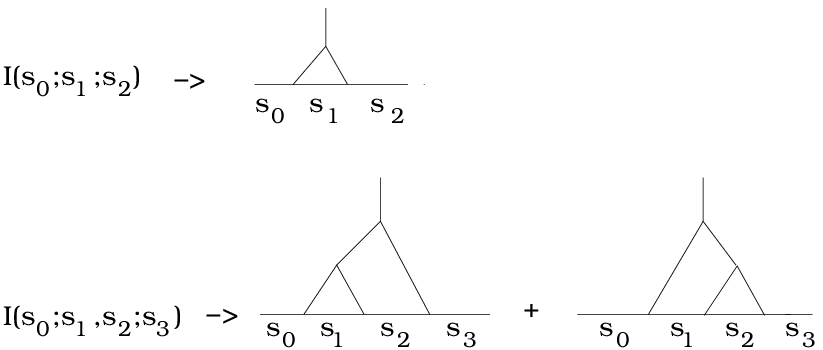}
\end{center}
We prove in Chapter 4 that the map $t$ provides a Hopf algebra homomorphism
$$t: \widetilde {\cal I}_{\bullet}(S) \lra {\cal T}_{\bullet}(S)$$ 

{\bf Remark}. This map does not descend  to a Hopf algebra map from 
$ {\cal I}_{\bullet}(S)$ to ${\cal T}_{\bullet}(S)$. 
I do not know whether there exists a quotient 
of the Hopf algebra ${\cal T}_{\bullet}(S)$ such that 
$t$ maps $ {\cal I}_{\bullet}(S)$ to this quotient. 

{\bf 7. An application: motivic multiple polylogarithm Hopf algebras}. 
Recall the multiple polylogarithms ([G0], [G6]) where 
$n_i \in {\Bbb N}_+$, $|x_i|<1$:
\begin{equation} \label{zhe5}
{\rm Li}_{n_{1},...,n_{m}}(x_{1},...,x_{m})   
\quad = \quad 
\sum_{0 < k_{1} < k_{2} < ... < k_{m} } \frac{x_{1}^{k_{1}}x_{2}^{k_{2}}
... x_{m}^{k_{m}}}{k_{1}^{n_{1}}k_{2}^{n_{2}}...k_{m}^{n_{m}}}
\end{equation}
Theorem 2.2 in [G3] allows to write them as iterated integrals. Namely, if 
 $|x_{i}| < 1$, then setting
$$
a_1 := (x_1...x_m)^{-1}, a_2 := (x_2...x_m)^{-1}, ..., a_m:= x_m^{-1}
$$
we get for a path $\gamma$ inside of the unit disc
$$
{\rm Li}_{n_{1},...,n_{m}}(x_{1},...,x_{m})  = 
(-1)^m {\rm I}_{n_1, ..., n_m}(a_1, ..., a_m):= 
$$
\begin{equation} \label{5*}
(-1)^m {\rm I}_{\gamma}(0; a_1,  \underbrace {0, ... ,0}
_{ n_{1}-1 }, a_2, \underbrace 
{0, ..., 0}_{ n_{2}-1 }, a_3, ..., a_m, \underbrace 
{0, ..., 0}_{ n_{m}-1 }; 1)
\end{equation}
Let ${\cal M}$ be the category of mixed Tate motives over a 
number field $F$ or  one of the mixed Tate categories  
described in section 3.1. 
Let $G$ be a subgroup of the multiplicative group  $F^*$ 
of the field $F$. Upgrading iterated integrals (\ref{5*}) to their 
motivic counterparts (\ref{7.26.02.11}), 
and using identity (\ref{5*}) as a definition, we arrive 
at motivic multiple polylogarithms 
\begin{equation} \label{5s*}
{\rm Li}^{\cal M}_{n_{1},...,n_{m}}(x_{1},...,x_{m}), \qquad x_i \in G \subset F^* 
\end{equation}
Adding the motivic logarithms $\log^{\cal M}(x)$ 
to them, we get a graded, commutative Hopf algebra 
${\cal Z}_{\bullet}^{\cal M}(G)$, see Theorem \ref{7.26.02.111}.
In particular if $G = \mu_N$ 
 we get the cyclotomic Hopf algebra 
${\cal Z}_{\bullet}^{\cal M}(\mu_N)$. According to (\ref{6.17.04.11}), there is an inclusion of 
graded Hopf algebras
$$
{\cal Z}_{\bullet}^{\cal M}(\mu_N) \subset {\cal A}_{\bullet}^{\cal M}(\Z[\zeta_N]\frac{1}{N}])
$$The structure of the 
corresponding cyclotomic Lie coalgebra 
$$
{\cal C}_{\bullet}^{\cal M}(\mu_N):= {\cal Z}_{>0}^{\cal M}(\mu_N)/({\cal Z}_{>0}^{\cal M}(\mu_N))^2
$$
is related to the geometry of  modular varieties,
see [G1-2] and the Section 6.7 for an example.  
 
%The cyclotomic Hopf algebra 
%${\cal Z}_{\bullet}^{\cal M}(\mu_N)$ is the algebra of regular functions on 
%the group scheme $G_{\cal M}(0 \cup \mu_N)$. We show that the semidirect product of 
%${\Bbb G}_m$ and $G_{\cal M}(0 \cup \mu_N)$ 
%is isomorphic to the image of the motivic Galois group acting 
%on the motivic fundamental group
%$\pi_1^{\cal M}({\Bbb G}_m - \mu_N, v_0)$. 
According to Conjecture \ref{1.10.01.3} 
the category of all finite dimensional graded comodules over the Hopf algebra 
${\cal Z}_{\bullet}^{\cal M}(F^*)$ should be equivalent to the conjectural abelian 
category of all mixed Tate motives over $F$. It follows from Theorem
\ref{8.3.02.1} that the 
Hopf algebra ${\cal Z}_{\bullet}^{\cal M}(F^*)$ has an additional
structure: the depth filtration, such that the depth of (\ref{5s*}) is
$m$. 
We suggest a conjectural intrinsic description of the depth filtration 
in Chapter 5.4.

{\bf Acknowledgment}. 
I am grateful to Maxim Kontsevich for several stimulating discussions.  
This paper was written during my stay at  
IHES (Bures sur Yvette) in July of 2002 and at MPI(Bonn) 
in August of 2002. 
I am grateful to IHES and MPI for hospitality and support. 
I was supported by the NSF grant DMS-0099390. 
I would like to thank the referees for careful reading of the manuscript and 
useful comments.

\section{Hopf algebras of iterated integrals}

In this chapter we introduce the Hopf algebras 
${\cal I}_{\bullet}(S)$ and $\widetilde {\cal I}_{\bullet}(S)$ 
and show that they appear naturally 
as the  algebras of functions on  the automorphism groups of 
certain non--commutative varieties. 

{\bf 1. The Hopf algebra $\widetilde {\cal I}_{\bullet}(S)$ and 
${\cal I}_{\bullet}(S)$}. 
Let $S$ be a set. We will define a commutative Hopf 
 algebra  ${\cal I}_{\bullet}(S)$ over $\Q$, graded in an obvious way 
by the integers $n \geq 0$. 

As a {\it commutative $\Q$-algebra},  ${\cal I}_{\bullet}(S)$ is generated 
by the elements 
\begin{equation} \label{10.16.00.1}
{\Bbb I}(s_0; s_1, ..., s_m; s_{m+1}), \quad s_i \in S, \quad m\geq 0
\end{equation}
The generator (\ref{10.16.00.1}) is homogeneous of degree $m$. 
The relations are the following. For any $s_i, a, b, x \in S$ one has:

i) {\it The unit}: for any $a,b \in S$ one has 
$ {\Bbb I}(a;  b) := {\Bbb I}(a; \emptyset ; b) = 1 
$.

ii) {\it The shuffle product formula}: for any $m,n\geq 0$ one has 
$$
{\Bbb I}(a; s_1, ..., s_m; b)\cdot {\Bbb I}(a; s_{m+1}, ..., s_{m+n}; b) = 
\sum_{\sigma \in \Sigma_{m,n}}{\Bbb I}(a; s_{\sigma(1)}, ..., s_{\sigma(m+n)}; b)
$$

iii) {\it The path composition formula}: for any $m \geq 0$ one has 
$$
{\Bbb I}(s_0; s_1, ..., s_m; s_{m+1}) = 
\sum_{k=0}^{m}{\Bbb I}(s_0; s_1, ..., s_k; x) \cdot {\Bbb I}(x; s_{k+1}, ..., s_m; s_{m+1}) 
$$

iv) ${\Bbb I}(a; s_1, ..., s_m; a) = 0$ for $m>0$. 

To define a Hopf algebra structure on ${\cal I}_{\bullet}(S)$
we use the following formula for the coproduct $\Delta$ on the generators:  
for $m\geq 0$ and
$a_0,\dots,a_{m+1}\in S$, we have
\begin{equation} \label{CP2**}
\Delta {\Bbb I}(a_0; a_1, a_2, ..., a_m; a_{m+1})= 
\end{equation}
$$
\sum_{0 = i_0 < i_1 < ... < i_k < i_{k+1} = m+1} 
{\Bbb  I}(a_0; a_{i_1}, ..., a_{i_k}; a_{m+1}) \otimes \prod_{p =0}^k
{\Bbb  I}(a_{i_{p}}; a_{i_{p}+1}, ..., a_{i_{p+1}-1}; a_{i_{p+1}})
$$

Thanks to i) the empty subset and the subset $\{a_1, ..., a_m\}$ contribute 
the terms
$$
1 \otimes {\Bbb  I}(a_0; a_{1}, ..., a_{m}; a_{m+1}) \quad \mbox{and} \quad    
{\Bbb  I}(a_0; a_{1}, ..., a_{m}; a_{m+1}) \otimes 1
$$
The counit is determined by the condition that it kills ${\cal I}_{>0}(S)$.

\begin{proposition} \label{8.15.02.1} In $I_\bullet(S)$ one has, for $m\geq 0$ and
$a_0,\ldots,a_{m+1}\in S$:
$$
{\Bbb  I}(a_0; a_{1}, ..., a_{m}; a_{m+1})  = (-1)^m
{\Bbb  I}(a_{m+1}; a_{m}, ..., a_{1}; a_0)
$$
\end{proposition}

{\bf Proof}. We use induction on $m$. When $m=1$ the path composition 
formula iii) 
plus i) 
gives
$$
0= {\Bbb  I}(a_0; a_{1}; a_{0}) = {\Bbb  I}(a_0; a_{1}; a_{2}) + 
{\Bbb  I}(a_2; a_{1}; a_{0})
$$
Let us assume we proved the claim for all $k<m$. The path composition formula 
with $x:= a_{m+1}$ 
gives
$$
0= {\Bbb  I}(a_0; a_{1}, ..., a_{m}; a_{0}) = 
{\Bbb  I}(a_0; a_{1}, ..., a_{m}; a_{m+1}) + 
$$
$$
\sum_{1 \leq k \leq m-1}{\Bbb  I}(a_0; a_{1}, ..., a_{k}; a_{m+1}) \cdot
{\Bbb  I}(a_{m+1}; a_{k+1}, ..., a_{m}; a_0) + 
{\Bbb  I}(a_{m+1}; a_{1}, ..., a_{m}; a_0) 
$$
Applying the induction assumption to the second factors in the sum we get
$$
0 = {\Bbb  I}(a_0; a_{1}, ..., a_{m}; a_{m+1}) + 
$$
$$
\sum_{1 \leq k \leq m-1}(-1)^{m-k}{\Bbb  I}(a_0; a_{1}, ..., a_{k}; a_{m+1}) \cdot
{\Bbb  I}(a_{0}; a_{m}, ..., a_{k+1}; a_{m+1}) +
$$
$$ 
{\Bbb  I}(a_{m+1}; a_{1}, ..., a_{m}; a_0) 
$$
We claim that 
\begin{equation} \label{8.15.02.2}
{\Bbb  I}(a_0; a_{1}, ..., a_{m}; a_{m+1}) + 
\end{equation} 
$$
\sum_{1 \leq k \leq m-1}(-1)^{m-k}{\Bbb  I}(a_0; a_{1}, ..., a_{k}; a_{m+1}) \cdot
{\Bbb  I}(a_0; a_{m}, ..., a_{k+1}; a_{m+1}) 
$$
$$= (-1)^{m-1}
{\Bbb  I}(a_{0}; a_{m}, ..., a_{1}; a_{m+1})
$$
Indeed, this equality is equivalent to 
\begin{equation} \label{8.15.02.3}
{\Bbb  I}(a_0; a_{1}, ..., a_{m}; a_{m+1}) + 
\end{equation}
$$
\sum_{0 \leq k \leq m-1}(-1)^{m-k}{\Bbb  I}(a_0; a_{1}, ..., a_{k}; a_{m+1}) \cdot
{\Bbb  I}(a_0; a_{m}, ..., a_{k+1}; a_{m+1}) =0
$$
To prove (\ref{8.15.02.3}) we use the shuffle product formula 
to rewrite the sum. Then the claim is rather obvious 
(see the beginning of the proof of Theorem 4.1 in [G1] for details).  
It remains to notice that the proposition follows immediately from (\ref{8.15.02.2}).
The proposition is proved.

Let $\widetilde {\cal I}_{\bullet}(S)$ be 
defined in the same way as
${\cal I}_\bullet(S)$, except that
the relations 
 ii)-- iv) are omitted. By abuse of notation, we will denote the generators of 
$\widetilde {\cal I}_{\bullet}(S)$ by the same symbol as for 
${\cal I}_{\bullet}(S)$. 

\begin{proposition} \label{10.15.00.1}
The coproduct $\Delta$ provides ${\cal I}_{\bullet}(S)$, as well as  
$\widetilde {\cal I}_{\bullet}(S)$, with the structure of a 
commutative, graded  Hopf algebra. 
\end{proposition}

The algebra of regular functions ${\cal O}(G)$ on 
a group scheme $G$ is a commutative Hopf algebra with the coproduct 
$\Delta: {\cal O}(G) \lra {\cal O}(G)\otimes {\cal O}(G)$ induced by the group 
multiplication map $G \times G \lra G$. 
 
To prove Proposition  \ref{10.15.00.1}, which will follow from Theorem 2.5 below, 
 we interpret  
${\cal I}_{\bullet}(S)$ and  $\widetilde {\cal I}_{\bullet}(S)$  
as the algebras of regular functions 
on  certain pro--unipotent group schemes, defined as automorphism groups  
of certain non--commutative objects.

{\bf 2. The path algebra $P(S)$}. Let $S$ be a  set. Let $K$ be a field. 
Let $P(S)$ be the $K$-vector space with  basis 
\begin{equation} \label{7.19.02.1}
p_{s_0, ..., s_{n}}, \quad n \geq 1, s_k \in S \mbox{ for $k=0, ..., n$}
\end{equation}
It has a grading such that the degree of a typical basis element 
$p_{s_0, ..., s_{n}}$ is $-2(n-1)$. 
Let us equip it with the following structures.

i) {\it The $\circ$-product}. There is an associative product 
$$
\circ: P(S) \otimes_K P(S) \lra P(S)
$$
$$
p_{a, X , b}\circ p_{c, Y, d} := \left\{ \begin{array}{ll}
p_{a, X ,  Y, d}&: \quad b=c;
\\ 
0 &: \quad  b\not =c\end{array}\right.
$$
Here the small letters denote elements, and the capital ones subsequences, possibly empty, 
 of the set  $S$. 
In particular 
$p_{a,b} =  p_{a,x} \circ p_{x,b}$. 
The element $e_0:= \sum_{i \in S} p_{i,i}$ is the unit for this product.
The algebra $P(S)$ is decomposed into a sum
\begin{equation} \label{8.31.02.4}
P(S) = \oplus_{i,j\in S}P(S)_{i,j}
\end{equation}
where $P(S)_{i,j}$ is spanned by the elements (\ref{7.19.02.1}) with $s_0=i, s_n=j$. 
Below we consider only the automorphisms $F$ of  $P(S)$ 
preserving this decomposition:  
\begin{equation} \label{8.31.02.5}
F(P(S)_{i,j}) \subset P(S)_{i,j}
\end{equation}
The algebra $P(S)$ 
has an interpretation as a tensor 
algebra in a certain monoidal category, see section 2.3, which makes 
this restriction on $F$ 
natural.

ii) {\it The $\ast$-product}. We define {\it another} associative   product 
$$
\ast: P(S)(1) \otimes P(S)(1) \lra P(S)(1) 
$$ 
by the formula
\begin{equation} \label{8.31.02.3}
p_{X , b}\ast p_{c, Y} = \left\{ \begin{array}{ll}
p_{X , b , Y}&: \quad b=c;
\\ 
0 &: \quad  b\not =c\end{array}\right.
\end{equation}
where $b,c$ are elements, and $X$, $Y$ are ordered collections of elements of $S$. 
Here $M(1)$ means $M$ with the grading shifted down by $2$. 
Then $P(S)(1) $ is an associative algebra  
generated by the elements $p_{i,j}$.

One can add to this algebra  the elements $p_i$, $i \in S$, 
whose composition with the other elements is given by  formula  (\ref{8.31.02.3}) 
where $X$ or/and $Y$ can be empty. 
Then the $p_i$ are orthogonal projectors: $p_i^2=p_i$, 
$p_i \ast p_j =0$ if $i \not = j$, and $e:= \sum_{i\in S} p_{i}$ is  
the unit with respect to $*$.   
Using these projectors we can describe decomposition (\ref{8.31.02.4}) as 
$$
P(S) = \oplus_{i,j \in S} P(S)_{i,j};  \quad P(S)_{i,j}:= p_i \ast P(S) \ast p_j
$$
The automorphisms of the algebra $P(S)$ preserving the projectors 
$p_i$ respect this decomposition. 

Recall the graph  $\Gamma_S$ from s. 1.2. 
The $\ast$--algebra ${P}(S)$ is the path algebra of this graph. Namely, 
$p_{i_0, ..., i_{n}}$ corresponds to the path $i_0 \to i_1 \to ... \to i_{n}$ 
in  $\Gamma_S$. 
It turns out to be isomorphic to the free product, denoted $*_{\Q}$, and not to be confused with $*$,  of the 
algebra $\Q[S]$ of functions on $S$ with finite support with the 
polynomial algebra $\Q[x]$:
$$
{P}(S) \stackrel{\sim}{\lra} \Q[S]\ast_{\Q}\Q[x]
$$

iii) {\it The coproduct $\delta$}. Let us define a coproduct
$$
\delta: P(S) \lra P(S) \otimes_{K} P(S)
$$
such that
\begin{equation} \label{8.19.02.4}
\delta: P(S)_{a,b} \lms
P(S)_{a,b} \otimes P(S)_{a,b}
\end{equation} 
by the formula 
\begin{equation} \label{8.19.02.4ca}
\delta(p_{a, x_1, ..., x_n, b}) := \sum p_{a, x_{i_1}, ..., x_{i_k}, b} \otimes 
p_{a, x_{j_1}, ..., x_{j_{n-k}}, b}
\end{equation} 
where $n \geq 0$ and the sum is over all decompositions 
$$
\{1, ..., n\}  = \{i_1, ..., i_k\} \cup \{j_1, ..., j_{n-k}\}, \quad i_1 < ... < i_k; 
\quad j_1 < ... < j_{n-k}
$$
It is easy to see that it is cocommutative and coassociative.

{\it The compatibilities}. 
The expression $A \ast B \circ C \ast D \circ E$ does not depend on 
the bracketing. 

$\delta$ is an algebra morphism for the $\circ$--algebra structure:
\begin{equation} \label{8.19.02.2}
\delta (X \circ Y) = \delta (X)  \circ \delta (Y)
\end{equation}  

The compatibility with the $\ast$--structure is this:
\begin{equation} \label{8.19.02.3}
\delta (X \ast Y) = \delta (X)  \ast_2 \delta (Y); 
\end{equation} 
where  the operation $*_2$ is defined by 
$$
(X_1 \otimes Y_1) \ast_2 
(X_2 \otimes Y_2):= (X_1 \ast Y_1) \otimes 
(X_2 \circ Y_2) + (X_1\circ  Y_1) \otimes 
(X_2 \ast Y_2)
$$

It follows that $\delta$ is uniquely determined by (\ref{8.19.02.2}), 
(\ref{8.19.02.3}) and 
\begin{equation} \label{8.19.02.1}
\delta(p_{a,b}) = p_{a,b} \otimes p_{a,b}
\end{equation}

Observe that 
$$
\delta({P}(S)_{a,b}) = \delta(p_{a,a} \circ {P}(S) \circ p_{b,b}) \subset 
$$
$$
(p_{a,a} \otimes p_{a,a}) \circ({P}(S) \otimes {P}(S)) \circ 
(p_{b,b} \otimes p_{b,b}) = {P}(S)_{a,b} \otimes {P}(S)_{a,b}
$$
So $(\ref{8.19.02.2})$ plus $ (\ref{8.19.02.1})$ imply $ (\ref{8.19.02.4})$.

{\bf Remark}.  
$\delta (e_0) \not = e_0 \otimes e_0$.  
Apart from this, the coproduct $\delta$ 
has all the properties of a coproduct in a Hopf $\circ$--algebra.

For every
 $a \in S$ the $\circ$--algebra  ${P}(S)_{a,a}$ 
is a cocommutative  Hopf algebra with the unit $p_{a,a}$. 
It is the universal enveloping algebra of the Lie algebra  
${L}(S)_a$, 
reconstructed as the subspace of primitives in ${P}(S)_{a,a}$. 
The Lie algebra ${L}(S)_a$ is free 
 with generators labeled by the set $S$. 
The canonical Lie algebra isomorphism 
$i_{a,b}: {L}(S)_a \lra {L}(S)_b$  is given 
by 
$
l \lms p_{b,a} \circ l \circ p_{a, b}
$.

{\bf 3. The path algebra as a tensor algebra in a monoidal category}.  
Let $S$ be a  finite set. Let ${\cal C}$ be a tensor category. 
Consider the category $Q_{\cal C}(S)$ 
whose objects 
$$
V = \{V_{i, j}\}, \quad (i,j) \in S \times S
$$
are objects $V_{i,j}$ of ${\cal C}$
labeled by elements of  $S \times S$. 
Let $V = \{V_{i, j}\}$ and $W = \{W_{i, j}\}$. Then 
$$
{\rm Hom}_{Q_{\cal C}(S)}(V, W):= \oplus_{(i,j) \in S \times S}
{\rm Hom}_{{\cal C}}(V_{i, j}, W_{i,j})
$$
We define  $V \otimes W$ by setting
$$
(V \otimes W)_{i, j} = \oplus_{k \in S}V_{i, k} \otimes_{\cal C} W_{k, j}
$$

From now on we assume that ${\cal C}$ is the category of $K$--vector spaces,  
and denote by  $Q(S)$ the corresponding monoidal category. 
There is a unit object ${\Bbb I}$  defined as 
$$
{\Bbb I}_{i, j}:= \left\{ \begin{array}{ll}
0&: \quad i \not = j;
\\ 
K &: \quad  i = j\end{array}\right.   
$$

Consider an object ${\Bbb E}_S$ such that 
$$
 {\rm dim} ({{\Bbb E}_S})_{i,j} = 1 \quad \mbox{for any $(i,j) \in S \times S$} 
$$
Let ${T}({\Bbb E}_S)$ be the 
 free associative algebra with unit  
in the category $Q(S)$ 
generated by  ${\Bbb E}_S$.

\begin{lemma} \label{7.21.02.1}
${T}({\Bbb E}_S)$ is isomorphic to the path $\ast$--algebra 
${P}(S)(1)$. 
\end{lemma} 

{\bf Proof}. For any $i,j \in S$ choose a nonzero vector  $p_{i,j}$ in $({\Bbb E}_S)_{i,j}$. Set 
\begin{equation} \label{7.21.02.2}
p_{i_0, ..., i_{n}} := p_{i_0, i_1} \otimes p_{i_1, i_2} \otimes ... \otimes 
p_{i_{n-1},  i_{n}}
\end{equation} 
Then 
$$
{T}({\Bbb E}_S) = \oplus_{n \geq 0}{\Bbb E}_S^{\otimes n} = {\Bbb I} \oplus
\oplus_{n \geq 1} \Bigl(\oplus_{i_0, ..., i_{n} \in S} K 
\cdot p_{i_0, ..., i_{n}} \Bigr)
$$
Writing ${\Bbb I}_{i,i} = K \cdot p_i$ we get the needed isomorphism. 
The lemma is proved.

Below we identify $p_{i,j}$ with the object of $Q(S)$ whose 
$(i',j')$--component is zero unless $i=i', j=j'$, when it is $K$. 
Then  we have the analogous objects to (\ref{7.21.02.2}) in the category 
$Q(S)$.

There is the second product
$$
\circ: {T}({\Bbb E}_S) \otimes {T}({\Bbb E}_S) \lra {T}({\Bbb E}_S); 
\quad {\rm deg}(\circ) = -2
$$

{\bf 4. The pro--unipotent algebraic group scheme $G(S)$}. 
\begin{definition} \label{7.18.02.1} a) ${\rm Aut}_0({P}(S))$ is 
the pro--unipotent group scheme of automorphisms 
of the $\ast$--algebra ${P}(S)$ which act as the identity on the quotient 
\begin{equation} \label{7.18.02.6}
{P}_+(S)/({ P}_+(S))^2
\end{equation}

b) $G(S) \subset {\rm Aut}_0({P}(S))$ is the subgroup  of 
all automorphisms $F$ such that: 

i) $F$ commutes with $\delta$. 

ii) $F$ is an automorphism of the $\circ$--algebra structure.  
\end{definition}

The grading of ${P}(S)$ provides natural gradings  
of the algebras of regular functions on 
$G(S)$ and ${\rm Aut}_0({P}(S))$.

\begin{theorem} \label{7.18.02.2} a) The  commutative, graded algebra 
${\cal O}\Bigl({\rm Aut}_0({P}(S))\Bigr)$ is 
identified with the polynomial 
algebra in an infinite number of variables 
$I_{s_0, ..., s_{n+1}}$ with $n \geq 0$ and $s_i \in S$:
$$
{\cal O}\Bigl({\rm Aut}_0({P}(S))\Bigr) = K[I_{s_0, ..., s_{n+1}}]; 
\quad \mbox{${\rm deg}(I_{s_0, ..., s_{n+1}}) = -2n$}
$$

b) The identification in a) provides an isomorphism of 
 commutative, graded Hopf algebras
\begin{equation} \label{7.18.02.4}
{\cal O}\Bigl({\rm Aut}_0({ P}(S))\Bigr) \stackrel{\sim}{\lra} \widetilde 
{\cal I}_{\bullet}(S)
\end{equation}

c) The isomorphism (\ref{7.18.02.4}) induces an isomorphism of 
 commutative, graded Hopf algebras
\begin{equation} \label{7.18.02.5}
{\cal O}(G(S)) \stackrel{\sim}{\lra} {\cal I}_{\bullet}(S)
\end{equation}
\end{theorem}

{\bf Proof}. We consider first the case when $S$ is finite, and then 
take the inductive limit over finite subsets of $S$. 

An automorphism $F$ of $P(S)$ satisfying (\ref{8.31.02.5}), 
being a $\ast$--algebra automorphism, 
is uniquely determined by its values on the generators $p_{a,b}$. 
Indeed, condition (\ref{8.31.02.5}) plus Lemma \ref{7.21.02.1} imply 
 that it is 
an automorphism of the free $\ast$--algebra in 
the monoidal category $Q(S)$ generated by 
the elements $p_{a,b}$.  
 
Let us write 
$$
F(p_{a,b}) = p_{a, b} + \sum I_{a, s_1, ..., s_m, b}(F) \cdot p_{a, s_1, ..., s_m, b}
$$
where the summation is over all nonempty ordered collections of elements  
$s_1, ..., s_m$ of $S$. Here  $I_{a, s_1, ..., s_m, b}(F) \in K$ are 
the coefficients,  providing a regular function 
 $$
I_{a, s_1, ..., s_m, b} \in {\cal O}({\rm Aut}_0({P}(S)))
$$    
Observe that $I_{a,  b}(F) =1$ just means that
$F$ acts as the identity on (\ref{7.18.02.6}). 

We claim that the map 
$$
I_{a, s_1, ..., s_m, b} \lra {\Bbb I}(a; s_1, ..., s_m; b)
$$
provides the  isomorphism (\ref{7.18.02.4}). Indeed, it is 
obviously an isomorphism of commutative, graded $K$--algebras. So we get a).

The crucial fact that the coproduct in 
${\cal O}\Bigl({\rm Aut}_0({P}(S))\Bigr)$ is identified with $\Delta$ as 
given in 
(\ref{CP2**}) 
follows easily from the very definitions.  So we have b).

c) The condition that $F$ commutes with the $\circ$--product 
 is equivalent to  iii) plus iv) in the definition of the Hopf algebra 
${\cal I}_{\bullet}(S)$. 
Indeed, the path composition formula iii) is equivalent to $F(p_{a,b})  = F(p_{a,x}) \circ F(p_{x,b})$. 
The $\circ$--unit is given by $\sum_{i\in S}  p_{i,i}$. 
The condition iv)  
just means that $F$ preserves the $\circ$--unit: 
$F(\sum  p_{i,i}) =\sum  p_{i,i}$.  

Given that $F$ preserves both products, 
the fact that $F$ commutes with $\delta$  is 
equivalent to the shuffle product formula ii). 
Indeed,  
the condition 
$$
\delta F (p_{a,b})  = F \delta (p_{a,b} ) \stackrel{(\ref{8.19.02.1})}{=} 
F (p_{a,b} \otimes p_{a,b} ) = F (p_{a,b}) \otimes F (p_{a,b})
$$
is just equivalent to the shuffle product formula. 
Since  $\delta$ is completely determined by (\ref{8.19.02.1}) 
and compatibilities with $\ast$ and $\circ$, the statement follows.

The part c) and hence the theorem are proved. 

Proposition \ref{10.15.00.1} follows from Theorem \ref{7.18.02.2}.

\section{The motivic fundamental groupoid and its Galois group}

In this chapter we show that the basic properties of the motivic 
fundamental groupoid ${\cal P}^{\cal C}({\Bbb A}^1_S;S)$ 
in one of the mixed Tate categories ${\cal C}$ described in s. 3.2 
imply that the canonical fiber functor on ${\cal C}$ sends it to 
the path algebra ${P}(S)$. 
This immediately implies Theorem \ref{8.3.02.1}. 

{\bf 1. The Betti realization of the fundamental groupoid}. 
Let $S$ be a subset of $\C$. Let $t$ be a standard coordinate on $\C$. 
Choose a tangent vector $v_s$ at every  
$s \in S$ such that $dt(v_s)=1$. 

For $a \in S$, let $\C_S:= \C - S$.
Let us recall the pronilpotent completion 
of the topological  torsors of paths associated with $\C_S$. 
Let $I_a$ be the augmentation 
ideal of the group ring $\Z[\pi_1(\C_S; a)]$ of the topological 
fundamental group 
$\pi_1(\C_S; a)$. Denote by ${\cal P}(\C_S; a,b)$, where $a,b \in S$, 
 the free 
abelian group generated by the set of homotopy classes 
of paths between the tangential base points  $v_a, v_b$. Then
$$
\pi_1^{\rm nil}(\C_S; a):= 
\lim_{\longleftarrow}\Z[\pi_1(\C_S; a)]/I_a^n; \quad 
{\cal P}^{\rm nil}(\C_S; a,b):= 
\lim_{\longleftarrow} I_a^n \backslash {\cal P}(\C_S; a,b)
$$

There are the   path composition morphisms: for any $a,b,c \in S$
$$
\circ: {\cal P}^{\rm nil}(\C_S; a,b) \otimes  
{\cal P}^{\rm nil}(\C_S; b, c) \lra 
{\cal P}^{\rm nil}(\C_S; a,c)
$$
They provide ${\cal P}^{\rm nil}(\C_S; a,b)$ 
with the structure of a principal homogeneous space over 
 $\pi^{\rm nil}_1(\C_S; a)$.  There is a coproduct map
$$
\delta: {\cal P}^{\rm nil}(\C_S; a,b) \lra {\cal P}^{\rm nil}(\C_S; a,b) 
\otimes {\cal P}^{\rm nil}(\C_S; a,b)
$$
provided by the map given on the generators by 
$\gamma \lms \gamma \otimes \gamma$. 
It is obviously compatible with the composition of paths. 
It makes $\pi_1^{\rm nil}(\C_S; a)$ into a cocommutative Hopf algebra. 

We have an increasing  filtration 
$W_{\bullet}$ of ${\cal P}^{\rm nil}(\C_S; a,b)$ indexed by $0, -2, -4, ...$:
$$
W_{-2n}{\cal P}^{\rm nil}(\C_S; a,b):= I_a^n \circ {\cal P}^{\rm nil}(\C_S; a,b)
$$ 
There are canonical isomorphisms
$$
{\rm gr}^W_{0}{\cal P}^{\rm nil}(\C_S; a,b) = \Z(0); \quad 
{\rm gr}^W_{-2}{\cal P}^{\rm nil}(\C_S; a,b) = H_1(\C_S, \Z) = \Z[S]
$$
Denote by  $p_{a,b}$ the canonical generator of 
${\rm gr}^W_0{\cal P}^{\rm nil}(\C_S; a,b)$. 
The second isomorphism is provided by the map
$$
s\in S \lms p_{a,s, b}:= p_{a,s} \circ ([\gamma_s] - [1]) \circ 
p_{s, b} \in {\rm gr}^W_{-2}{\cal P}^{\rm nil}(\C_S; a,b)
$$
 where $\gamma_s$ is 
a simple loop around $s$ based at $v_s$, and 
$1$ is the identity loop at $v_s$, see the left picture below.

Set
$$
{\cal P}^{\rm nil}(\C_S; S) := \oplus_{a,b \in S}{\cal P}^{\rm nil}(\C_S; a,b)
$$
The $\ast$--algebra structure on ${\cal P}^{\rm nil}(\C_S; S)(1)$ is given by 
$$
\alpha_{a,b} \ast \alpha_{b, c}:= 
\alpha_{a,b} \circ ([\gamma_b] -[1]) \circ \alpha_{b, c}
$$

\begin{lemma} \label{8.29.02.2} There is a canonical isomorphism
$$
 {P}(S)\stackrel{\sim}{\lra}{\rm gr}^W_{\bullet}{\cal P}^{\rm nil}(\C_S; S)
$$
 It respects the grading, the $\ast$--algebra and $\circ$--algebra structures 
and the coproduct $\delta$.  
\end{lemma}

{\bf Proof}. The isomorphism is given by 
$$
p_{a, s_1, ..., s_n, b} \lms p_{a,s_1} 
\circ ([\gamma_{s_1}] - [1]) \circ p_{s_1, s_2} 
\circ ([\gamma_{s_2}] - [1]) \circ ... \circ p_{s_n, b} 
$$
The picture below illustrates this formula in the case when $n=1$ and $n=2$. 
\begin{center}
\hspace{4.0cm}
\epsffile{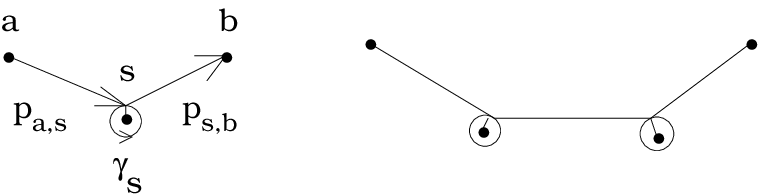}
\end{center}

The proof follows easily from the fact that, if $S$ is finite,  $\pi_1(\C_S; a)$ 
is a free group with $|S|$ generators, and hence 
$\pi^{\rm nil}_1(\C_S; a)$ is the completion of the tensor algebra 
generated by the elements $p_{a,s, a}$, $s\in S$. 
The lemma is proved. 

Now let us discuss the motivic fundamental groupoids.

{\bf 2. The set--up}. Let $F$ be a field. In the rest of this Chapter  we  
work in one of the following categories  ${\cal C}$. Their properties 
and the  corresponding 
formalism are presented in the Appendix.

i)  the abelian category of mixed Tate 
motives over a number field $F$ ([L1], [G5], [DG]).

ii) $F=\C$, and  ${\cal C}$ is the category of Hodge-Tate structures. 
(Appendix, s. 3.2).

iii) $F$ is an arbitrary field such that $\mu_{l^{\infty}} \not \in F$, and 
 ${\cal C}$ 
is 
 the mixed Tate category of 
$l$-adic Tate ${\rm Gal}(\overline F/F)$--modules. 
(Appendix, s. 3.3,  and the references there).

There is also one hypothetical set--up: 

iv) $F$ is an arbitrary field, ${\cal C}$ is the 
  hypothetical abelian category of mixed Tate 
motives over $F$. 

Any   category ${\cal C}$ from the list above  is a mixed Tate $K$--category, 
where $K = \Q$ in i), ii) and $K=\Q_l$ in iii), see [BD] or the 
Appendix, subsections 1 and 2, for the background. 
Let us briefly recall the main features which will be used. 
${\cal C}$ is generated as a tensor category by a simple object $K(1)$.  
 Each object carries a canonical weight filtration 
$W_{\bullet}$, and morphisms in ${\cal C}$ are strictly 
compatible with this filtration.  
There is a canonical fiber functor 
to the category of finite dimensional graded vector spaces:
$$
\omega: {\cal C} 
\lra {\rm Vect}_{\bullet}, \quad X \lms \oplus_n Hom_{\cal C}(K(-n), 
{\rm gr}^W_{2n}X)
$$ 
Forgetting the grading we get a fiber functor $\widetilde \omega$. The space ${\rm End}(\widetilde \omega)$
of its endomorphisms is a  graded Hopf algebra. 
Its graded dual is canonically isomorphic to the commutative graded Hopf algebra 
${\cal A}_{\bullet}(\cal C)$ 
of the framed objects in ${\cal C}$
$$
{\cal A}_{\bullet}(\cal C) \stackrel{\sim}{=} \quad \mbox{graded dual 
${\rm End}(\widetilde \omega)^{\vee}$},
$$ 
see Theorem \ref{8-23.3/99} in the Appendix. 
The functor $\omega$ provides a 
canonical equivalence between the category ${\cal C}$ and the 
category of finite dimensional graded ${\cal A}_{\bullet}(\cal C)$--comodules. 
${\rm Spec}({\cal A}_{\bullet}({\cal C}))$ is a 
 pro--unipotent group scheme. The grading on  ${\cal A}_{\bullet}(\cal C)$ 
encodes a natural semidirect product of ${\Bbb G}_m$ and 
${\rm Spec}({\cal A}_{\bullet}({\cal C}))$.

{\bf 3. The fundamental groupoid}. 
Let $S$ be any subset of $F = {\Bbb A}^1(F)$. A particular interesting 
case is  $S =F$. Choose a tangent vector $v_s$ at every point 
$s \in S$. We assume $v_s$ is defined over  $F$. The  
 differential $dt$  provides  a canonical 
choice of $v_s$ for all $s \in F$ which is used below. 

Let 
${\cal P}^{\cal C}({\Bbb A}^1_S; S)$ be the fundamental groupoid 
of paths on ${\Bbb A}^1_S:= {\Bbb A}^1 - S$ between the tangential base 
points $v_s$. 
It is a pro--object in  ${\cal C}$.

The  fundamental groupoid 
in the situations ii) and iii) above was defined by Deligne [D]. In the situation i) 
it is defined in [DG] (another 
 construction of the same object is given in  Chapter 6  of [G8]).
In particular there are  path composition morphisms in ${\cal C}$
\begin{equation} \label{7.25.02.20}
\circ: {\cal P}^{\cal C}({\Bbb A}^1_S; a,b) \otimes 
{\cal P}^{\cal C}({\Bbb A}^1_S; b, c) \lra 
{\cal P}^{\cal C}({\Bbb A}^1_S; a,c)
\end{equation}
They provide ${\cal P}^{\cal C}({\Bbb A}^1_S; a,b)$ 
with a structure of principal homogeneous space over 
the fundamental group ${\cal P}^{\cal C}({\Bbb A}^1_S; a, a)$, 
understood as a Hopf algebra in ${\cal C}$.

{\bf 4. The  structures on the  fundamental groupoid}. 
Let us assume that $S$ is finite. 
Recall the monoidal category $Q_{\cal C}(S)$, see s. 2.3. 
We define an object ${\cal P}^{\cal C}(S)$ in  $Q_{\cal C}(S)$ by 
$$
{\cal P}^{\cal C}(S)_{a,b}:= {\cal P}^{\cal C}({\Bbb A}^1_S; a,b)
$$
There are the following structures on this object ([DG]).

i) The path composition morphisms (\ref{7.25.02.20}) provide ${\cal P}^{\cal C}(S)$ 
with  the structure of an algebra in $Q_{\cal C}(S)$, called 
the  $\circ$--algebra structure.

ii) There are canonical ``loop around $s$'' morphisms
$$
\gamma_s: \Q(1) \lra {\cal P}^{\cal C}({\Bbb A}^1_S; s,s)
$$
They provide   morphisms
$$
\ast: {\cal P}^{\cal C}({\Bbb A}^1_S; a,b)(1) \otimes {\cal P}^{\cal C}({\Bbb A}^1_S; b, c)(1) \lra 
{\cal P}^{\cal C}({\Bbb A}^1_S; a,c)(1)
$$ 
$$
 \alpha_{a,b}\ast  \alpha_{b,c}:= \alpha_{a,b} \circ \gamma_b \circ \alpha_{b, c}
$$
These morphisms make ${\cal P}^{\cal C}(S)(1)$ into an  algebra in 
the category $Q_{\cal C}(S)$. We call it the $\ast$--algebra structure.

iii) For any $a,b \in S$ there is a coproduct  given by a 
morphism in  ${\cal C}$ 
$$
\delta: {\cal P}^{\cal C}({\Bbb A}^1_S; a,b) \lra 
{\cal P}^{\cal C}({\Bbb A}^1_S; a,b)\otimes {\cal P}^{\cal C}({\Bbb A}^1_S; a,b)
$$
It is
 a $\circ$--algebra morphism. 
One has 
\begin{equation} \label{8.20.02.1}
\delta(\gamma_b) = \gamma_b \otimes 1 + 1 \otimes \gamma_b; 
\end{equation} 
  It follows from this that 
the compatibility of $\delta$ with the $\ast$--algebra 
product is given by (\ref{8.19.02.3}).

{\bf 5. The  path algebra provided by the fundamental groupoid and 
$G_{\cal C}(S)$}. 
Let us apply the fiber functor $\omega$ to ${\cal P}^{\cal C}(S)$, getting  
$$
{\cal P}_{\omega}(S):= \omega({\cal P}^{\cal C}(S))
$$

{\bf Remark}. Since ${\cal P}^{\cal C}(S)$ is a pro--object in ${\cal C}$, 
${\cal P}_{\omega}(S)$ is a 
projective limit of  finite dimensional $K$--vector spaces. However 
its graded components are finite dimensional. 
Denote by ${\cal P}^{\rm gr}_{\omega}(S)$ the direct sum of the 
graded components of 
${\cal P}_{\omega}(S)$.

Then ${\cal P}^{\rm gr}_{\omega}(S)$ is an algebra in the monoidal category 
$Q_{{\cal A}_{\bullet}({\cal C})-mod}(S)$, and  
${\cal P}^{\rm gr}_{\omega}(S)(1)$ is a $\ast$--algebra in the same category. 

The space 
$$
{\rm Hom}(K(0), {\rm gr}^W_0{\cal P}^{\cal C}({\Bbb A}^1_S; a,b)) 
$$
is one--dimensional. It has a natural generator  $p_{a,b}$ such that 
\begin{equation} \label{8.20.02.2}
\delta(p_{a,b}) = p_{a,b}\otimes p_{a,b}
\end{equation}
Set 
\begin{equation} \label{8.31.02.1}
p_{a; s_1, ..., s_m; b}:= p_{a,s_1} \ast p_{s_1, s_2} \ast  ... \ast 
p_{s_{m-1}, s_m} \ast
p_{s_m, b} 
\end{equation}
It follows from (\ref{8.20.02.1}) and (\ref{8.20.02.2}) 
that $\delta$ is given by formula (\ref{8.19.02.4ca}).

\begin{proposition} \label{7.24.02.11} There is a natural isomorphism 
of $K$-vector spaces 
$$
{\cal P}^{\rm gr}_{\omega}(S) \lra {P}(S)
$$
respecting the grading, the $\ast$--algebra and $\circ$--algebra structures, 
and the coproduct on both objects. 
\end{proposition}

{\bf Proof}. 
The statement boils down to the fact that the elements (\ref{8.31.02.1}), 
when the set $\{s_1, ..., s_m\}$ runs through all elements of $S^m$, 
form a basis in 
$
{\rm gr}^W_{2m}{\cal P}^{\cal C}({\Bbb A}^1_S; a, b)
$. 
Since  ${\cal P}^{\rm nil}(\C_S; S)$ 
is the Betti realization of the Hodge version ${\cal P}^{\cal H}({\Bbb A}^1_S; S)$,  
there is an  isomorphism  
$$
{\cal P}_{\omega}(S)_{(n)} \stackrel{\sim}{=} 
{\rm gr}^W_{-2n}{\cal P}^{\rm nil}(\C_{S}; S)
$$
This plus  
Lemma \ref{8.29.02.2} implies the statement in the  Hodge setting. 
Hence (by the injectivity of the regulators, see, for example, [DG]) 
it is true in the motivic situation i). The $l$-adic case is handled similarly 
using the comparison theorem with the Betti realization. 
The proposition is proved. 

We define an element 
\begin{equation} \label{7.26.02.1}
{\rm I}^{\cal C}(a_0; a_1, a_2, ..., a_m; a_{m+1}) \in {\cal A}_m({\cal C})
\end{equation}
as the linear functional on ${\rm End}(\omega)$ 
given by the matrix element
$$
F\in {\rm End}(\omega) \lms <F(p_{a,b}), p_{a; s_1, ..., s_m; b}>
$$
One can describe the elements of ${\cal A}_{\bullet}({\cal C})$ 
by framed objects in ${\cal C}$, see Subsection 2 of the Appendix. 
Then element (\ref{7.26.02.1}) is represented by 
the following framed object in ${\cal C}$:
$$
\Bigl({\cal P}^{\cal C}({\Bbb A}^1_S; a, b); p_{a,b}, 
p^*_{a; s_1, ..., s_m; b}\Bigr)
$$ 
Here $p^*_{...}$ are the elements of the basis dual to $p_{...}$.

Recall the group $G(S)$ defined in Chapter 2. 
The group scheme ${\rm Spec}({\cal A}_{\bullet}({\cal C}))$ 
acts on ${\cal P}_{\omega}(S)$ through its quotient, denoted 
$G_{\cal C}(S)$. The semidirect product of ${\Bbb G}_m$ and 
 $G_{\cal C}(S)$ 
 is called the Galois 
group of the fundamental groupoid ${\cal P}^{\cal C}({\Bbb A}^1_S; S)$.

\begin{theorem} \label{7.24.02.12} 
a)  $G(S)$ is the group of all automorphisms of 
the $\ast$--algebra ${\cal P}^{\rm gr}_{\omega}(S)(1)$ in the monoidal category 
$Q_{{\cal A}_{\bullet}({\cal C})-mod}(S)$  
preserving the $\circ$--algebra structure on ${\cal P}^{\rm gr}_{\omega}(S)$, 
commuting with the coproduct $\delta$   
 and acting as the identity  on 
\begin{equation} \label{8.31.02.2}
{\cal P}_{\omega}(S)(1)/({\cal P}_{\omega}(S)(1))^2
\end{equation}

b) 
The map 
\begin{equation} \label{7.25.02.1}
{\Bbb I}(a_0; a_1, a_2, ..., a_m; a_{m+1}) \lra 
{\rm I}^{\cal C}(a_0; a_1, a_2, ..., a_m; a_{m+1})
\end{equation}
provides a  surjective morphism of Hopf algebras 
${\cal O}(G(S)) \twoheadrightarrow 
{\cal O}(G_{\cal C}(S))$, and hence an inclusion of 
pro--unipotent group schemes 
$$
G_{\cal C}(S) \hookrightarrow G(S)
$$  

c) The coproduct is computed by the formula 
\begin{equation} \label{CP2*}
\Delta {\rm I}^{\cal C}(a_0; a_1, a_2, ..., a_m; a_{m+1})= 
\end{equation}
$$
\sum_{0 = i_0 < i_1 < ... < i_k < i_{k+1} = m+1} 
{\rm  I}^{\cal C}(a_0; a_{i_1}, ..., a_{i_k}; a_{m+1}) \otimes \prod_{p =0}^k
{\rm  I}^{\cal C}(a_{i_{p}}; a_{i_{p}+1}, ..., a_{i_{p+1}-1}; a_{i_{p+1}})
$$
\end{theorem}

{\bf Proof}. Part c) means simply that (\ref{7.25.02.1}) 
commutes with the coproduct. Part b) follows immediately from a). 
And a) is a direct consequence of Theorem \ref{7.18.02.2} and 
Proposition \ref{7.24.02.11}. 
The theorem is proved.

{\bf Example}. When $m=3$ the formula (\ref{CP2*}) gives  
$$
\Delta {\rm I}^{\cal C}(a_0; a_1, a_2, a_3; a_4) \quad = \quad 1 
\otimes {\rm I}^{\cal C}(a_0; a_1, a_2, a_3; a_4) +
$$
$$
{\rm I}^{\cal C}(a_0; a_1; a_4) \otimes {\rm I}^{\cal C}(a_1; a_2, a_3; a_4) +  
$$
$$
{\rm I}^{\cal C}(a_0; a_2; a_4) \otimes {\rm I}^{\cal C}(a_0; a_1; a_2) \cdot {\rm I}^{\cal C}(a_2; a_3; a_4) +
$$
$$
{\rm I}^{\cal C}(a_0; a_3; a_4) \otimes {\rm I}^{\cal C}(a_0; a_1, a_2; a_3) +  
{\rm I}^{\cal C}(a_0; a_1, a_2; a_4) \otimes {\rm I}^{\cal C}(a_2; a_3; a_4) +
$$
$$
{\rm I}^{\cal C}(a_0; a_1, a_3; a_4) \otimes {\rm I}^{\cal C}(a_1; a_2; a_3) +
{\rm I}^{\cal C}(a_0; a_2, a_3; a_4) \otimes {\rm I}^{\cal C}(a_0; a_1; a_2) +
$$
$$
{\rm I}^{\cal C}(a_0; a_1, a_2, a_3; a_4) \otimes 1 
$$

There is an equivalent version of Theorem \ref{7.24.02.12}a):

{\it $G(S)$ is the group of all automorphisms of 
the $\circ$--algebra ${\cal P}^{\rm gr}_{\omega}(S)$ in the monoidal 
category $Q_{{\cal A}_{\bullet}({\cal C})-mod}(S)$ preserving the 
elements $\gamma_s$, the coproduct $\delta$,  
and acting as the identity on (\ref{8.31.02.2})}. 

Indeed, the $\ast$--product is determined by the 
$\circ$--product and the elements $\gamma_s$.

\section{Iterated integrals and plane trivalent rooted trees}

In this chapter we show how the Hopf algebra 
of decorated rooted plane trivalent trees 
encodes the properties of the motivic iterated integrals.   

{\bf 1. Terminology}. A tree can have two types of edges: 
internal edges and legs. Both vertices of an internal edge 
are of valency $\geq 3$. 
A valency $1$ vertex of a leg is called an external vertex.

%A {\it plane  tree} is a tree 
%located on the plane. 
We may picture  plane trees inscribed into a circle such that 
the  external vertices divide this circle into a union of arcs. 
Let $S$ be a set. An {\it $S$--decoration} of a plane  tree is 
an $S$--valued function on the set of open arcs.  

A {\it rooted tree} is a tree with one distinguished leg called the root. 
We picture plane rooted trees growing down from the root,
and put the external vertices of all the legs but the root  on a line. These 
vertices  divide this line into a union of arcs. 
 They correspond to the arcs defined above, so we can talk about 
$S$--decorated plane rooted trees.

\begin{center}
\hspace{4.0cm}
\epsffile{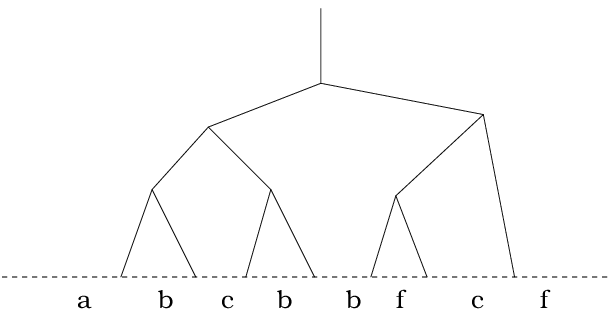}
\end{center}

We say that a plane rooted tree is {\it decorated by an ordered set 
$\{s_0, s_1, ..., s_m, s_{m+1}\}$} if the tree has $m+1$ bottom legs, 
and the corresponding $m+2$ arcs are decorated, in their natural order from 
left to  right,  by the set above. 
We  decorate the arcs, not legs -- see 
however the 
remark after Theorem 4.2.

 {\bf 2. The 
Hopf algebra ${\cal T}_{\bullet}(S)$ of $S$--decorated plane rooted 
trivalent trees}. 
It is a  commutative, graded Hopf algebra. 
As a vector space it is generated by disjoint unions of $S$--decorated plane 
rooted trivalent trees. 
More precisely, consider a graded vector space with  a basis given by $S$--decorated
plane rooted trivalent trees with $n+2$ arcs, where $n \geq 1$ provides the grading. The vector space  ${\cal T}_{\bullet}(S)$ 
is its symmetric algebra. 
In other words, ${\cal T}_{\bullet}(S)$ is the commutative algebra generated by 
$S$--decorated plane rooted trivalent trees  with the only relation 
that a decorated plane tree with just one edge equals  $1$ 
(regardless of the decoration). 

Let us define the coproduct $\Delta_{\cal T}$ on ${\cal T}_{\bullet}(S)$. 
Since $\Delta_{\cal T}$ has to be an algebra morphism, we 
have  to define it only on the generators. 
Let $T$ be an $S$--decorated 
plane rooted trivalent tree. Consider the  set
$$
{\cal E}_{T}:= \quad \mbox{$\{$ all internal edges of $T$ $\}$  $\cup$  $\{$ 
the root of $T$ $\}$}
$$
An element $E \in {\cal E}_{T}$
 determines a plane trivalent rooted tree $T_{E}$ 
growing down from $E$: the edge $E$ serves as the root of this tree. 
The tree $T_{E}$ inherits a natural $S$--decoration:  take the arcs 
containing the endpoints of the tree $T_{E}$ and keep their decorations.  
A subset  $\{E_1, ..., E_k\}$  of ${\cal E}_{T}$  
is called {\it admissible} if for any  $i \not = j$ the edges $E_i$ and $E_j$ do not lie on the same path going down from the root.  In other words, $E_i$ is not contained 
in the tree $T_{E_j}$ if $i \not = j$. Such an admissible subset determines 
a connected $S$--decorated plane trivalent tree 
$T/(T_{E_1} \cup ... \cup  T_{E_k})$. Namely, this tree is obtained by 
shrinking each of the rooted trees $T_{E_1}, ...,   T_{E_k}$ into a 
leg of a new tree: we shrink the domain encompassed by each of these trees and 
the bottom line. In particular we shrink to points all the arcs 
on the bottom line located under these trees. The remaining arcs 
with the inherited decorations are the arcs of the new tree. 
See the picture where we shrank the two trees under the thick  edges:

\begin{center}
\hspace{4.0cm}
\epsffile{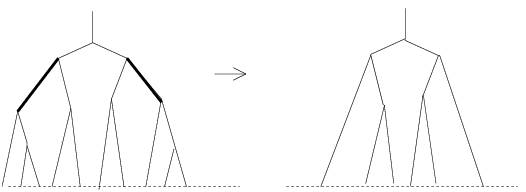}
\end{center}
Now we set, for a $T$ as above, 
$$
\Delta_{\cal T}(T) := \sum \frac{T}{T_{E_1} \cup ... \cup  T_{E_k}} \otimes \prod_{i=1}^k 
T_{E_i}
$$
where the sum is over all admissible
 subsets of ${\cal E}(T)$, including  the empty subset
and the  subset formed by the root. 

\begin{lemma} \label{7.23.02.1}
$\Delta_{\cal T}$ provides ${\cal T}_{\bullet}(S)$ with the structure of a commutative, 
graded 
Hopf algebra.
\end{lemma}

{\bf Proof}. Left to the reader as an easy exercise. 
It is essentially the same as in  [CK]. 

{\bf Remark}. The Hopf algebra ${\cal T}_{\bullet}(S)$ is 
a Hopf subalgebra of the similar Hopf algebra 
of all  decorated rooted plane trees (not necessarily trivalent). 
The latter is the plane decorated version of 
the one defined by Connes and Kreimer [CK]. 
On the other hand, as was pointed out to me by J.-L. Loday, 
the Hopf algebra ${\cal T}_{\bullet}(S)$ 
is closely related to the structures studied in [Lo1], [Lo2].

\begin{theorem} \label{7.23.02.2}
 The map
$$
t: {\Bbb I}(s_0; s_1, ..., s_m; s_{m+1})  \longmapsto  \mbox{\rm sum of all 
plane rooted trivalent trees}
$$
$$
\mbox{\rm decorated by the ordered set 
$\{s_0, s_1, ..., s_m, s_{m+1}\}$}
$$
provides an injective homomorphism of commutative, graded
Hopf algebras 
\begin{equation} \label{7.24.02.1}
t: \widetilde {\cal I}_{\bullet}(S) \hookrightarrow {\cal T}_{\bullet}(S)
\end{equation}
\end{theorem}

{\bf Proof}. We need to check only that our map commutes with the coproducts. 
Let $T$ be a generator of ${\cal T}_{\bullet}(S)$, decorated by an ordered set 
$\{s_0, ..., s_{m+1}\}$. 

Observe that an edge $E$ determines 
a subset $S_E \subset \{s_0, ..., s_{m+1}\}$ consisting of the labels 
of all arcs located 
between the very left and right bottom legs 
of the tree $T_{E}$, i.e. just under this tree. 
Different edges produce different subsets. 

Let 
   $\{E_1, ..., E_k\}$ be an admissible subset of ${\cal E}_T$. 
The subsets $S_i$ corresponding 
to the edges $E_i$ are disjoint -- this is a reformulation of the 
definition of an admissible  subset. In fact  they 
are even separated: there is at least one arc between the arcs labeled 
by $S_i$ and $S_{i+1}$. 

We claim that the sum over all 
trees decorated by $\{s_0, ..., s_{m+1}\}$, with such subsets 
$S_1, .., S_k$ given to us, corresponds to a single term of the formula 
$\Delta {\Bbb I}(s_0; s_1, ..., s_m; s_{m+1})$. Namely, if we 
enumerate  the terms in the formula  using  
 (\ref{7.23.2.6}),  
this term corresponds to the subset 
$\{s_0, ..., s_{m+1}\} - (S_1 \cup ... \cup S_k)$.  

\begin{center}
\hspace{4.0cm}
\epsffile{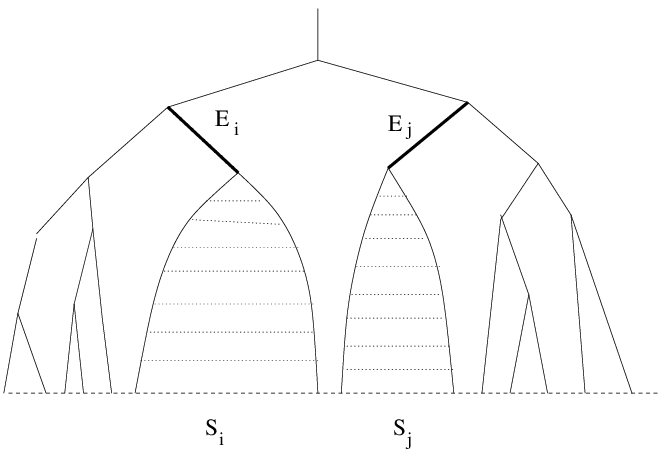}
\end{center}

To verify this claim observe two things: 

i) given such a tree $T$ we can alter the 
subtree $T_{E_i}$ to any other subtree with the same root $E_i$ and  the same 
set of arcs ``under'' this tree,  labeled by $S_i$. 

ii) if we shrink, for any admissible subset $(E_1, ..., E_k)$ of ${\cal E}_T$, each of the trees $T_{E_1}, ...,  T_{E_k}$ into a leg, then  the result can be any tree 
labeled by the ordered set $\{s_0, ..., s_{m+1}\} - (S_1 \cup ... \cup S_k)$. 

The theorem is proved.

{\bf Remark}. 
The element ${\rm I}^{\cal C}(s_0; s_1, ..., s_m; s_{m+1})$ 
is invariant under the action of the  translation group of the affine line.
Therefore it is natural to encode it by orbits of the action 
$$
\{s_0, s_1, ..., s_m, s_{m+1}\} \lra \{s_0+a, s_1+a, ..., s_m+a, s_{m+1}+a\}
$$
of the additive group 
${\Bbb A}^1$ on the set of $(m+2)$--tuples $
\{s_0, s_1, ..., s_m, s_{m+1}\}$.  
These orbits are described by 
the $(m+2)$--tuples 
$$
\{s_0 - s_{m+1}, s_1-s_0, s_2-s_1, ..., s_{m+1} - s_m\}
 $$
which naturally sit not at the arcs, but rather 
at the legs of the corresponding trees.

{\bf 3. The $\otimes^m$--invariant of variations of mixed Tate
  structures ([G5], s. 5.1)}. Let ${\cal V}$ be a variation of
  mixed Tate objects over a smooth base $B$ in one of our set--ups ii) --
  iv). So ${\cal V}$ is a unipotent variation of Hodge--Tate structures
  in ii), a lisse $l$-adic mixed Tate sheaf in iii), and a (yet
  hypothetical) variation of mixed Tate motives in iv). 
We assume that $\mu_{l^{\infty}} \not \subset 
{\cal O}^*(B)$ in iii). 

Then  ${\cal V}$ is an object of an appropriate mixed Tate 
category ${\cal C}_B$ of 
variations of mixed Tate objects over $B$. 
The Tate object $K(m)_B$ is given by 
$K(m)_B:= p^*K(m)$ where  $p: B \longrightarrow  {\rm
  Spec}(F)$ is  the structure morphism. 
We can use the standard formalism 
of mixed Tate categories (see [BD] or Appendix). 
In particular there is a commutative,
  graded Hopf $K$--algebra 
${\cal A}_{\bullet}^{\cal C}(B)$ with a coproduct $\Delta$, called the 
fundamental Hopf algebra of this
  category. One has 
$$
{\cal A}_1^{\cal C}(B) = {\cal O}^*(B)_ K
:= {\cal O}^*(B)\otimes_{\Q} K 
$$

For any positively graded (coassociative) 
Hopf algebra ${\cal
  A}_{\bullet}$ there is a canonical map
$$
\Delta^{[m]}: {\cal A}_m \longrightarrow \otimes^m {\cal A}_1
$$
Namely, it is dual to the multiplication map $\otimes^m {\cal
  A}^{\vee}_1 \longrightarrow {\cal A}_m^{\vee}$, and can be defined 
as the composition
$$
{\cal A}_m \stackrel{\Delta}{\longrightarrow}{\cal A}_{m-1} \otimes {\cal A}_1 
\stackrel{\Delta \otimes {\rm Id}}{\longrightarrow}{\cal A}_{m-2} \otimes {\cal
  A}_1 \otimes {\cal A}_1 \stackrel{\Delta \otimes {\rm Id}\otimes
  {\rm Id}}{\longrightarrow } ... \stackrel{\Delta \otimes {\rm
    Id}\otimes ...\otimes  
  {\rm Id}}{\longrightarrow }\otimes^m {\cal A}_1
$$

Let us choose a $K(m)_B$--framing of the variation ${\cal V}$.  
Then, by the very definition, we get an element
$$
[{\cal V}] \in {\cal A}_m^{\cal C}(B)
$$ 
 
\begin{definition} \label{5.8.02.1}
The element 
$$
\Delta^{[m]}\Bigl([{\cal V}] \Bigr) \in \otimes^m{\cal O}^*(B)_K
$$ 
is called the $\otimes^m$--invariant of a framed variation of mixed Tate
objects over $B$. 
\end{definition}

{\bf 4. The $\otimes^m$--invariant of the multiple logarithm variation 
and plane rooted trivalent trees}. 
Let ${\cal M}_{m+2}({\Bbb A}^1)$ be the space of $m+2$ ordered 
distinct
points $(a_0, ..., a_{m+1})$ on the affine line over a field. 
Every internal vertex $v$ of a plane trivalent rooted
tree $T$ provides an invertible 
 function $f^T_v$ on ${\cal M}_{m+2}({\Bbb A}^1)$. 
Indeed, if we picture a plane trivalent tree inscribed
into a circle such that the external vertices lie on this circle, 
the complement to the tree in the disc bounded by this
circle is a union of several connected domains. 
These domains are in bijective correspondence with the arcs on the
circle defined by the tree.

An internal  vertex $v$ determines 
 three such domains $D^v_1, D^v_2, D^v_3$, so that  
$v$ shares the boundaries of these domains. The plane structure of the
tree plus the orientation of the plane provide a cyclic order of
these domains, and the 
root of the tree provides a natural order of these domains. 
Namely, among  the edges sharing the vertex $v$ one edge is closer
to the root than the others. We count the domains going counterclockwise from
this edge. We
will assume that the enumeration $D^v_1, D^v_2, D^v_3$ reflects this order. 
 Let $a_1^v, a_2^v, a_3^v$ be the labels of the arcs assigned to the
 domains $D^v_1, D^v_2, D^v_3$. Set
$$
f^T_v:= \frac{a^v_3 - a^v_2}{a^v_1 - a^v_2} \in {\cal O}^*
({\cal M}_{m+2}({\Bbb A}^1))
$$

\begin{center}
\hspace{4.0cm}
\epsffile{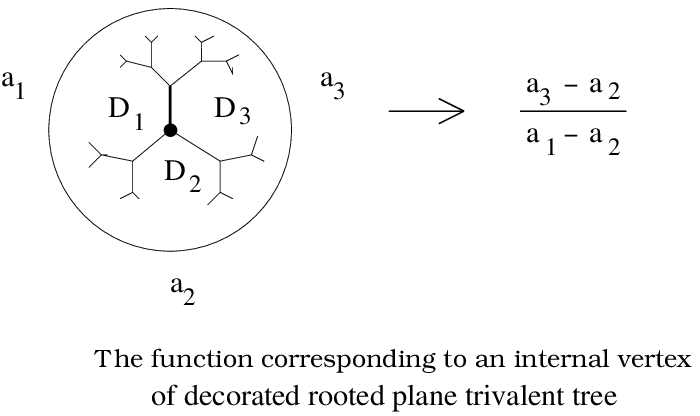}
\end{center}
On the internal vertices of a plane trivalent rooted tree
we have a {\it canonical partial order} $<$: we have $v_1< v_2$ 
if and only if 
there is a path going down from the root such that both  vertices
are located on this path, and  $v_1$ is closer to the
root than $v_2$. 
We say that an ordering $(v_1, ..., v_m)$ 
of internal vertices $v_i$ of a plane rooted tree
$T$ 
is {\it compatible} with the canonical order if $v_i < v_j$ implies $i<j$. 

\begin{definition} \label{4.8.02.2}
$$
\Omega_m:= \sum_T \sum _{\{v_1, ..., v_m\}} f^T_{v_1} \otimes ... \otimes f^T_{v_m}\in \otimes^m{\cal O}^*
({\cal M}_{m+2}({\Bbb A}^1))
$$
Here the first sum is over all different plane trivalent rooted trees $T$ 
with $m+1$ leaves, and the second sum is over all orderings 
$\{v_1, ..., v_m\}$ of the set of
internal vertices of the tree $T$ compatible with the canonical
partial order on this set.  
\end{definition}

\begin{proposition} \label{8.4.02.10} In the notations as above, we have 
$$
\Omega_m = \Delta^{[m]}{\rm I}^{\cal C}(a_0; a_1, ..., a_m; a_{m+1})
$$
\end{proposition}

{\bf Proof}. Follows easily from formula (\ref{CP2*}) by induction on $m$.

\section{The motivic multiple polylogarithm Hopf algebra}

 In this section we use formula (\ref{CP2*}) 
to calculate explicitly the coproduct for the motivic multiple polylogarithms. 
Then we define the  
cyclotomic Hopf and Lie algebras.

{\bf 1. The basic formula for the coproduct}. 
In this section we work with the Hopf algebra ${\cal I}_{\bullet}(S)$ where 
$S$ is any set containing a distinguished element $0$. Set, for $m \geq 0$ 
and $a_i \in S$, 
\begin{equation}  \label{12.18.00.2}
{\Bbb  I}_{n_0, n_1+1, ..., n_m+1}
(a_0; a_1, ..., a_m; a_{m+1}):= 
\end{equation} 
$$
{\Bbb  I}(a_0; \underbrace{0, ..., 0,}_{\mbox{$n_0$ times}} a_1, 
\underbrace{0, ..., 0,}_{\mbox{$n_1$ times}} a_2,\quad ...,  \quad 
\underbrace{0, ..., 0}_{\mbox{$n_m$ times}}; a_{m+1})
 $$
 
Let $V$ be a vector space. 
Denote by $V[[t_1, ..., t_m]]$ the vector space of formal power 
series in the variables $t_i$ whose coefficients are vectors of $V$. 
Let us form the generating series
\begin{equation} \label{CO1}
{\Bbb  I}(a_0; a_1, ... , a_m ; a_{m+1}
| t_0; t_1; ... ; t_m) := 
\end{equation}
$$
\sum_{n_i \geq 0}{\Bbb  I}_{n_0, n_1+1, ..., n_m+1}
(a_0; a_1, ... , a_m ; a_{m+1}) t_0^{n_0} ... t_m^{n_m} 
 \in \quad {\cal I}_{\bullet}(S)[[ t_0, ..., t_m ]]
$$
To visualize them  
consider a line segment with the following additional data, called 
{\it decoration}:

i) The beginning of the segment is labeled by $a_0$, the end by $a_{m+1}$.

ii) There are $m$ points {\it inside} of the segment 
labeled by $a_1, ..., a_m$ from 
 left to  right.

iii) These points cut the segment into $m+1$ arcs labeled  by 
$t_0, t_1, ..., t_m$.

\begin{center}
\hspace{4.0cm}
\epsffile{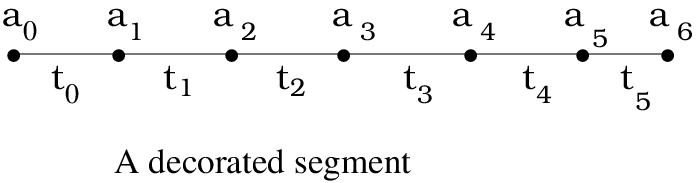}
\end{center}

{\bf Remark}. The way in which the $t$'s sit between the $a$'s reflects  
the shape of the iterated integral 
${\rm I}_{n_0, n_1+1, ..., n_m+1}(a_0; a_1, ..., a_m; a_{m+1})$. 

Our goal is to calculate the coproduct of the generating series (\ref{CO1}). 
Here, by definition, the coproduct acts on the coefficients of the series, 
leaving the $t$'s untouched:
$$
\Delta\Bigl(\sum_{n_i \geq 0}{\Bbb  I}_{n_0, n_1+1, ..., n_m+1}
(a_0; a_1, ... , a_m ; a_{m+1}) t_0^{n_0} ... t_m^{n_m} \Bigr):= 
$$
$$
\sum_{n_i \geq 0}\Delta\Bigl({\Bbb  I}_{n_0, n_1+1, ..., n_m+1}
(a_0; a_1, ... , a_m ; a_{m+1})\Bigr) t_0^{n_0} ... t_m^{n_m} 
$$
 As we  show in Theorem \ref{ur8-4,3} below, the 
terms of the coproduct of the element (\ref{CO1}) correspond to 
the decorated segments equipped with the following 
additional data, called marking:      

a) Mark (by making them fat points on the picture) 
points $a_0; a_{i_1}, ..., a_{i_{k}}; 
a_{m+1}$  
so that  
\begin{equation}  \label{ur11}
0 = i_0 < i_1 < ... < i_k < i_{k+1} = m+1
\end{equation} 

b) Mark (by crosses) segments $t_{j_0}, ..., t_{j_k}$ 
such that there is just one marked segment between any two 
neighboring marked points.

\begin{center}
\hspace{4.0cm}
\epsffile{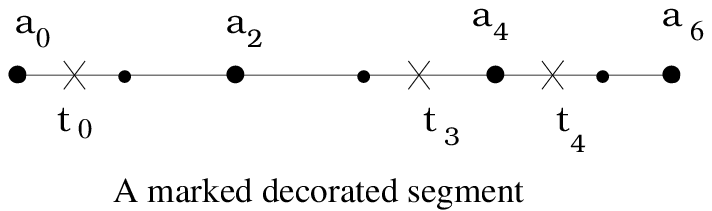}
\end{center}

The conditions on the crosses for a marked decorated segment just mean that 
\begin{equation}  \label{ur10}
i_{\alpha } \leq j_{\alpha } < i_{\alpha +1} \quad \mbox{ for any 
$0\leq \alpha \leq k$}
\end{equation}
The marks provide a new decorated segment:

\begin{equation}  \label{ur106}
(a_0|t_{j_0} |a_{i_0} | t_{j_1} |
a_{i_1} |\quad ... \quad | a_{i_k} |t_{j_k}|a_{m+1})
\end{equation}

\begin{theorem} \label{ur8-4,3} Let $a_i \in S, m\geq 0$. Then 
\begin{equation}  \label{ur100}
\Delta {\Bbb I}(a_0; a_1, ..., a_m; a_{m+1}| t_0; ...; t_{m}) 
= 
\end{equation}
$$
\sum {\Bbb I}(a_0; a_{i_1}, ..., a_{i_{k}}; 
a_{m+1}| t_{j_0}; t_{j_1}; ... ;
  t_{j_{k}}) \otimes
$$
$$
\prod_{\alpha =0}^k \Bigl(
{\Bbb I}(a_{i_{\alpha}}; a_{i_{\alpha}+1}, 
..., a_{j_{\alpha}};0| 
t_{i_{\alpha}}; ... ; t_{j_{\alpha}}) \cdot 
 $$
$$
{\Bbb I}(0; a_{j_{\alpha}+1}, ..., 
a_{i_{\alpha+1}-1}; a_{i_{\alpha+1}}|
t_{j_{\alpha}}; t_{j_{\alpha}+1}; ... ; t_{i_{\alpha+1}-1} )\Bigr)
$$
where the sum is over all marked decorated segments, i.e. over all 
sequences $\{i_{\alpha}\}$ and $\{j_{\alpha}\}$ 
satisfying inequality (\ref{ur10}). 
\end{theorem}

{\bf Proof}. Recall  that by iv) in section 2.1 one has 
\begin{equation} \label{12.18.00.1}
{\Bbb I}(0; a_1, ..., a_m; 0) =0
\end{equation}
To calculate (\ref{ur100}) we apply the coproduct formula (\ref{CP2**}) 
to the element 
(\ref{12.18.00.2}) 
and then keep track of the nonzero terms using (\ref{12.18.00.1}). 
 
The left hand side factors of the nonzero terms in the formula for the coproduct 
correspond to certain  subsequences
$$
A \subset \{a_0; \underbrace{0, ..., 0,}_{\mbox{$n_0$ times}} a_1, 
\underbrace{0, ..., 0,}_{\mbox{$n_1$ times}} \quad ...  \quad , a_m,
 \underbrace{0, ..., 0}_{\mbox{$n_m$ times}}; a_{m+1}   \}
$$ 
containing $a_0$ and $a_{m+1}$, and called the admissible subsequences. Such a 
subset $A$ determines the subsequences 
$
I = \{i_1 < ... < i_k\}\quad 
$
where $a_0, a_{i_1}, ..., a_{i_k}, a_{m+1}$ are precisely the set of all  $a_i$'s 
contained in $A$. 
A subsequence $A$ is called  {\it admissible} if 
 it satisfies the following properties:

i) $A$  contains $a_0$ and $a_{m+1}$.
 
ii) The sequence of $0$'s in $A$ located between 
$a_{i_{\alpha}}$ and $a_{i_{\alpha+1}}$ must be a string of {\it consecutive} $0$'s located 
between $a_{j_{\alpha}}$ and $a_{j_{\alpha}+1}$ for some $i_{\alpha} \leq j_{\alpha} 
<  i_{{\alpha}+1}$.

In other words, the factors in the coproduct are parametrized by 

{\it a marked decorated segment} and 
{\it a connected string of $0$'s in each of the crossed arcs}. 

The connected string of $0$'s in some  of the crossed arcs might be empty. 

The string of zeros 
between $a_{i_{\alpha}}$ and $a_{i_{\alpha}+1}$ satisfying ii) 
looks as follows: 
\begin{equation}  \label{ur8-4,5}
\{a_{i_{\alpha}}, \widehat 
a_{j_{\alpha}}, \underbrace{\widehat 0, ..., \widehat 0, 
}_{\mbox{$p_{j_{\alpha}}$}} \underbrace{0, ..., 0, 
}_{\mbox{$s_{j_{\alpha}}$}} \underbrace{\widehat 0, ..., 
\widehat 0, }_{\mbox{$q_{j_{\alpha}}$}} \widehat 
a_{j_{\alpha}+1}, a_{i_{\alpha+1}} \}
\end{equation}
where  $p_{j_{\alpha}} + q_{j_{\alpha}} + s_{j_{\alpha}} = n_{j_{\alpha}}$.
This notation emphasizes that all $0$'s located between $a_{i_{\alpha}}$ and 
$a_{i_{\alpha+1}}$ are in fact  located between $a_{j_{\alpha}}$ and 
$a_{j_{\alpha}+1}$, and form a connected segment of length $s_{j_{\alpha}}$. 

\begin{center}
\hspace{4.0cm}
\epsffile{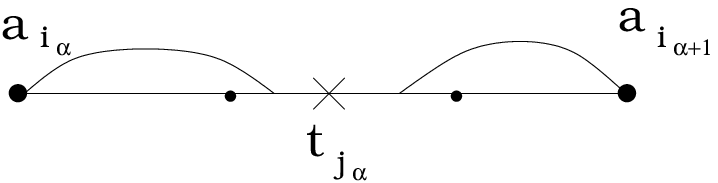}
\end{center}

 An admissible subset $A$ provides  the following element, 
which is the left hand side of the corresponding term in the coproduct: 
$$
{\Bbb I}(a_0; \widehat a_{j_{0}},  
\underbrace{0, ..., 0,}_{\mbox{$s_{j_0}$ times}}   \widehat a_{j_{0}+1}, 
a_{i_1}, \widehat a_{j_{1}},  \underbrace{0, ..., 0,}_{\mbox{$s_{j_1}$ times}}  
\widehat a_{j_{1}+1}, \quad ... \quad ; a_{m+1})
$$

The right hand side of the term in the coproduct corresponding to the subset $A$ 
is a product over $0 \leq \alpha \leq k$ of elements of the following shape:
$$
{\Bbb I}(a_{i_{\alpha}}; \underbrace{0, ..., 0, 
}_{\mbox{$n_{i_{\alpha}}$}} 
a_{i_{\alpha}+1}, \underbrace{0, ..., 0, }_{\mbox{$n_{
i_{\alpha}+1}$}} \quad ... \quad a_{j_{\alpha}}, 
\underbrace{0, ..., 0}_{\mbox{$p_{j_{\alpha}}$}};0) 
\cdot \underbrace{{\Bbb I}(0;0) \cdot ... 
\cdot {\Bbb I}(0;0)}_{s_{j_{\alpha}}}
$$
$$
{\Bbb I}(0; \underbrace{0, ..., 0, }_{\mbox{$q_{j_{\alpha}}$}} 
a_{j_{\alpha}+1}, \underbrace{0, ..., 0, }_{\mbox{$n_{j_{\alpha}+1}$}}  \quad ... \quad a_{i_{\alpha+1}-1}, \underbrace{0, ..., 0}_{\mbox{$n_{i_{\alpha+1}-1}$}};
a_{i_{\alpha+1}})
$$
where  the middle factor 
${\Bbb I}(0;0) \cdot ... \cdot {\Bbb I}(0;0)$ is equal to $1$ since ${\Bbb I}(0;0)= 1$ according to (\ref{12.18.00.1}). 
Translating this into the language of the generating series we get the promised 
formula for the coproduct.  The case $s_{j_{\alpha}}=0$ needs a special
treatment: the expression (\ref{ur8-4,5}) must be decomposed using the
path composition formula with $0$ as one of the endpoints. This
produces the correct formula for all admissible $j_{\alpha}$. 
 The theorem is proved. 

{\it A geometric interpretation of formula 
(\ref{ur100}) }. It is surprisingly similar to the 
 one  for the multiple logarithm element. 
Recall that the expression (\ref{CO1}) 
is encoded by a decorated segment 
\begin{equation}  \label{ur101} 
(a_0| t_0| a_1| t_1| \quad ... \quad |a_m|t_m|a_{m+1}) 
\end{equation}
The terms of the coproduct are in bijection with the marked decorated 
segments $C$ obtained from  a given one (\ref{ur101}). 
Denote by $L_C \otimes R_C$ the term in the coproduct corresponding to $C$. 
The left factor $L_C$ 
is encoded by the decorated segment (\ref{ur106}) obtained from the marked points and arcs. For example 
for the marked decorated segment on the picture  we get
$
L_C = {\Bbb I}(a_{0}|t_0|a_2|t_3|a_4|t_4|a_6)
$. 

The marks (which consist of $k+1$ crosses and $k+2$ boldface points) 
determine a decomposition of the segment (\ref{ur101}) 
into $2(k+1)$ 
little decorated segments in the following way. Cutting 
the initial segment in all the marked points and crosses we get 
$2(k+1)$ little segments. For instance, the very right one is 
the segment between the last cross and  point $a_m$, and so on.  
Each of these segments either starts
 from a marked point and ends by a cross, or starts from a 
cross and ends at a marked point. 

There is a natural way to make 
 a {\it decorated} segment out of each of these  
little segments: mark
 the ``cross endpoint'' of the little segment by the point $0$, and for the arc which 
is next  to this marked point use the 
letter originally attached to the arc containing it. 
For example the marked decorated segment on the picture above produces 
the following sequence of little decorated segments:
$$
(a_0|t_0|0), \quad (0|t_0|a_1|t_1|a_2), \quad (a_2|t_2|a_3|t_3|0), \quad (0|t_3| 
a_4) \quad (a_4| t_4|0), \quad (0|t_4|a_5|t_5|a_6)
$$

Then the factor $R_C$ is the product of the generating series for the elements 
corresponding to these little decorated segments. 
For example for the marked decorated segment on the picture we get 
$$
R_C = {\Bbb I}(a_0;  0|t_0) \cdot {\Bbb I}(0; a_1; a_2|t_0; t_1) 
\cdot 
$$
$$
{\Bbb I}(a_2; a_3; 0|t_2;t_3) \cdot {\Bbb I}(0; a_4|t_3)  \cdot 
{\Bbb I}(a_4; 0|t_4) \cdot {\Bbb I}(0; a_5; a_6|t_4; t_5)
$$

To check that the formula for the coproduct of the multiple logarithms fits 
into this description, we use the path composition formula iii) 
in  Section 2.1,  
provided by Theorem \ref{7.24.02.12}b),   
together with  the fact that 
each term  on  the right hand side of this formula corresponds to a marked 
colored segment shown in the picture:
\begin{center}
\hspace{4.0cm}
\epsffile{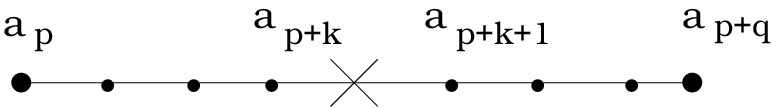}
\end{center}

{\bf 2. The multiple polylogarithm Hopf algebra of a subgroup $G \subset F^*$}. 
Now let us suppose that we work in one of the set--ups i)--iv), 
and $F$ is the corresponding field.  Let $S:= F$. 
We define the ${\rm I}^{\cal C}$--elements 
as the images of the corresponding ${\Bbb I}$--elements  
by the homomorphism established in Theorem \ref{7.24.02.12}. 
Below we choose a coordinate $t$ on the affine line and  use the standard 
tangent vectors $v_s$ dual to  $dt$. 
By Theorem   \ref{7.24.02.12}, 
the results of the previous section are valid for the corresponding 
${\rm I}^{\cal C}$--elements. 

The ${\rm I}^{\cal C}$--analogs of 
elements (\ref{12.18.00.2}) are clearly invariant under translation. 
In general  they are 
not  
invariant under the 
action of ${\Bbb G}_m$ given by 
 $a_i \lms \lambda a_i$. However, if the iterated integral 
${\rm I}_{n_0, n_1+1, ..., n_m+1}(a_0; a_1, ... , a_m ; a_{m+1})$ is convergent, then 
the corresponding element is invariant under the 
action of the affine group, see Lemma \ref{9.30.00.2} below. 

Let $a \in F^*$. We set 
$$
[a]:= \log^{\cal C}(a) := {\rm I}^{\cal C}(0; 0; a) 
\in {\cal A}_{1}({\cal C})
$$
In any of the set-ups i) - iii), and hypothetically in 
iv),  
we have an   isomorphism $F^{*}\otimes_{\Q}
 K\stackrel{\sim}{=}{\cal A}_{1}({\cal C})$ 
materialized by the ``motivic logarithm'':
$$
\log^{{\cal C}}:  F^* \otimes_{\Q} K \lra  {\rm Ext}_{{\cal C}}^1(K(0), K(1)) 
 \stackrel{\sim}{=}  {\cal A}_{1}({\cal C}); \qquad a \lms 
\log^{{\cal C}}(a)
$$  
Its $l$-adic realization is given by the Kummer extension, 
and if $F \subset \C$ its Hodge 
realization is  $\log^{\cal H} (a)$. See Appendix, subsection 3, especially formula 
(\ref{8-23/99gs}). Here it is essential that we use the tangential
base 
point $\partial/\partial t$, where $t$ is the standard parameter on
the line. Then the fact that we got the above isomorphism is a basic
standard fact in each of the set-ups. One also has $\log^{\cal C}(a) =
{\rm I}^{\cal C}(1;0;a)$, but we will not use it. Observe that we are
working here in realizations, rather then with the formal elements ${\Bbb
  I}$. Set ($\widehat \otimes$ stands for the completed tensor product)
$$
a^t:= {\rm exp}(\log^{\cal C}(a)\cdot t) \in 
{\cal A}_{\bullet}({\cal C})\widehat \otimes k[[t]]
$$
The shuffle product formula implies 
$$
{\rm I}^{\cal C}(0; a|t) = a^{t}
$$
One has $$
\Delta(1) = 1 \otimes 1, \quad \Delta([a]) = 
[a] \otimes 1 +  1 \otimes [a]
 $$
It follows that  
\begin{equation} \label{d0f}
\Delta(a^t) =   a^t\otimes a^t
\end{equation}

Suppose that  $a_i \not = 0$. Set 
$$
{\rm I}^{\cal C}_{ n_1, ..., n_m}(a_1, ... , a_m):=
$$
$$
{\rm I}^{\cal C}(0; a_1, \underbrace{0, ..., 0}_{\mbox{$n_1-1$}}, a_2, 
\underbrace{0, ..., 0}_{\mbox{$n_2-1$}}, ... , a_m, 
\underbrace{0, ..., 0}_{\mbox{$n_m-1$}}; 1)
$$
Let us 
package them into the generating series
\begin{equation} \label{CO1e}
{\rm I}^{\cal C}( a_1, ... , a_m|  t_1, ... ,t_m) := 
\end{equation}
$$
\sum_{n_i \geq 1}{\rm I}^{\cal C}_{n_1, ..., n_m}(a_1, ... , a_m) 
t_1^{n_1-1} ... t_m^{n_m-1} 
 \in \quad {\cal A}_{\bullet}({\cal C})[[t_1, ..., t_m]] 
$$

\begin{lemma} \label{9.30.00.2u} Suppose that $a_i \not = 0$. Then  
\begin{equation}  \label{3.7.01.19}
{\rm I}^{\cal C}(0; a_1, ... , a_m; a_{m+1}| t_0; ... ;t_{m}) = 
\end{equation}
$$
a_{m+1}^{t_0} {\rm I}^{\cal C}(\frac{a_1}{a_{m+1}}, ... , \frac{a_m}{a_{m+1}}| t_1-t_0, ..., t_m-t_0)
$$
\end{lemma}

{\bf Proof}. It follows by induction 
using the shuffle product formula 
for 
$$
\sum_{n \geq 0} 
{\rm I}^{\cal C}(0; \underbrace{0, ..., 0}_{\mbox{$n$ times}};
 a_{m+1})t_0^{n}\cdot \sum_{n_i > 0}
{\rm I}_{n_1, ..., n_m}^{\cal C}(0; a_1, ..., a_m; a_{m+1})t^{n_1-1}_1
 ... t^{n_m-1}_m
$$
just as in the second proof of Proposition 2.15 in [G3].
The lemma is proved. 

\begin{lemma} \label{9.30.00.2} One has  
$$
{\rm I}^{\cal C}(a_1; a_2, ..., a_m; 0| t_1; ...; t_m) = 
$$
\begin{equation}  \label{ur8-9}
(-1)^{m-1}{\rm I}^{\cal C}(0; a_m, ..., a_2; a_1| -t_m; -t_{m-1}; 
  ...;  - t_1)  
\end{equation}
In particular $
{\rm I}^{\cal C}(a; 0|t) = a^{-t}$.
\end{lemma}

{\bf Proof}. This is a special case of the equality 
$$
{\rm I}^{\cal C}(a_1; a_2, ..., a_m; 0) = (-1)^{m-1} 
{\rm I}^{\cal C}(0; a_m, ..., a_2; a_1)
$$
provided by Proposition \ref{8.15.02.1} and Theorem \ref{7.24.02.12}. The lemma is proved.

A marked decorated segment is {\it special} if the first cross is marking 
the segment $t_0$. 
A decorated segment with such a  data is called a {\it special marked 
decorated segment}. 
See an example on the picture after formula \ref{ur11}. 
The conditions on the crosses for a marked decorated segment just mean that 
\begin{equation}  \label{ur10*}
i_{\alpha } \leq j_{\alpha } < i_{\alpha +1} \quad \mbox{ for any 
$0\leq \alpha \leq k$}, \qquad j_0 = i_0 = 0
\end{equation}

We will employ the notation
$$
{\rm I}^{\cal C}(a_1: ... : a_{m+1}| t_0: ... :t_{m}):= \quad 
{\rm I}^{\cal C}(0; a_1, ... , a_m; a_{m+1}| t_0; ... ;t_{m})
$$

\begin{theorem} \label{7.26.02.111} One has, for any $a_i$'s,  
$$
\Delta {\rm I}^{\cal C}(a_1: ... : a_{m+1}| t_0: ... :t_{m}) =
$$
$$
\sum {\rm I}^{\cal C}( a_{i_1}: ...: a_{i_{k}}: a_{m+1}| t_{j_0}: t_{j_1}: ... :
  t_{j_{k}}) \otimes
$$
\begin{equation}  \label{ur100100}
\prod_{\alpha =0}^k \Bigl(
(-1)^{j_{\alpha}-i_{\alpha}} {\rm I}^{\cal C}(a_{j_{\alpha}}: a_{j_{\alpha}-1}: ... : 
 a_{i_{\alpha}}| 
-t_{j_{\alpha}}: -t_{j_{\alpha}-1}: ... : - t_{i_{\alpha}}) \cdot 
\end{equation}
$$
   {\rm I}^{\cal C}(a_{j_{\alpha}+1}: ... : a_{i_{\alpha+1}-1}: a_{i_{\alpha+1}}|
t_{j_{\alpha}}: t_{j_{\alpha}+1}: ... : t_{i_{\alpha+1}-1} )\Bigr)
$$
where the sum is over all special marked 
decorated segments, i.e. over all 
sequences $\{i_{\alpha}\}$ and $\{j_{\alpha}\}$ 
such that 
\begin{equation}  \label{ur10*}
i_{\alpha } \leq j_{\alpha } < i_{\alpha +1} \quad \mbox{ for any 
$0\leq \alpha \leq k$}, \qquad j_0 = i_0 = 0, i_{k+1} = i_{m+1} 
\end{equation}
\end{theorem}

{\bf Proof}. It follows immediately from Theorem \ref{ur8-4,3} using Lemmas 
\ref{9.30.00.2u} and \ref{9.30.00.2}. 

Now let  $G$ be a subgroup of $F^*$. We assume that $F \subset \C$ in the Hodge case,
and $F$ is the base field otherwise. 
Denote by ${\cal Z}^{\cal C}_{w}(G) \subset {\cal A}_{w}({\cal C})$ 
the $\Q$-vector subspace 
generated by the elements
\begin{equation} \label{1.11.01.1}
{\rm I}^{\cal C}_{n_1, ..., n_m} (a_1, ..., a_m), \quad a_i \in G, 
\quad w = n_1 + ... + n_m
\end{equation}
Set 
$$
{\cal Z}^{\cal C}_{\bullet}(G):= 
\oplus_{ w \geq 1} {\cal Z}^{\cal C}_{w}(G)
$$

\begin{definition} \label{8.4.03.1}
A depth filtration ${\cal F}_{ \bullet}^{{\cal D}}$ on the space 
${\cal Z}^{\cal C}_{\bullet}(G)$ is 
defined as follows: 

${\cal F}_{ 0}^{{\cal D}}{\cal Z}^
{{\cal C}}_{\bullet}(G) $ is spanned by  products of 
$\log^{\cal C} (a)$, $a \in G$, and 

$
{\cal F}_{ k}^{{\cal D}}{\cal Z}^{{\cal C}}_{\bullet}(G) 
$ for $k \geq 1$ is spanned by the elements  (\ref{1.11.01.1}) with 
$m \leq k$. 
\end{definition}
Since 
$\log^{{\cal C}} (a) = {\rm I}^{{\cal C}}(0; 0; a)$ is 
the depth zero multiple polylogarithm, this agrees with this definition.

\begin{theorem} \label{9.30.00.3}
Let $G$ be any subgroup of  $F^*$. Then 
${\cal Z}^{\cal C}_{\bullet}(G)$ is a graded 
Hopf subalgebra of ${\cal A}_{\bullet}({\cal C})$. 
The depth provides a filtration on this Hopf algebra.
\end{theorem}

{\bf Proof}. 
The graded vector space ${\cal Z}^{\cal C}_{\bullet}(G)$ is 
closed under the coproduct by  Theorem \ref{7.26.02.111}.
The statement about the depth 
filtration is evident from the formula for the coproduct given in 
Theorem \ref{7.26.02.111}. 
It is a graded algebra  by the shuffle product formula. 
The theorem  is proved. 

{\bf Remark}. The depth filtration is not defined by a grading of the algebra
${\cal Z}^{\cal C}_{\bullet}(G)$ because of  
relations like 
$$
{\rm Li}^{\cal C}_{n}(x) \cdot  {\rm Li}^{\cal C}_{m}(y) \quad = \quad 
{\rm Li}^{\cal C}_{n,m}(x,y) + 
{\rm Li}^{\cal C}_{m,n}(y,x) +  {\rm Li}^{\cal C}_{n+m}(xy)
$$

{\bf 3. An example: the cyclotomic Lie algebras}. Let $\mu_N$ be the group 
of $N$-th roots of unity. In this subsection we work with the abelian category ${\cal M}_T(\Q(\zeta_N))$ of mixed Tate motives over the cyclotomic field. Thus we use 
the superscript 
${\cal M}$ instead of ${\cal C}$.
Set 
$$
 {\cal C}^{{\cal M}}_{\bullet}(\mu_N):= \quad 
\frac{{\cal Z}^{{\cal M}}_{>0}(\mu_N)}{
{\cal Z}^{{\cal M}}_{>0}(\mu_N) \cdot {\cal Z}^{{\cal M}}_{>0}(\mu_N)}
$$
 \begin{corollary} \label{ur8-11} 
a) ${\cal C}^{{\cal M}}_{\bullet}(\mu_N)$ is a graded Lie coalgebra. 
The depth provides a filtration on this Lie coalgebra.

b).
${\cal C}^{{\cal M}}_{1}(\mu_N) = {\cal Z}^{{\cal M}}_{1}(\mu_N) 
\stackrel{\sim}{=} 
\Bigl(\mbox{the group of cyclotomic units 
in $\Z[\zeta_N][1/N]\Bigr) \otimes \Q$}$.
\end{corollary}

{\bf Proof}. a). Follows from theorem \ref{9.30.00.3}.  

b) 
Since $$
{\rm I}^{{\cal M}}_1(a) =  \log^{{\cal M}} (1-a^{-1})\quad \mbox{and} \quad  
\log^{{\cal M}} (1-a) -  \log^{{\cal M}} (1-a^{-1}) =  \log^{{\cal M}} (a)
$$
 the  weight $1$ component ${\cal Z}^{{\cal M}}_{1}(G)$
is generated by $\log^{{\cal M}} (1-a)$ and $ \log^{{\cal M}} (a)$. 
Notice that if $a^N=1$ then 
$N \cdot  \log^{{\cal M}} (a) = 0$. This proves b). The corollary is proved.

We call 
${\cal C}^{{\cal M}}_{\bullet}(\mu_N)$ the 
{\it cyclotomic Lie coalgebra}.   Its dual 
${C}^{{\cal M}}_{\bullet}(\mu_N)$ is the {\it cyclotomic Lie algebra}. 
The dual to the 
universal enveloping algebra  of the cyclotomic Lie algebra is 
isomorphic to ${\cal Z}^{{\cal M}}_{\bullet}(\mu_N)$.

Let ${\cal C}_{m}^{\Delta}(\mu_N)$ be the $\Q$--subspace of  
${\cal C}^{{\cal M}}_{m}(\mu_N)$ generated by all 
$
{\rm I}^{{\cal M}}(0; a_1, ..., a_m; 1)$ where $ a_i^N =1$. 

\begin{corollary} \label{ur8-11/1} ${\cal C}_{\bullet}^{\Delta}(\mu_N):= \oplus_{m\geq 1}{\cal C}_{m}^{\Delta}(\mu_N)$ 
is a graded Lie coalgebra. 
\end{corollary}

We call it the  
 {\it diagonal cyclotomic Lie coalgebra}. 

{\bf Proof}. Clear from the  Theorem \ref{7.26.02.111}.

The depth filtration on the cyclotomic Lie algebra plays a central 
role in the mysterious correspondence between the structure of this 
Lie algebra and the 
geometry of certain modular varieties for $GL_m, m=1,2,3,4, ...$. 
([G1-2]). However our definition of the depth filtration depends on  the
realization of the Hopf algebra ${\cal Z}^{\cal M}_{\bullet}(G)$ 
related to the projective line.  In the next section we suggest 
that the depth filtration is induced, via the canonical embedding 
${\cal Z}^{\cal M}_{\bullet}(G) \hookrightarrow {\cal A}_{\bullet}(F)$, 
 by a natural 
filtration on the motivic Tate Hopf algebra of $F$,  
given intrinsically in terms of the corresponding Lie algebra. 

{\bf 4. A hypothetical  intrinsic definition of the depth filtration}. 
In this subsection we 
assume the existence of the {\it abelian} category ${\cal M}_T(F)$ 
of mixed Tate motives over
 an arbitrary field $F$ with all the standard properties. In particular 
${\cal M}_T(F)$ is a mixed Tate category, so 
there are the corresponding fundamental Lie algebra
 $L_{\bullet}(F)$, its graded dual 
${\cal L}_{\bullet}(F)$, and 
${\cal A}_{\bullet}(F)$ is isomorphic 
to the universal enveloping algebra of $L_{\bullet}(F)$, see the Appendix. 

We assume that for any  $a_1, ..., a_m \in F^*$ there exists an 
object ${\rm I}^{{\cal M}}_{n_1, ..., n_m}(a_1, ..., a_m)$ 
of the category ${\cal M}_T(F)$ framed by $\Q(0)$ and $\Q(w)$, 
where $w = n_1+...+n_m$, called the {\em motivic 
multiple polylogarithm}, and its coproduct  $\Delta
{\rm I}^{{\cal M}}_{n_1, ..., n_m}(a_1, ..., a_m)$ is given by the 
formula from Theorem \ref{7.26.02.111}.  
When $F$ is a number field all this is already available.

{\it The universality conjecture}. 
It is the  Conjecture 17a) in [G6], which tells us  that 
every 
framed mixed Tate motive over $F$ is equivalent to a $\Q$-linear combination of 
the ones ${\rm I}^{{\cal M}}_{n_1, ..., n_m}(a_1, ..., a_m)$ 
with $a_i \in F^*$. 
The universality conjecture  can be reformulated as follows:

\begin{conjecture} \label{1.10.01.3} The Hopf algebra ${\cal Z}^{\cal M}_{\bullet}(F^*)$ 
is isomorphic to the fundamental 
Hopf algebra ${\cal A}_{\bullet}(F)$ of the abelian 
category of mixed Tate motives 
over $F$. 
\end{conjecture}

In other words 
 every framed mixed Tate motive over $F$ is a $\Q$--linear 
combination of the ones coming from 
the motivic torsor of paths ${\cal P}^{{\cal M}}({\Bbb A}^1_{/F}; v_0, v_1)$ between 
the standard tangential base points at $0$ and $1$.

{\it The main conjecture}. 
Since ${\cal Z}_{\bullet}^{\cal M}(F^*)$ is a Hopf subalgebra in 
${\cal A}_{\bullet}(F)$, we can consider the depth filtration 
provided by Definition \ref{8.4.03.1} on 
the former as a filtration on the latter. 
It follows from the 
formula in Theorem \ref{7.26.02.111} that the depth filtration is a filtration 
by coideals  in  ${\cal A}_{\bullet}(F)$. 
The depth filtration on ${\cal A}_{\bullet}(F)$ induces the 
depth filtration on the 
Lie coalgebra 
${\cal L}_{\bullet}(F)$. In particular ${\cal F}_{0}^{\cal D}{\cal L}_{\bullet}(F) = 
{\cal L}_{1}(F) = F^* \otimes \Q$. 

Let us define another filtration on the Lie algebra ${L}_{\bullet}(F)$. 
Consider the following ideal of $L_{\bullet}(F)$:
$$
 I_{\bullet}(F):= \oplus^{\infty}_{n=2}  L_{-n}(F)
$$
We define the {\it depth filtration} $\widetilde {\cal F}^D$ 
on the Lie algebra $L_{\bullet}(F)$ as an increasing  filtration  
indexed by integers $m \leq 0$ and given by the powers of the ideal $I_{\bullet}(F)$
\begin{equation} \label{1.10.01.1}
\widetilde {\cal F}^{D}_0L_{\bullet}(F) = L_{\bullet}(F); \quad 
\widetilde {\cal F}^{D}_{-1}L(F)_{\bullet} = I_{\bullet}(F) 
; \quad \widetilde {\cal F}^{D}_{-m-1}L_{\bullet}(F) = [I_{\bullet}(F), \widetilde {\cal F}^{D}_{-m}L(F)_{\bullet}]
\end{equation}
Thus there are two filtrations on the Lie coalgebra ${\cal L}_{\bullet}(F)$:  the 
dual to the  depth 
filtration (\ref{1.10.01.1}) and the 
filtration  by the depth of multiple polylogarithms provided by the Definition 
\ref{8.4.03.1}. 
 
\begin{conjecture} \label{13.1.01.1}
The dual to the depth filtration (\ref{1.10.01.1}) coincides with the 
filtration provided by the Definition 
\ref{8.4.03.1}. 
\end{conjecture}
This conjecture, of course, implies Conjecture \ref{1.10.01.3}. 
%When $F$ is a number field, it implies  Zagier's conjecture [Z]. 

{\it The role of the classical polylogarithms}. Let  ${\cal B}_{n}(F)$ be 
the $\Q$-vector space in 
${\cal L}_{n}(F)$ spanned by the classical $n$-logarithm framed motives 
${\rm I}^{{\cal M}}_n(a)$, $a \in F^*$. (This definition differs 
from the usual one given in [G4], although
 the corresponding groups  are expected to be isomorphic.) 
Let $H^{i}_{(n)}$ be the degree $n$ part of $H^{i}$.
The following conjecture was stated in [G4]:

\begin{conjecture} \label{7.3}
a)  $H^{1}_{(n)}I_{\bullet}(F) \cong {\cal
B}_{n}(F)$ for $n\geq 2$, i.e. $ I_{\bullet}(F)$ 
is generated as a graded Lie algebra by the spaces ${\cal
B}_{n}(F)^{\vee}$ sitting in degree $-n$.

b) $I_{\bullet}(F)$  is a free graded (pro) -
Lie algebra.
\end{conjecture}

Conjecture \ref{13.1.01.1} obviously implies the part a) of Conjecture \ref{7.3}.

The following result, which is a motivic version of Theorem 2.22 in [G3], 
 shows  that all  $n$-framed mixed Tate motives but perhaps a countable set are 
given by  motivic multiple polylogarithms. 
This is  strong support for the universality conjecture. 
\begin{proposition} \label{3.2} Let us assume the existence of the category 
of mixed motivic sheaves with all the expected properties. 
Let ${\Bbb V}$ be a variation of n-framed mixed Tate motives over a 
connected rational variety $Y$. For an $F$-point $y \in Y$ 
the fiber $V_y$ provides  a framed 
mixed Tate motive, and hence an element of $[V_y] \in {\cal A}_{\bullet}(F)$. 
Then for any two points $y_1,y_2 \in Y$ the difference $[V_{y_1}] - [V_{y_2}]$ 
is a sum of motivic multiple polylogarithms.
\end{proposition}

{\bf Proof}. One can suppose without loss of generality that there 
is a rational curve $X$ passing through $y_1$ and $y_2$. Let 
${\cal P}^{{\cal M}}(X; y_1,y_2)$ be the motivic torsor of  paths 
from $y_1$ to $y_2$ on $X$. 
During the proof of Theorem \ref{7.26.02.111} 
we calculated explicitly the mixed Tate motive corresponding to an arbitrary 
framing on 
${\cal P}^{{\cal M}}(X; y_1,y_2)$, and showed that it is equivalent to  
a sum of motivic multiple polylogarithms. This sum is given by the 
right hand side of the main formula in Theorem \ref{7.26.02.111}. Observe that product of motivic multiple polylogarithms is expressible as a sum of 
motivic multiple polylogarithms.  The 
parallel transport along paths from $y_1$ to $y_2$
on $X$ provides a morphism of 
mixed Tate motives $p_{y_1,y_2}^X :V_{y_1}\otimes
 {\cal P}^{{\cal M}}(X;y_1,y_2)
\rightarrow V_{y_2}$. Let ${\cal A}$ be the kernel of the action of ${L}_{\bullet}(F)$
on   ${\cal P}^{{\cal M}}(P^1; v_0, v_1)$. Then ${\cal P}^{{\cal M}}(X; y_1,y_2)$
  is a trivial  ${\cal A}$ - module. Therefore any $p \in  {\cal P}^{{\cal M}}
(X; y_1,y_2)$ defines an isomorphism of ${\cal A}$ - modules $V_{y_1}\otimes p 
\rightarrow V_{y_2}$. This is equivalent to the statement of the proposition.

\section{Examples and applications}

{\bf 1. The ${\rm Li}^{\cal C}$--generators of the 
multiple polylogarithm Hopf algebra}. 
Let us define several other generating series 
for the  multiple polylogarithm elements. 
Consider the following two pairs of  sets of variables:
$$
i) \quad (x_0, ... , x_m) \quad \mbox{such that $x_0 ...  x_m=1$}; \qquad ii) 
\quad (a_1: ... : a_{m+1})
$$
$$
iii) \quad (t_0: ... : t_m); \qquad iv) 
\quad (u_1, ..., u_{m+1}) \quad \mbox{such that $u_1 +  ... + u_{m+1} = 0$}
$$
The relationship between them is given by 

$$
x_i := \frac{a_{i+1}}{a_i}, \quad i = 1, ..., m; \quad 
x_0 = \frac{a_1}{a_{m+1}}
$$
$$
u_i := t_i -t_{i-1}, \quad u_{m+1} := t_0 - t_m
$$
where the indices are taken modulo $m+1$, and is illustrated on the picture below

\begin{center}
\hspace{4.0cm}
\epsffile{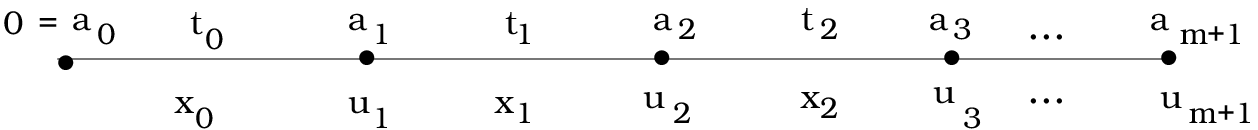}
\end{center}
Observe that $x_i, a_i$ are multiplicative variables, 
and $t_i, u_i$ are additive variables. 

We introduce the ${\rm Li}^{\cal C}$-generating series
\begin{equation} \label{11.16.00.1}
{\rm Li}^{\cal C}(*, x_1, ..., x_m| t_0: ... : t_m):= \quad 
{\rm Li}^{\cal C}(x_0, ..., x_m| t_0 : ... : t_m) := 
\end{equation}
$$
(-1)^m 
{\rm I}^{\cal C}( (x_1 ... x_m)^{-1}: 
(x_2 ... x_m)^{-1}: ... : x^{-1}_m :1 | t_0 : ... : t_m) 
$$
$$
=: (-1)^m 
{\rm I}^{\cal C}( a_1: a_2: ... : a_{m+1}| t_0 : ... : t_m)  \quad =: 
{\rm I}^{\cal C}( a_1 : a_2 : ... : a_{m+1}| u_1,  ... u_{m+1})
$$
{\bf Remark}. The $(,)$-notation is used for the variables which sum to zero 
(under the appropriate group structure), and the $(:)$-notation is used for those 
sets of variables which are essentially homogeneous with respect to 
 the multiplication by a common factor, 
see Lemma \ref{9.30.00.2}. The $*$ in the left expression in (\ref{11.16.00.1}) 
denotes $(x_1 ... x_m)^{-1}$. This notation is handy in some situations. 

{\bf 2.  The coproduct in terms of the ${\rm Li}^{\cal C}$--generating series}. 
In this section we rewrite the formula 
for the coproduct using the ${\rm Li}^{\cal C}$--generators instead of 
the ${\rm I}^{\cal C}$--generators. Set
$$
X_{a\to b}:= \quad \prod_{s=a}^{b-1} x_{s}
$$

\begin{proposition} \label{ur8-4,3**} Let us suppose that $x_i \not = 0$. Then 
$$
\Delta {\rm Li}^{\cal C}(x_0, x_1, ..., x_m | t_0: t_1: ... :t_{m}) \quad = \quad 
$$
\begin{equation}  \label{ur10000}
\sum 
{\rm Li}^{\cal C}(X_{i_0 \to i_{1}}, X_{i_1 \to i_{2}}, ... , X_{i_k \to m}| 
t_{j_0}: t_{j_1}:  ... :
  t_{j_{k}}) \otimes
\end{equation}
\begin{equation}  \label{12-15.00.1fr}
\prod_{p = 0}^k \Bigl( (-1)^{j_p-i_p}X_{i_p \to i_{p+1}}^{t_{j_p}}
{\rm Li}^{\cal C} (*, x^{-1}_{{j_p}-1}, x^{-1}_{{j_p}-2}, ..., x^{-1}_{i_p}|-t_{j_{p}}:  
-t_{j_{p}-1}:  ... : -t_{i_{p}}) \cdot 
\end{equation}
\begin{equation}  \label{12-15.00.2}
{\rm Li}^{\cal C} (*, x_{j_{p}+1}, x_{j_{p}+2}, ..., x_{i_{p+1}-1}| 
t_{j_{p}}: t_{j_{p}+1}:  ... : t_{i_{p+1}-1} )\Bigr)  
\end{equation}
Here the sum is over special marked decorated segments, i.e.
sequences $\{i_{p}\}$, $\{j_{p}\}$ 
satisfying (\ref{ur10*}). 
\end{proposition}

The $p$-th factor in the product on the right is encoded by the data 
on the $p$-th segment:

\begin{center}
\hspace{4.0cm}
\epsffile{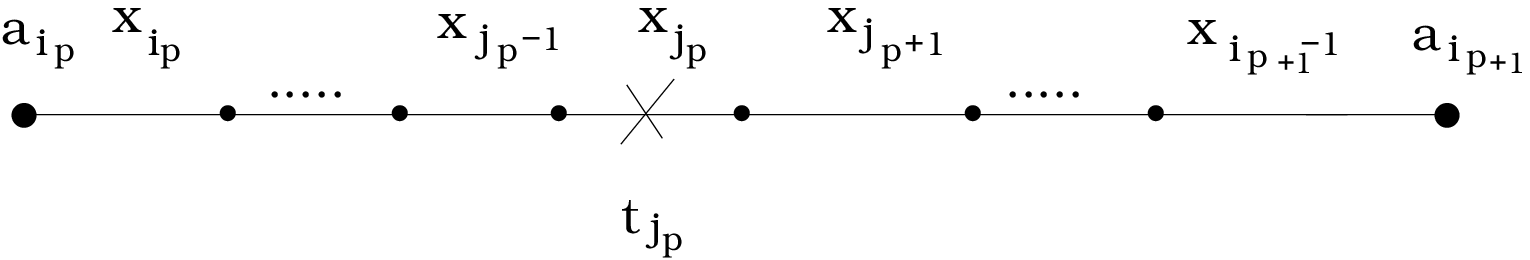}
\end{center}
 Namely, $X_{i_p \to i_{p+1}}$ is the product of all $x_i$ on this segment. 
The factor (\ref{12-15.00.1fr}) 
is encoded by the segment  between $t_{j_p}$ and $t_{i_p}$, which we read from 
right to  
left. The factor (\ref{12-15.00.2})  is encoded by 
the segment  between $t_{j_p}$ and 
$t_{i_{p+1}-1}$ which we read from left to right.

{\bf Proof}. Follows from Lemma \ref{9.30.00.2u} and 
Theorem \ref{7.26.02.111}.

{\bf 3. The coproduct in the classical polylogarithm case}. 
Recall that ${\cal A}_{\bullet}({\cal C})$ is the commutative 
Hopf algebra of the framed mixed  objects in a mixed Tate category ${\cal C}$ 
with the coproduct $\Delta$ and the product $\ast$, see 
the Appendix, s.2 (where the notation $\mu$ for the product is used). 
Recall 
the restricted coproduct: $\Delta'(X): = \Delta(X)
- (X\otimes 1 + 1\otimes X)$. Notice that $\Delta$ is a
homomorphism of algebras and $\Delta'$ is not.

\begin{corollary}  \label{9.30.00.1}
\begin{equation} \label{d1f}
\Delta: {{\rm Li}}^{\cal C}(x|t) \longmapsto { {\rm Li}}^{\cal C}(x|t
)\otimes x^{t} + 1 \otimes  {\rm Li}^{\cal C}(x|t)
\end{equation}
\end{corollary}

{\bf Proof}. This is a special case of Proposition \ref{ur8-4,3**}. 

This formula just means that 
$$
\Delta' { {\rm Li}}^{\cal C}_n(x) = 
{ {\rm Li}}^{\cal C}_{n-1}(x) \otimes 
\log^{\cal C} x + {{\rm Li}}^{\cal C}_{n-2}(x) 
\otimes \frac{(\log^{\cal C} x)^2}{2} + ... + 
{{\rm Li}}^{\cal C}_{1}(x) \otimes 
\frac{(\log^{\cal C} x)^{n-1}}{(n-1)!}
$$

{\bf 4. The coproduct for the depth two
 multiple polylogarithms}. 
We will use both types of the ${\rm I}$--notations 
for multiple polylogarithm elements, so for instance 
$$
{\rm I}^{\cal C}(a_1:a_2:1|t_1,t_2) \quad = \quad 
{\rm I}^{\cal C}(0; a_1, a_2; 1|t_1;t_2)
$$ 
Set $\zeta^{\cal C}(t_1, ..., t_m) :=  {\rm I}^{\cal C}
(1: ... : 1 | t_1, ..., t_m)$.

 \begin{proposition} \label{ur107} a) One has 
$$
\Delta {\rm I}^{\cal C}(a_1:a_2:1|t_1,t_2) \quad = 1 \otimes 
{\rm I}^{\cal C}(a_1:a_2:1|t_1,t_2) +
$$
$$
{\rm I}^{\cal C}(a_1:a_2:1|t_1,t_2) \otimes a_1^{-t_1}\ast a_2^{t_1-t_2} \quad +
 \quad {\rm I}^{\cal C}(a_1:1|t_1) \otimes a_1^{-t_1}\ast  
{\rm I}^{\cal C}(a_2:1|t_2-t_1)
$$
$$
-{\rm I}^{\cal C}(a_1:1|t_2) \otimes a_1^{-t_2} \ast 
{\rm I}^{\cal C}(a_2:a_1|t_2-t_1)
 \quad + \quad  {\rm I}^{\cal C}(a_2:1|t_2) \otimes  
{\rm I}^{\cal C}(a_1:a_2|t_1)\ast   a_2^{-t_2}
$$

b) Let us suppose that $a_1^N = a_2^N = 1$. Then modulo  $N$-torsion 
one has 
$$
\Delta' {\rm I}^{\cal C}(a_1:a_2:1|t_1,t_2) \quad = \quad 
{\rm I}^{\cal C}(a_1:1|t_1) \otimes   {\rm I}^{\cal C}(a_2:1|t_2-t_1)
$$
$$
- {\rm I}^{\cal C}(a_1:1|t_2) \otimes {\rm I}^{\cal C}(a_2:a_1|t_2-t_1)
 \quad + \quad  {\rm I}^{\cal C}(a_2:1|t_2) \otimes 
{\rm I}^{\cal C}(a_1:a_2|t_1)
$$
In particular
\begin{equation} \label{9.4.02.2}
\Delta' \zeta^{\cal C}(t_1,t_2) \quad = \quad  \zeta^{\cal C}(t_1) \otimes 
\zeta^{\cal C}(t_2-t_1) - \zeta^{\cal C}(t_2) \otimes 
 \zeta^{\cal C}(t_2-t_1) + 
\zeta^{\cal C}(t_2) \otimes 
\zeta^{\cal C}(t_1)
\end{equation}
\end{proposition}

{\bf Proof}. a) Since in our case $a_0=0$ and  $t_0=0$, 
the nonzero contribution can be obtained only from 
those marked decorated segments where the $t_0$-arc is not marked. 
Let us call pictures where $a_0=0$, $a_{m+1} =1$,  $t_0=0$ and 
the $t_0$-arc is not marked, {\it special} marked decorated segments. 

The five  terms in the formula above  correspond to the five special 
marked decorated segments presented in the picture. 
\begin{center}
\hspace{4.0cm}
\epsffile{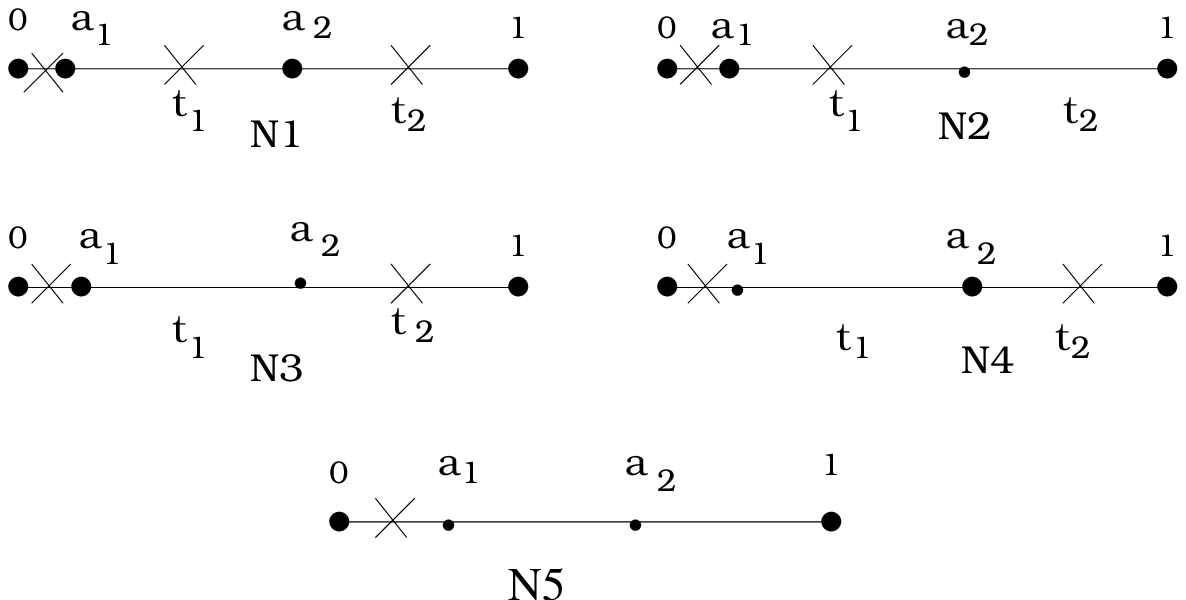}
\end{center}

Using the formulas from Lemma \ref{9.30.00.2}
we get the following four terms corresponding to the terms N1-N4 on the 
picture: 
$$
 {\rm I}^{\cal C}(a_1:a_2:1|t_1,t_2) \otimes 
 {\rm I}^{\cal C}(a_1;0|t_1) \cdot  {\rm I}^{\cal C}(0;a_2|t_1) \cdot 
{\rm I}^{\cal C}(a_2;0|t_2)\cdot  {\rm I}^{\cal C}(0;1|t_2) = 
$$
$$
{\rm I}^{\cal C}(a_1:a_2:1|t_1,t_2) \otimes a_1^{-t_1} \cdot a_2^{t_1-t_2}
$$

$$
{\rm I}^{\cal C}(a_1:1|t_1) \otimes {\rm I}(a_1;0|t_1) \cdot 
  {\rm I}^{\cal C}(0; a_2; 1|t_1;t_2) \quad = \quad 
{\rm I}^{\cal C}(a_1:1|t_1)\otimes a_1^{-t_1} \cdot  {\rm I}^{\cal C}(a_2:1|t_2-t_1)
$$

$$
{\rm I}^{\cal C} (a_1:1|t_2) \otimes {\rm I}^{\cal C}(a_1;a_2;0|t_1;t_2) 
\cdot  {\rm I}^{\cal C}(0;  1|t_2) = 
$$
$$
- {\rm I}^{\cal C}(a_1:1|t_2) \otimes {\rm I}(0;a_2;a_1|-t_2;-t_1) \quad = \quad - 
{\rm I}^{\cal C}(a_1:1|t_2) \otimes a_1^{-t_2}\cdot {\rm I}^{\cal C}(a_2: a_1|t_2-t_1) 
$$

$$
{\rm I}^{\cal C}(0;a_2;1|0;t_2) \otimes {\rm I}^{\cal C}(0; a_1;a_2|0;t_1) \cdot  {\rm I}^{\cal C}(a_2; 0|t_2) 
\cdot  {\rm I}^{\cal C}(0;1|t_2) =
$$
$$
{\rm I}^{\cal C}(a_2: 1|t_2) \otimes {\rm I}^{\cal C}(a_1 : a_2|t_1) \cdot a_2^{-t_2}
$$

Part b) follows from a) if we notice that $a^N =1$ 
provides $a^t =1$ modulo  $N$-torsion. 
The proposition is proved. 

{\bf Remark}. In Theorem 4.5 of [G9] 
the reader can find a different way to write the 
formulas for the coproduct in  depth $2$ case. It is easy to see that the formulas 
given there are equivalent to the formulas above.

{\bf 5. Explicit formulas for  the coproduct of the weight three, 
 depth two multiple  polylogarithm elements}.  
Applying Proposition \ref{ur107} or Proposition \ref{ur8-4,3**} we get 
$$
\Delta': {{\rm Li}}^{\cal C}_{2,1}(x, y) \quad  \longmapsto 
\quad{ {\rm Li}}^{\cal C}_{1,1}(x, y)\otimes [x] +
{{\rm Li}}^{\cal C}_{1}(y) \otimes {{\rm Li}}^{\cal C}_{2}(x)
 +{{\rm Li}}^{\cal C}_{2}(x y) \otimes {{\rm Li}}^{\cal C}_{1}(y)
$$
$$
- {{\rm Li}}^{\cal C}_{1}(x y) \otimes 
\Bigl({{\rm Li}}^{\cal C}_{2}(x) + 
{{\rm Li}}^{\cal C}_{2}(y) - 
{{\rm Li}}^{\cal C}_{1}(y)\cdot ([xy] + \frac{[x]^2}{2})\Bigr) 
$$

$$
\Delta': { {\rm Li}}^{\cal C}_{1,2}(x, y)  \quad   
\longmapsto \quad { {\rm Li}}^{\cal C}_{1,1}(x, y)\otimes [y] -
{ {\rm Li}}^{\cal C}_{2}(xy)\otimes [x] - 
{ {\rm Li}}^{\cal C}_{1}(xy)\otimes [x y]
$$
$$
+ {{\rm Li}}^{\cal C}_{2}(y) \otimes { {\rm Li}}^{\cal C}_{1}(x) + {{\rm Li}}^{\cal C}_{1}(y) \otimes {{\rm Li}}^{\cal C}_{1}(x)\cdot [y] 
+ {{\rm Li}}^{\cal C}_{1}(x y) \otimes  {{\rm Li}}^{\cal C}_{2}(y)  
$$
$$
- {{\rm Li}}^{\cal C}_{2}(x y) \otimes {{\rm Li}}^{\cal C}_{1}(x) - 
{{\rm Li}}^{\cal C}_{1}(x y) \otimes  {{\rm Li}}^{\cal C}_{1}(x)\cdot [xy] - 
{{\rm Li}}^{\cal C}_{1}(x y) \otimes  {{\rm Li}}^{\cal C}_{2}(x)   
$$

{\bf 6. The depth two 
Lie coalgebra structure}. 
Recall that 
the space of the indecomposables
$$
{\cal L}_{\bullet}({\cal C}):= \quad \frac{{\cal A}_{>0}({\cal C})}
{{\cal A}_{>0}({\cal C}) \cdot {\cal A}_{>0}({\cal C})}
$$ 
inherits  a structure of  graded Lie coalgebra with the cobracket $\delta$.  
I will describe
$$%\begin{equation}  \label{12e}
\sum _{n,  m > 0} \delta  {\rm I}^{\cal C}_{n, m}
(a_1, a_2) \cdot t_1^{n - 1 } 
t_2^{m - 1 } \quad \in \quad \Lambda^2
{\cal A}_{\bullet}({\cal C}) \widehat \otimes {\Z} [[t_1, t_2]]
$$
In the formulas below $\delta$ acts  on the first factor in 
${\cal A}_{\bullet}({\cal C}) \widehat \otimes {\Z} [[t_1, t_2]]$.
Using Proposition \ref{ur107} we have
 
$$%\begin{equation}  \label{12*}
 \delta \bigr(\sum_{m> 0,n> 0}  {{\rm I}}_{m,n}(a:b:c)\cdot t_{1}^{m- 1}t_{2}^{n-1}\bigl) 
 \quad = \quad
$$%\end{equation}
$$
\sum_{m> 0, n > 0} \biggr({{\rm I}}_{m,n}(a:b:c)\cdot t_{1}^{m- 1}t_{2}^{n-1})
\wedge ([\frac{b}{a}]\cdot t_{1} + [\frac{c}{b}]\cdot t_{2}) \quad - \quad  
{ {\rm I}}_{m}(a:b) \cdot t_{1}^{m- 1} 
\wedge 
{{\rm I}}_{n}(b:c) 
\cdot t_{2}^{n-1} 
$$
$$
 + \quad  {{\rm I}}_{m}(a:c) \cdot t_{1}^{m- 1} \wedge {{\rm I}}_{n}(b:c) 
\cdot (t_{2} - t_{1})^{n-1} 
\quad - \quad {{\rm I}}_{n}(a:c) \cdot t_{2}^{n- 1} \wedge 
{{\rm I}}_{m}(b:a) 
\cdot (t_{2} - t_{1})^{m - 1} \biggl) 
$$

Set  ${\rm I}_{0,n} = {\rm I}_{n,0} =0$.  
Here is a more concrete formula for  $\delta$:  
$$%\begin{equation}  \label{14*}
\delta  { {\rm I}}_{m,n}(a,b)
\quad = 
$$%\end{equation}
$$
\quad {{\rm I}}_{m-1,n}(a,b) \wedge [\frac{b}{a}] \quad + \quad
{ {\rm I}}_{m,n-1}(a,b) \wedge [\frac{1}{b}] \quad 
- \quad {{\rm I}}_{m}(\frac{a}{b}) \wedge {{\rm I}}_{n}(b) \quad +
$$
$$
\sum_{i=0}^{m-1} (-1)^{i}  { n+ i - 1 \choose i }  {{\rm I}}_{m - i}(a) \wedge { {\rm I}}_{n + i}(b) 
$$
$$
 - (-1)^{m-1} \sum_{j=0}^{n-1} (-1)^{j}  { m+ j - 1 \choose j }  {{\rm I}}_{n - j}(a) \wedge {{\rm I}}_{m + j}(\frac{b}{a})
$$

{\bf 7. An application: the motivic double $\zeta$'s}. Below we use the notation 
${\cal Z}_{\bullet}^{\cal M}$ for the cyclotomic Hopf algebra 
${\cal Z}_{\bullet}^{\cal M}(\mu_N)$ in the case $N=1$. 
Recall the restricted coproduct  map 
\begin{equation} \label{8.8.03.2}
\Delta': {\cal Z}_{\bullet}^{\cal M} \lra {\cal Z}_{\bullet}^{\cal M}\otimes 
{\cal Z}_{\bullet}^{\cal M}
\end{equation} 
\begin{theorem} \label{8.8.03.1}
The degree $n$ part of ${\rm Ker}\Delta'$ is zero if $n=1$ or if $n$ is even, 
and it is one dimensional otherwise. It is spanned by $\zeta^{\cal M}(n)$.
\end{theorem}

{\bf Proof}. The motivic multiple zeta values are defined as 
the matrix coefficients of the motivic torsor of path 
${\cal P}^{\cal M}(P^1-\{0,1,\infty\}; v_0, v_1)$, which  itself is
a mixed Tate motive over ${\rm Spec}(\Z)$ ([DG]). Therefore 
motivic multiple zeta's belong to ${\cal A}_{\bullet}(\Z)$. 
 So the kernel of the map (\ref{8.8.03.2}) is a subspace of the 
kernel of the map $\Delta': {\cal A}_{\bullet}(\Z) \lra 
{\cal A}_{\bullet}(\Z)\otimes {\cal A}_{\bullet}(\Z)$. The degree $n$ part 
of the  latter 
is isomorphic to 
\begin{equation} \label{SEE}
{\rm Ext}^1_{{\cal M}_T(\Z)}(\Q(0), \Q(n)) = K_{2n-1}(\Z)\otimes \Q = 
\left\{ \begin{array}{ll}
0 &  \quad \mbox{$n$: even} \\
 1 &    \quad \mbox{$n>1$: odd} \end{array} \right. 
\end{equation}
On the other hand, the element  $\zeta^{\cal M}(n)$ is non-zero 
for the same values of $n$, and $\Delta'\zeta^{\cal M}(n)=0$. So it 
generates the degree $n$ part of ${\rm Ker}\Delta'$ in (\ref{8.8.03.2}). 
The theorem is proved. 

Recall the depth filtration 
 ${\cal F}^{\cal D}$ on the Hopf algebra ${\cal Z}_{\bullet}^{\cal M}$. 
The $m$-th 
associate graded for the depth filtration ${\cal Z}_{\bullet, m}^{\cal M}:= 
{\rm gr}^{\cal D}_m{\cal Z}_{\bullet}^{\cal M}$ and 
${\cal C}_{\bullet, m}^{\cal M}:= 
{\rm gr}^{\cal D}_m{\cal C}_{\bullet}^{\cal M}$ are graded by the weight. 

Restricting the coproduct to the subspace of the depth two motivic multiple zeta's 
and using Theorem \ref{8.8.03.1} and formula (\ref{9.4.02.2}) we get an embedding
$$
\Delta': {\cal Z}_{\bullet, 2}^{\cal M}  \hookrightarrow 
{\cal Z}_{\bullet, 1}^{\cal M}
\otimes {\cal Z}_{\bullet, 1}^{\cal M}
$$ 
Observe that if $\Delta'(x) = \Delta'(y)=0$ then $\Delta'(xy)= 
x\otimes y + y \otimes x$. 
Therefore the product 
provides an 
inclusion 
$
S^2{\cal F}^{\cal D}_1{\cal Z}_{\bullet}^{\cal M} \hookrightarrow
{\cal F}^{\cal D}_{\leq 2}{\cal Z}_{\bullet}^{\cal M}$, as well as 
an inclusion $
S^2{\cal Z}_{\bullet, 1}^{\cal M} \hookrightarrow 
{\cal Z}_{\bullet, 2}^{\cal M}$. 
 Recall $$
{\cal C}_{\bullet, 2}^{\cal M}= 
{\cal Z}_{\bullet, 2}^{\cal M}/ S^2{\cal Z}_{\bullet, 1}^{\cal M}, \qquad {\cal C}_{\bullet, 1}^{\cal M}:= {\cal Z}_{\bullet, 1}^{\cal M} 
$$
Then $\Delta'$ induces an injective map 
\begin{equation} \label{8.8.03.6}
\delta: {\cal C}_{\bullet, 2}^{\cal M} \hookrightarrow \Lambda^2
{\cal C}_{\bullet, 1}^{\cal M}
\end{equation}
Denote by $\overline \zeta^{\cal M}(m,n)$ the projection of 
$\zeta^{\cal M}(m,n)$ onto 
${\cal C}_{m+n, 2}^{\cal M}$.
 The formula (\ref{9.4.02.2}) tells us that:
\begin{equation} \label{9.4.02.2sa}
\delta \overline \zeta^{\cal M}(t_1,t_2) \quad = \quad  \zeta^{\cal M}(t_1) \wedge
\zeta^{\cal M}(t_2-t_1) - \zeta^{\cal M}(t_2) \wedge
 \zeta^{\cal M}(t_2-t_1) + 
\zeta^{\cal M}(t_2) \wedge
\zeta^{\cal M}(t_1) 
\end{equation}
Since $\zeta^{\cal M}(t) =  \zeta^{\cal M}(-t)$, we can write (\ref{9.4.02.2sa}), 
employing a linear map $ U: (t_1, t_2) \lms (t_1-t_2, t_1)$, as 
\begin{equation} \label{9.4.02.2sasa}
\delta:  \overline \zeta^{\cal M}(t_1,t_2) \lra
-\Bigl(I + U + U^2\Bigr)\zeta^{\cal M}(t_1) 
\wedge\zeta^{\cal M}(t_2)
\end{equation}
Consider the generating series 
$
\overline \zeta^{\cal M}(t_0, t_1, t_2) := \overline \zeta^{\cal M}(t_1, t_2), 
$ where $t_0+t_1+t_2 =0$.  

\begin{theorem} \label{8.8.03.4}
a) The generating series $\overline \zeta^{\cal M}(t_0, t_1, t_2)$ 
satisfy the dihedral symmetry relations:
\begin{equation} \label{DS}
\overline \zeta^{\cal M}(t_0, t_1, t_2) = \overline \zeta^{\cal M}(t_1, t_2, t_0) =
-\overline \zeta^{\cal M}(t_0, t_2, t_1) = \overline \zeta^{\cal M}(-t_0, -t_1, -t_2)
\end{equation}

b) There are no other relations between the coefficients of the generating series 
$\overline \zeta^{\cal M}(t_0, t_1, t_2)$. 
\end{theorem}

{\bf Proof}. a) The dihedral symmetry 
relations are evidently in the kernel of the map (\ref{8.8.03.6}). 
Since this map is injective we are done. 

b). Consider the classical modular triangulation 
of the hyperbolic  plane obtained by the action of the group $GL_2(\Z)$ 
on the ideal triangle $T_{0, 1, \infty}$ with the vertices at $0, 1, \infty$. 
Let $$
M_{(2)}^*:= M_{(2)}^1 \stackrel{\partial}{\lra} M_{(2)}^2
 $$
be the chain complex of this triangulation sitting in the degrees $[1,2]$. 
It is a complex of 
$\Z[GL_2(\Z)]$--modules 
where the first group is generated by the class $[T_{0, 1, \infty}]$ 
of the oriented triangle $T_{0, 1, \infty}$, and the second one 
by the class $[E_{0, \infty}]$ of the oriented geodesic from $0$ to $\infty$. 

It is handy to introduce a new bigraded Lie coalgebra defined as a direct sum 
of Lie coalgebras
$$
\widehat {\cal C}_{\bullet, \bullet}^{\cal M}:= {\cal C}_{\bullet, \bullet}^{\cal M}
\oplus \Q^{\cal M}_{1,1}
$$
where $\Q^{\cal M}_{1,1}$ is the one-dimensional Lie coalgebra 
in bidegree $(1,1)$ 
spanned by a new formal element  
$\zeta^{\cal M}(1)$. This element does not have any motivic meaning. 
The map (\ref{9.4.02.2sa}) provides a map
\begin{equation} \label{DS1}
{\cal C}_{\bullet, 2}^{\cal M} \lra \Lambda^2 \widehat {\cal C}_{\bullet, 1}^{\cal M}
\end{equation}

\begin{lemma} \label{CD4} 
The formulas
$$
\mu^1: [T_{0, 1, \infty}] \otimes t_1^{m-1}t_2^{n-1} \lms \zeta^{\cal M}(m,n) 
t_1^{m-1}t_2^{n-1}; 
\qquad 
\mu^2: [E_{0, \infty}] \otimes t_1^{m-1}t_2^{n-1} \lms \zeta^{\cal M}(m) \wedge 
\zeta^{\cal M}(n) t_1^{m-1}t_2^{n-1}
$$
provide  an isomorphism of complexes (\ref{DS2}), where 
$\Lambda^2_w$ means the degree $w$ part of $\Lambda^2$:  
\begin{equation} \label{DS2}
\begin{array}{ccc}
M_{(2)}^1 \otimes_{GL_2(\Z)}S^{w-2}V_2 & \stackrel{\partial}{\lra} &
M_{(2)}^2\otimes_{GL_2(\Z)}S^{w-2}V_2 \\
&&\\
\mu^1\downarrow &&\downarrow \mu^2\\
&&\\
{\cal C}_{\bullet, 2}^{\cal M} &\stackrel{\delta}{\hookrightarrow} &\Lambda^2_w
\widehat {\cal C}_{\bullet, 1}^{\cal M}
\end{array}
\end{equation}\end{lemma}

{\bf Proof}. 
The dihedral 
group is the 
stabiliser of the triangle $T_{0,1,\infty}$ in $GL_2(\Z)$. 
Thus the relations (\ref{DS}) just mean 
that the first formula gives a well defined map of $M_{(2)}^1
\otimes_{GL_2(\Z)}S^{w-2}V_2$. The relation 
$\zeta^{\cal M}(t) = \zeta^{\cal M}(-t) $ is equivalent to a 
similar fact about the second formula. And formula (\ref{9.4.02.2sasa}) just means 
that the constructed maps commute with the differentials. 
The lemma is proved. 

By its very definition the maps $\mu^1$ and $\mu^2$ are surjective. 
The map $\mu^2$ 
 is an isomorphism. Indeed, the stabilizer of $E_{0, \infty}$ in
 $GL_2(\Z)$ is the
 subgroup generated by $\left (\matrix{\pm 1&
     0\cr 0& \pm 1\cr}\right )$ and $\left (\matrix{0& 1\cr
     -1&0\cr}\right )$.  So  taking the coinvariants of the action of this
 subgroup on $S^{w-2}V_2$ amounts to the relation $\zeta^{\cal M}(t) =
 \zeta^{\cal M}(-t)$ and the skew-symmetry of $\Lambda^2$. 
Since the map (\ref{DS1}) is injective, 
to check that $\mu^1$ is injective it is 
sufficient to show that the top arrow $\partial$ is injective. Let us
prove this. 
Let $\varepsilon_2$ be the determinant representation of $GL_2(\Z)$. By Lemma 2.3 from [G1] one has 
$$
H^i\Bigl(M_{(2)}^*\otimes_{GL_2(\Z)}S^{w-2}V_2 \Bigr) = H^{i-1}(GL_2(\Z), 
S^{w-2}V_2 \otimes \varepsilon_2), \quad i = 1,2
$$
Since $H^{0}(GL_2(\Z), 
S^{w-2}V_2 \otimes \varepsilon_2) =0$, the map $\partial$ is injective. The theorem is proved. 

{\bf Remark}. The relations (\ref{DS}) are the motivic counterparts 
of the (regularized) double shuffle relations for the 
double zeta values considered modulo products and depth one terms, 
see Section 2 of [G2] for the form of the double shuffle relations making this 
statement obvious.  Therefore Theorem \ref{8.8.03.4} plus 
Conjecture \ref{8.04.03.2} imply 
that there should be no other relations 
between the double zeta's. This was conjectured by Zagier [Z]. 

{\bf Example}. Formula (\ref{9.4.02.2})
 tells us that 
\begin{equation} \label{9.4.02.4}
\Delta' \zeta^{\cal M}(3,5)  = -5\cdot \zeta^{\cal M}(3) \otimes \zeta^{\cal M}(5)
\end{equation}
Since  $\zeta^{\cal M}(2n+1) \not = 0$ this implies that 
$\zeta^{\cal M}(3,5)\not = 0$.  Moreover since 
$\delta' \overline \zeta^{\cal M}(3,5) = -5\cdot \zeta^{\cal M}(3) \wedge
\zeta^{\cal M}(5) $ is also non-zero, $\zeta^{\cal M}(3,5)$ is irreducible, i.e. it is not a product of the classical motivic zeta's. Thus it is a 
non-zero, non-classical, irreducible, motivic period over $\Z$ of the  
smallest weight.  

{\bf 8. When does a motivic itegrated integral unramified at ${\cal P}$?}. 
Let ${\cal O}_F$ be the ring of integers in a number field $F$. 
Let $S$ be a collection of ideals of ${\cal O}_F$. Recall from Section 8.4 
the fundamental Hopf algebra ${\cal A}_{\bullet}({\cal O}_{F, S})$ of the category 
of mixed Tate motives over ${\rm Spec}({\cal O}_{F, S})$. 

Recall that our definition of the motivic iterated integrals 
\begin{equation} \label{6.12.04.2}
{\rm I}^{\cal M}(a_0; a_1, ..., a_n; a_{n+1}) \in {\cal A}_n(F)
\end{equation}
assumes a choice of a tangent vector at every point $a_i$. Below we assume that 
all these tangent vectors equal $\partial/\partial t$, where $t$ is 
the chosen coordinate on the affine line ${\Bbb A}^1$. 
To emphasise this choice we use the notation ${\rm I}_{\partial/\partial t}^{\cal M}(a_0; a_1, ..., a_n; a_{n+1})$. These elements are invariant under the shift 
$a_i \to a_i +c$. Denote by $v_{\cal P}: F^* \to \Z$ the valuation at the prime ideal ${\cal P}$. 

\begin{theorem} \label{6.12.04.1}
Let $a_0, ..., a_{n+1}$ be elements of a number field $F$. 
Then ${\rm I}_{\partial/\partial t}^{\cal M}(a_0; a_1, ..., a_n; 
a_{n+1})$ is unramified at a place ${\cal P}$ of $F$ if 
for any $0 \leq i < j < k \leq n+1$ one has 
$
v_{\cal P}\widetilde r(a_i, a_j, a_k) =0
$. 
\end{theorem}

{\bf Proof}. Recall the inductive definition of the 
space ${\cal A}_n({\cal O}_{F, S})$ given by (\ref{8-24.7/99fr}). 
We prove the theorem by the induction on $n$. Let $n=1$. Then one has 
$$
{\rm I}_{\partial/\partial t}^{\cal M}(a_0;a_1;a_2) = \widetilde r(a_0; a_1; a_2) 
$$
(The right hand side is defined in the Example in Chapter 1.2). 
Since ${\cal A}_1({\cal O}_{F, S}) \stackrel{\sim}{=} {\cal O}_{F, S}^* \otimes \Q$, 
this proves the theorem for $n=1$. Let us show using the induction assumption 
that our  condition implies 
\begin{equation}\label{6.12.04.3}
\Delta'{\rm I}_{\partial/\partial t}^{\cal M}(a_0; a_1, ..., a_n; a_{n+1}) \in \oplus_{k=1}^{n-1}
{\cal A}_k({\cal O}_{F, S}) \otimes {\cal A}_{n-k}({\cal O}_{F, S}) 
\end{equation}
Indeed, it is clear from the formula (2) that every factor in (\ref{6.12.04.3}) 
is of type ${\rm I}_{\partial/\partial t}^{\cal M}(a_{i}; a_{p_1}, ..., a_{p_k}; a_{j})$ 
where $0 \leq i < p_1 < ... < p_k < j \leq n+1$. Moreover one has $k<n$ 
since we consider the restricted coproduct $\Delta'$. Therefore each of these factors belongs 
${\cal A}_n({\cal O}_{F, S})$ by the induction assumption. The theorem is proved.

\section{Feynman integrals, Feynman diagrams  and mixed motives}

{\bf 1. Motivic correlators}. 
Let us return to Chapter 4. 
Why do plane trivalent trees 
appear in the description of motivic iterated integrals, 
and what is the general framework for this relationship? 

Previously, we defined in Chapters 8 and 9 of [G2] 
a real--valued version of multiple polylogarithms as correlators 
of a 
Feynman integral. In fact we gave there a more general construction which 
provides a {\it definition} of a real valued version of 
multiple polylogarithms on an arbitrary smooth curve. 
Moreover the Feynman diagram 
 construction given in [G2] provides much more than just 
functions -- it  can be lifted to a construction of 
 framed mixed motives 
whose periods  
 are the correlators of the corresponding Feynman integral. 
 A version of Conjecture \ref{1.10.01.3} claims that 
all framed mixed Tate motives appear this way.

So why do Feynman integrals and Feynman 
diagrams appear 
in the theory of mixed motives, and what role do they play there? 
I think we see in the two  examples above
a manifestation of the following general principle.

Feynman integrals are described by their correlators. 
Correlators of Feynman integrals 
are often periods of framed mixed motives: see an example in Chapter 9 of [G2].  
In this case, the correlators 
can be upgraded to  more  
sophisticated objects: the equivalence classes of the 
corresponding framed mixed motives (see Appendix or 
Chapter 3 of [G7] for the background). 
These objects lie in a certain commutative Hopf algebra ${\cal H}_{\rm Mot}$ 
with the coproduct $\Delta_{\rm Mot}$. 
One can loosely think about ${\cal H}_{\rm Mot}$ 
as of the algebra of regular functions on the motivic Galois group. 
More precisely, ${\cal H}_{\rm Mot}$ is a Hopf algebra in the 
hypothetical  
abelian tensor category ${\cal P}_{\cal M}$ of all pure motives, see [G7].  
 The motivic Galois group is a pro--affine group scheme 
in this category, see [D2]. 
Although ${\cal H}_{\rm Mot}$ 
is still a hypothetical mathematical object, we can  see it 
in different existing realizations, e.g. in the Hodge realization. 

Climbing up  the road 
$$
\mbox{correlators of Feynman integrals (numbers)} \quad \lra 
\quad 
$$
$$
\mbox{motivic correlators (elements of the Hopf algebra ${\cal H}_{\rm Mot}$)}
$$
we gain a new perspective: one can now raise the question 
\begin{equation} \label{8.23.02.1}
\mbox{ what is the coproduct 
of motivic correlators in ${\cal H}_{\rm Mot}$}?
\end{equation}
Reflections on this theme occupy the rest of the Chapter.

{\bf 2. The correspondence principle}. Let us formulate question 
(\ref{8.23.02.1}) more precisely. 
Let 
\begin{equation} \label{8.23.02.5}
{\rm Cor}= \langle\varphi (s_1), ..., \varphi (s_{m+1})\rangle 
\end{equation}
be a correlator of a certain Feynman integral ${\cal F}$. 
(Here $\varphi$ is the traditional notation for fields 
over which we integrate). 
If we understand our Feynman integral 
by its perturbation series expansion then according to the  
Feynman rules correlator (\ref{8.23.02.5}) 
 is defined  
as a sum of finite dimensional integrals: 
\begin{equation} \label{8.23.02.7}
{\rm Cor} := 
\sum_{\Gamma \in {\cal S}_{\cal F}({Cor})} \int_{X_{\Gamma}}
\omega_{\Gamma}
\end{equation}
The sum is over a (finite) set ${\cal S}_{\cal F}({Cor})$ 
of certain combinatorial objects  
$\Gamma$, given by decorated graphs. It is  determined by 
the Feynman integral ${\cal F}$ and the type of the correlator we consider.  
Such a  $\Gamma$ 
provides a real algebraic variety $X_{\Gamma}$ and a differential form of top degree 
$\omega_{\Gamma}$ on $X_{\Gamma}(\R)$. 

Let us assume that the integrals in  (\ref{8.23.02.7}) are convergent. 
This is often not the case in physically interesting examples. 
However it is the case for the Feynman integral considered in [G2]. 
Then correlator (\ref{8.23.02.5}) is a well defined number.

Let us assume further that this number is a period of a mixed motive. 
This means that it is given by a sum of integrals of 
rational  differential forms on certain varieties (which may differ from 
$X_{\Gamma}\otimes \C$) 
over certain chains whose boundaries lie in 
a union of divisors. Below we consider only Feynman integrals whose correlators 
have this property, called Feynman integrals of algebraic-geometric type. 
A nontrivial example of such a Feynman integral is given in [G2].

Then we conjecture that one can uniquely  upgrade this number to 
the corresponding ``motivic correlator''
\begin{equation} \label{8.23.02.6}
{\rm Cor}_{\cal M} \in {\cal H}_{\rm Mot}
\end{equation}
 It is a framed mixed motive. The   period of its  
Hodge realization  is given by (\ref{8.23.02.5}). 
We ask in (\ref{8.23.02.1}) about the coproduct
\begin{equation} \label{8.23.02.8}
\Delta_{\cal M}\Bigl({\rm Cor}_{\cal M}\Bigr)
 \in {\cal H}_{\rm Mot}\otimes {\cal H}_{\rm Mot}
\end{equation}

We suggest that {\it the answer 
should be 
given combinatorially in  terms of the decorated graphs $\Gamma$ 
used in the Definition (\ref{8.23.02.7})}. Here is a more 
precise version of this guess.

Let $H_1$ and $H_2$ be Hopf algebras in tensor categories $T_1$ and $T_2$. 
 A Hopf algebra homomorphism 
$H_1 \to H_2$ is given by a tensor functor $F: T_1 \to T_2$ 
and a Hopf algebra homomorphism $F(H_1) \to H_2$.

{\bf The correspondence principle}. {\it 
For a Feynman integral ${\cal F}$ of algebraic-geometric type 
there should exist a combinatorially 
defined commutative Hopf algebra $({\cal H}_{\cal F}, \Delta_{\cal F})$
 in a tensor category $T_{\cal F}$ such that 
the right hand side of (\ref{8.23.02.7})  
provides an element 
\begin{equation} \label{8.24.02.1}
 \gamma:= 
\sum_{\Gamma \in {\cal S}_{\cal F}({\rm Cor})}[\Gamma] \in {\cal H}_{\cal F}
\end{equation}
The map  
\begin{equation} \label{8.23.02.9}
c_{\cal M}:\gamma \lms
 {\rm Cor}_{\cal M}
\end{equation}
 gives rise to a Hopf algebra homomorphism
$$
c_{\cal M}: {\cal H}_{\cal F} \lra {\cal H}_{\rm Mot}
$$ 
provided by  a tensor functor $F: T_{\cal F} \lra {\cal P}_{\cal M}$}

In particular  $c_{\cal M}$ 
 is compatible with the coproducts:
\begin{equation} \label{8.29.02.1}
(c_{\cal M} \otimes c_{\cal M}) \Bigl( \Delta_{\cal F}(\gamma)\Bigr) 
= \Delta_{\cal M} ({\rm Cor}_{\cal M})
\end{equation}
This allows us  to calculate the coproduct (\ref{8.23.02.8})
combinatorially as the left hand side in (\ref{8.29.02.1}), providing an  
answer to the question (\ref{8.23.02.1}).

{\bf Remarks}. 1. If the motivic correlator is a pure motive,  question 
(\ref{8.23.02.1}) is trivial: the coproduct is zero. 

2. As we will see below, the correspondence 
principle imposes very strong constraints on the motivic correlators. 
It is not clear to me whether we can expect it for 
any algebraic-geometric Feynman integral. 
See  s. 7.6.

{\bf 3. An example}. The story described in  Chapter 4 should serve as an 
example of such a situation. 
In this case the motivic iterated integrals 
\begin{equation} \label{8.25.01.1}
{\rm I}^{\cal M}(s_0; s_1, ..., s_m; s_{m+1}), \qquad s_i \in S
\end{equation}  
should be seen as 
the motivic correlators. 
So far there is no Feynman integral known providing these correlators, so we 
work directly  with Feynman diagrams, given by trivalent plane rooted trees.  
Thanks to 
Theorem \ref{7.23.02.2} there is a Hopf algebra map 
$$
t: \widetilde {\cal I}_{\bullet}(S) \hookrightarrow {\cal T}_{\bullet}(S)
$$ 
where  $t$ is as in (\ref{8.22.01.1}). The Hopf subalgebra 
$$
t(\widetilde {\cal I}_{\bullet}(S)) \subset {\cal T}_{\bullet}(S)  
$$   %(rather then ${\cal T}_{\bullet}(S)$ itself) 
should be considered as the combinatorially defined 
Hopf algebra ${\cal H}_{\cal F}$ responsible for  
the motivic correlators (\ref{8.25.01.1}). 
Then the map  $c_{\cal M}$ is given by 
$$
\sum\mbox{ 
plane rooted trivalent trees decorated by 
$\{s_0, s_1, ..., s_m, s_{m+1}\}$}\lms 
$$
$$
{\rm I}^{\cal M}(s_0; s_1, ..., s_m; s_{m+1})
$$

Since  ${\cal I}_{\bullet}(S)$ is a quotient of 
$\widetilde {\cal I}_{\bullet}(S)$, we arrive at the diagram 
$$
{\cal O}(G_{\cal M}(S))\twoheadleftarrow {\cal I}_{\bullet}(S) \twoheadleftarrow 
\widetilde 
{\cal I}_{\bullet}(S) =:{\cal H}_{\cal F}
$$
We define the  affine group scheme  $G_{\cal F}$ 
as the spectrum of ${\cal H}_{\cal F}$. Passing to the corresponding group schemes we get
$$
G_{\cal M}(S) \stackrel{(\ref{8.3.02.q12})}{\hookrightarrow} G(S) 
\stackrel{\rm def}{\hookrightarrow} {\rm Spec}(\widetilde 
{\cal I}_{\bullet}(S)) =:G_{\cal F}
%\stackrel{t^*}{\twoheadleftarrow} {\rm Spec}({\cal T}_{\bullet}(S)). 
$$

Observe that the coproduct in ${\cal H}_{\cal F}:= t(\widetilde {\cal I}_{\bullet}(S))$ 
is inherited from a bigger Hopf algebra 
 ${\cal T}_{\bullet}(S)$.

{\bf 4. The motivic Galois group of a Feynman integral and 
the correspondence principle}. 
Let us reformulate the correspondence principle as a relationship between 
the groups $G_{\cal F}$  and the motivic 
Galois groups. We start with a somewhat more precise formulation of the 
correspondence principle:

i) A Feynman integral ${\cal F}$ can be understood as an (infinite)
 collection of correlators.  
Correlators of Feynman integrals of algebraic--geometric type can be lifted to their motivic avatars. 
The commutative algebra they generate {\it should} form a  Hopf subalgebra 
${\cal H}_{\rm Mot}({\cal F})$  of the motivic Hopf algebra:  
\begin{equation} \label{8.24.02.10}
 {\cal H}_{\rm Mot}({\cal F}) \hookrightarrow {\cal H}_{\rm Mot}
\end{equation}
 We call it the {\it motivic Hopf algebra of the Feynman integral ${\cal F}$}.

ii) A Feynman integral ${\cal F}$ should determine 
its combinatorially defined Hopf algebra ${\cal H}_{\cal F}$. 
%We assume 
%(but do not insist) that 
%it is a Hopf subalgebra of certain more universal (renormalization) Hopf algebra
%${\cal H}_{R}$:
%\begin{equation} \label{8.24.02.11}
%{\cal H}_{R}({\cal F}) \hookrightarrow {\cal H}_{\rm R}
%\end{equation}

iii) These two Hopf algebras are related by a canonical 
homomorphism 
of Hopf algebras, the {\it motivic correlator homomorphism} 
$$
c_{\cal M}: {\cal H}_{\cal F} \twoheadrightarrow {\cal H}_{\rm Mot}({\cal F})
$$
provided by a tensor functor $F: T_{\cal F} \lra {\cal P}_{\cal M}$. 
The map $c_{\cal M}$ is surjective by the very definition of ${\cal H}_{\rm Mot}({\cal F})$.

 Just like in s. 7.3, the coproduct in ${\cal H}_{\cal F}$ 
is probably inherited from a bigger object: although only the sum 
$\sum[\Gamma]$ belongs to ${\cal H}_{\cal F}$, we should be able to make sense of 
$\Delta_{\cal F}[\Gamma]$ for every summand $[\Gamma]$.  

We define the {\it motivic Galois group $G_{\rm Mot}({\cal F})$  
of a Feynman integral ${\cal F}$} as the spectrum of the commutative algebra 
${\cal H}_{\rm Mot}({\cal F})$.  
Then $c_{\cal M}$ provides an injective 
homomorphism
$$
c_{\cal M}^*: G_{\rm Mot}({\cal F}) \hookrightarrow G_{\cal F}
$$

It seems that 
the Feynman integral itself should be upgraded to 
an infinite dimensional mixed motive. Its matrix coefficients 
corresponding to different framings should 
provide us all its correlators.

{\bf 5. The correspondence principle almost determines motivic correlators}. 
To simplify the exposition we will demonstrate this 
in the particularly interesting mixed Tate case. 
Similar arguments work in  general. 

Let $F$ be a number field. Denote by ${\cal H}_{\rm Mot}^T(F):= {\cal A}_{\bullet}(F) $ the fundamental Hopf algebra of the 
category of mixed Tate motives over $F$. It is a commutative, graded Hopf algebra. 

  If all motivic correlators 
of a Feynman integral ${\cal F}$ 
were mixed Tate motives,  
then 
$G_{\rm Mot}({\cal F})$ would be  a  semidirect product of ${\Bbb G}_m$ and 
${\rm Spec}({\cal H}_{\rm Mot}({\cal F}))$, with the action of ${\Bbb G}_m$ 
provided by the grading. 
In this case ${\cal H}_{\cal F}$ 
should also be a graded Hopf algebra, and 
the semidirect product of ${\Bbb G}_m$ and its spectrum 
provides us the  group 
$G_{\cal F}$ associated to 
 a Feynman integral.

Let ${\rm Cor}_{\cal M} = c_{\cal M}(\gamma)$, where $\gamma$ as in (\ref{8.24.02.1}). 
Suppose that ${\rm Cor}_{\cal M}$ is a mixed motive framed by $\Q(0)$ and $\Q(n)$ 
 defined over $F$. 
The calculation of  $\Delta_{{\cal H}_{\cal F}}(\gamma)$ is a combinatorial 
problem. So we may assume $\Delta_{{\cal H}_{\cal F}}(\gamma)$ is given to us. 
Let us stress that we do not assume that ${\rm Cor}_{\cal M}$ is a mixed Tate motive. 
 However let us assume that 
$$
(c_{\cal M} \otimes c_{\cal M})(\Delta'_{{\cal H}_{\cal F}}(\gamma)) \in 
{\cal H}_{\rm Mot}^T(F) \otimes {\cal H}_{\rm Mot}^T(F)
$$ 
is mixed Tate. Then thanks to (\ref{8.29.02.1}) 
\begin{equation} \label{9.4.02.1}
\Delta'_{\cal M}({\rm Cor}_{\cal M}) \in {\cal H}_{\rm Mot}^T(F)
 \otimes {\cal H}_{\rm Mot}^T(F)
\end{equation}

 Since we used the restricted coproduct in (\ref{9.4.02.1}), all factors 
there are of smaller weight, and are mixed Tate by assumption. 
Therefore the following lemma 
implies that 
if the framing 
on ${\rm Cor}_{\cal M}$ is Tate then ${\rm Cor}_{\cal M}$ 
is (equivalent to) a framed mixed Tate motive over $F$! 

\begin{lemma} Let $X\in {\cal H}_{\rm Mot}(F)$ 
be the equivalence class of a mixed motive over $F$ 
framed by $\Q(0)$ and $\Q(n)$.  Then  
$$
\Delta_{\cal M}' X 
\in {\cal H}^T_{\rm Mot}(F)
\otimes {\cal H}^T_{\rm Mot}(F)
$$
 implies $X \in {\cal H}^T_{\rm Mot}(F)$. 
\end{lemma}

{\bf Proof}. Observe that $\Delta_{\cal M}' (\Delta_{\cal M}'X) = 0$. Since 
${\rm Ext}^2_{\rm Mot/F}(\Q(0), \Q(n)) =0$ it follows that there exists an element 
$Y \in {\cal H}^T_{\rm Mot}(F)$ such that $\Delta_{\cal M}' Y = \Delta_{\cal M}'X$. 
Therefore the  lemma follows from the following  easy but basic  fact: 
\begin{equation} \label{9.9.02.1}
X, Y \in {\cal H}_{\rm Mot}(F), \quad 
\Delta_{\cal M}' X = \Delta_{\cal M}' Y \quad => \quad 
X - Y \in {\rm Ext}_{\rm Mot/F}^1(\Q(0), \Q(n))
\end{equation}

{\bf Remark}. There is a similar lemma where the category of Tate motives 
is replaced by any 
pure tensor category of motives. However its proof uses  a different idea. 

 If  ${\rm Cor}_{\cal M}$ is defined over a number field $F$ then thanks to (\ref{9.9.02.1}) 
it is determined 
by its restricted coproduct $\Delta_{\cal M}'({\rm Cor}_{\cal M})$ up to an element 
of 
\begin{equation} \label{9.9.02.2}
{\rm Ext}_{\rm Mot/F}^1(\Q(0), \Q(n)) = K_{2n-1}(F)\otimes \Q
\end{equation}
Therefore if 
a mixed motive representing ${\rm Cor}_{\cal M}$ is defined over $F$ 
(which is usually very easy to check), and the framing 
on it is Tate, 
 the knowledge of (\ref{9.4.02.1}) 
allows us to determine the equivalence class of 
${\rm Cor}_{\cal M}$ up to an element of (\ref{9.9.02.2}).  

For example let $F = \Q$. Then  (\ref{9.9.02.2}) is 
spanned by $\zeta^{\cal M}(n)$, see also (\ref{SEE}). 
It is zero 
for even $n$ and  one-dimensional for odd $n>1$. So in this case 
${\rm Cor}_{\cal M}$ is determined up to a rational multiple 
of $\zeta^{\cal M}(n)$. In particular it is  determined uniquely  if  $n$ is even. 

Recall the inductive definition of the Hopf algebra 
${\cal H}^T_{\rm Mot}(\Z) = {\cal A}_{\bullet}(\Z)$ 
given in Lemma \ref{8-24.2/99} and (\ref{5-28.00.1}) in the Appendix. 
Then induction on the weight shows that if we have 
the correspondence principle, 
and if the weight $2$ correlators are zero, 
then, thanks to (\ref{5-28.00.1}), in the situation above all motivic correlators must belong to 
 ${\cal H}^T_{\rm Mot}(\Z)$, i.e. must be framed mixed Tate motives over $\Z$. 
This plus Conjecture  17b) from [G6], which 
says that all framed mixed Tate motives over $\Z$ 
are multiple $\zeta$--motives,  imply that ${\rm Cor}_{\cal M}$
 is a multiple $\zeta$--motive. 

More generally, let ${\cal O}_S$ be the ring of $S$--integers in a number field 
$F$. Then if  all weight two motivic 
correlators come from  ${\cal O}^*_S$, 
 then in the situation above all motivic correlators must  
come from ${\cal H}^T_{\rm Mot}({\cal O}_S)$, i.e. are framed mixed Tate motives over $F$ 
with possible ramification at $S$.

{\bf 6. An example: $\zeta^{\cal M}(3,5)$}. Physicists computed many  correlators.  
D. Kreimer and D. Broadhurst discovered that, in small weights, 
these correlators are often  expressed by  the multiple $\zeta$--values ([BGK]).  
According to the example in Section 6.7, the simplest multiple 
$\zeta$--value which should not be expressed 
by the classical $\zeta$--values appears in  weight $16$ 
(the weights in the Tate case are even numbers, so it is customary to 
divide them by two; then the weight is $8$).  
For example one can take  $\zeta(3,5)$. Its motivic avatar has 
been discussed in the example
mentioned above. 
According to  [BGK] $\zeta(3,5)$ appears as the correlator corresponding to the 
following remarkable Feynman diagram: 

\begin{center}
\hspace{4.0cm}
\epsffile{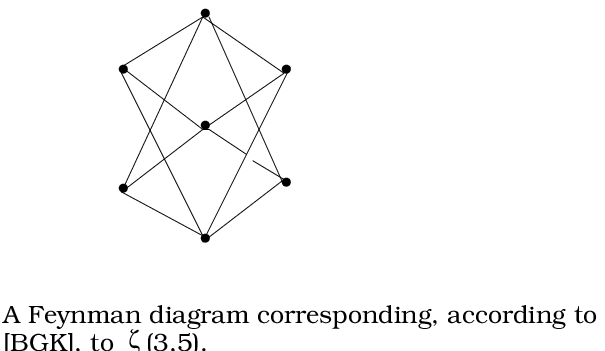}
\end{center}
This graph has some combinatorial properties 
which can not be found in any simpler graph. For instance it is the 
simplest graph which does not have Hamiltonian circles, i.e. cycles 
without selfintersections going through all vertices. 
So in this respect it is like $\zeta(3,5)$, which is the 
simplest non--classical $\zeta$--value. 
It would be very interesting to recover formula (\ref{9.4.02.4}) 
for the coproduct  from the combinatorics of this graph. 

According to  s. 7.5, 
the correspondence principle explains why the  multiple $\zeta$--values 
appear in calculations of correlators.
% and predicts the values of 
% such correlators modulo the ideal generated by $\pi^2$ in the multiple $\zeta$--algebra 
%  up to a rational multiple of $\zeta(2n+1)$.   
On the other hand,  results of the 
paper [BB]  suggest, although do not prove,   that 
correlators of certain physically interesting Feynman integral 
might  be any periods, and not only very specific multiple $\zeta$--values. 
It is not clear whether this, if true,    contradicts the correspondence principle:
the non--mixed Tate periods of smallest weight can appear because 
of non--Tate framing corresponding to the correlator, and then 
spread out into higher weights.  The subject needs further investigation.

%However it  suggests 
%that we can not expect that the correspondence principle 
%works for an arbitrary  Feynman integral.  

We assumed before that the correlators are given by convergent integrals, 
and thus are well defined numbers/motives. 
If these integrals are divergent, they  still have well-defined 
regularizations, given by convergent finite dimensional 
Feynman integrals, which deliver the interesting numbers/functions. 
The renormalization group enters into the game, acting 
on possible regularizations. The Hopf algebra ${\cal H}_{\cal F}$ 
should  extend the renormalization Hopf algebra of the type studied in [CK], 
since the latter concentrates on the divergent correlators, but 
conceptually there should be no difference between ${\cal H}_{\cal F}$  
and the renormalization algebra.

\section{ Appendix. Mixed Tate categories and framed objects}

We address   
 the Tannakian formalism for mixed Tate categories in the language 
of framed objects.
For a more  general set-up see chapter 3 of [G7]. 
 Let ${\rm Spec}{\cal O}_{F, S}$ be the scheme obtained by deleting 
an arbitrary set $S$  of closed points from the spectrum of the 
ring of integers of a number field $F$. 
The abelian category of mixed Tate motives 
over ${\rm Spec}{\cal O}_{F, S}$ was defined in [DG]. 
 We give a simple description 
of its fundamental Hopf algebra.

{\bf 1. A review of the Tannakian formalism for mixed Tate categories}. 
Let $K$ be a characteristic zero field. 
Let ${\cal M}$ be a Tannakian  $K$-category with an invertible object $K(1)$. 
So 
in particular ${\cal M}$ is an abelian tensor  category. 
Set $K(n):= K(1)^{\otimes n}$, where $K(-1)$ is dual to $K(1)$. 
Recall ([BD]) that a pair $({\cal M}, K(1))$ is called a {\it mixed 
Tate category} 
if the objects $K(n)$
are mutually nonisomorphic, any simple object is isomorphic to one of them, and
$Ext_{{\cal M}}^1(K(0),K(n)) = 0$ if $n \leq 0$. A pure functor 
between two mixed Tate categories $({\cal M}_1, K(1)_1)$ and $({\cal M}_2, 
K(1)_2)$ is a 
tensor functor $\varphi: {\cal M}_1 \lra {\cal M}_2$  equipped 
with an isomorphism $\varphi(K(1)_1) = K(1)_2$. 

It is easy to show that any object $M$ of a mixed Tate category 
has a canonical 
weight filtration $W_{\bullet}M$ indexed by $2\Z$ 
such that $gr^W_{2n}M$ is a direct sum of copies of $K(-n)$, and 
morphisms are strictly compatible with the weight filtration. The functor
$$
\omega = \omega_{{\cal M}}:  {\cal M} \longrightarrow Vect_{\bullet}, 
\quad M \longmapsto 
\oplus_{n} Hom_{{\cal M}}(K(-n), gr^W_{2n}M)
$$
  to the category of graded $K$-vector spaces is a fiber functor. Let 
$\widetilde \omega$ be the fiber functor to the category of finite 
dimensional 
$K$-vector spaces obtained from $\omega$ by forgetting the grading. Let $$
L_{\bullet}({\cal M}):= Der ( \widetilde \omega) : = 
\{F \in End \widetilde\omega | F_{X \otimes Y} = F_{X} 
\otimes id_{Y} +  id_{X} \otimes F_{Y}\}
$$  
be the space of its derivations. 
It is a pro-Lie algebra over $K$, called 
 the {\it fundamental Lie algebra} of the mixed Tate category ${\cal M}$. 
It has a natural grading by the integers $n<0$ 
provided by the original grading 
on $\omega({\cal M})$. 

The automorphisms of the fiber functor $\widetilde \omega$ 
respecting the tensor structure provide 
a pro-algebraic group scheme over $K$, denoted 
$ Aut^{\otimes} \widetilde \omega$. It is a semidirect product of 
the multiplicative group scheme ${\Bbb G}_m$ 
and a pro-unipotent group scheme 
$U({\cal M})$. The pro-Lie algebra $L_{\bullet}({\cal M})$ is the 
Lie algebra of  $U({\cal M})$. The action of ${\Bbb G}_m$ provides 
a grading on $L_{\bullet}({\cal M})$. 

According to the Tannakian formalism the functor $\widetilde \omega$ 
provides an equivalence 
between the category $\cal M$ and 
 the category of finite dimensional modules over the
pro-group scheme $ Aut^{\otimes} \widetilde \omega$.  
This category is naturally equivalent 
to the category of graded finite dimensional modules over the 
group scheme $U({\cal M})$.
 Since $U({\cal M})$ is  pro-unipotent,  the last category is 
equivalent  to the category of 
graded finite dimensional modules over the
graded pro-Lie algebra $L_{\bullet}({\cal M})$.  

Pure functors  $({\cal M}_1, K(1)_1)  \lra ({\cal M}_2, K(1)_2)$ between the 
mixed Tate categories 
 are in bijective correspondence with the graded Lie algebra 
morphisms  $L_{\bullet}({\cal M}_2) \lra L_{\bullet}({\cal M}_1)$. 

 Let  ${\cal U}_{\bullet}({\cal M}):= End(\widetilde \omega)$ be the 
space of all endomorphisms of the 
fiber functor $\widetilde \omega$. It is a graded 
Hopf algebra isomorphic to the universal enveloping algebra 
of the Lie algebra $L_{\bullet}({\cal M})$. 

Recall that a {\it Lie coalgebra} is a vector 
space ${\cal D}$ equipped with a linear map 
$
\delta:  {\cal D} \lra \Lambda^2{\cal D} 
$
such that the composition
$
{\cal D} \stackrel{\delta}{\lra}\Lambda^2{\cal D} \stackrel{\delta \otimes 1 - 1 \otimes 
\delta}{\lra}
\Lambda^3{\cal D} 
$ 
is zero. If ${\cal D} $ is finite dimensional then it is a Lie coalgebra if 
and only if its dual is a Lie algebra.

 Let $G$ be a 
unipotent algebraic group over $\Q$. Then the ring of regular functions 
$\Q[G]$ is a Hopf algebra with the coproduct induced by the multiplication in $G$.
The (continuous) dual of its completion at the group unit $e$ is isomorphic 
to the universal enveloping algebra
of the Lie algebra ${\rm Lie}(G)$ of $G$. Let $\Q[G]_0$ be the ideal of functions
equal to
zero at $e$. Then the coproduct induces on $\Q[G]_0/\Q[G]_0 \cdot \Q[G]_0$ the 
structure of a Lie coalgebra  dual to ${\rm Lie}(G)$.

Recall the  duality $V \lms V^{\vee}$ between the Ind- and pro-objects 
in the category of finite dimensional $K$-vector spaces. The graded dual Hopf algebra ${\cal U}_{\bullet}({\cal M})^{\vee}:= 
\oplus_{k \geq 0} {\cal U}_{-k}({\cal M})^{\vee}$ can be identified
with the Hopf algebra    
of regular functions on the pro-group scheme $U({\cal M})$. 
Therefore its  quotient by the square of the augmentation ideal 
$$
{\cal L}_{\bullet}({\cal M}):= \frac{{\cal U}_{>0}({\cal M})^{\vee}}
{{\cal U}_{>0}({\cal M})^{\vee} \cdot {\cal U}_{>0}({\cal M})^{\vee}}
$$
is a Lie coalgebra. 
The cobracket $\delta$ on ${\cal L}_{\bullet}({\cal M})$ 
is induced by the restricted coproduct 
\begin{equation} \label{8-23.2/99}
 \Delta'(X):= \quad \Delta (X)-  (X \otimes 1 + 1 \otimes X)
\end{equation}
on ${\cal U}_{\bullet}({\cal M})^{\vee}$. 
The  graded dual of the Lie coalgebra ${\cal L}_{\bullet}({\cal M})$ 
is identified with the fundamental Lie algebra 
$L_{\bullet}({\cal M})$. Below we  recall a more efficient 
  way to think about  it. 

{\bf 2. The Hopf algebra of framed objects in a mixed Tate category 
(cf. [BGSV], [G7])}. Recall that 
a Hopf algebra over a field $K$ is an associative $K$-algebra $A$ with 
a product $\mu: A \otimes A \to A$ and a unit $i: A \to K$,  
equipped in addition with 
a comultiplication $\Delta: A \to A \otimes  A$, 
a counit $\varepsilon: K \to A$, and an antipode $S: A \to A$.
They must obey the following properties:

1. The maps $\Delta$ and $\varepsilon$ define 
a structure of a coassociative coalgebra on $A$.

2. The maps $\Delta: A \to A \otimes A$ and $\varepsilon: A \to K$ 
are homomorphisms of algebras.

3. The map $S$ is a linear isomorphism satisfying the relations
\begin{equation} \label{4.16.01.1}
\mu \circ (S \otimes {\rm id}) \circ \Delta  = \mu \circ ({\rm id} \otimes S) 
\circ \Delta  = i \circ \varepsilon
\end{equation}
Removing the antipode $S$ and condition (\ref{4.16.01.1}) from the  list 
of axioms we get a definition of a  bialgebra. 

\begin{lemma} \label{4.15.01.1}
Let $A_{\bullet}$ be a commutative bialgebra 
graded by the integers $n \geq 0$ such that 
$A_0 =k$. Then there exists a unique 
antipode map $S$ on $A_{\bullet}$. 
\end{lemma}
 
{\bf Proof}. The condition (\ref{4.16.01.1}) 
determines uniquely the restriction of the map 
$S$ to $A_n$ by induction. One has $S|A_0 = {\rm id}$. Then, for example 
$S|A_1 = -{\rm id}$ and  $S|A_2 = -{\rm id} + \mu \Delta_{1,1}$ 
where $\Delta_{1,1}$ is the 
$A_1 \otimes A_1$ component of the coproduct. The lemma is proved.

Let $n \geq 0$. Say
that   $M$ is an {\it $n$-framed  object} of ${\cal M}$, denoted 
$(M,v_0,f_n)$,   
if it is supplied with  
non-zero morphisms $v_{0} : K(0) \longrightarrow 
gr^W_{0}M$ and   $f_{n}: gr^W_{-2n} M \longrightarrow K(n)$.  
Consider the  coarsest equivalence
relation on the set of all $n$-framed  objects for which $M_1 \sim M_2$
if there is a map $ M_1 \to  M_2$ respecting the frames.  For example 
replacing $M$ by 
$W_0M/W_{-2p-2}M$ we see that any $n$-framed
 object is equivalent to a one $ M$ with $W_{-2-2n} M = 0$, $W_{0} M= M$.
Let ${\cal A}_n({\cal M})$ be the set of equivalence classes.
It  has the structure of an abelian group  
with the composition law defined as follows:
$$
[ M, v_{0},f_{n}] + [ M', v_{0}',f_{n}'] \quad := \quad [M \oplus M', (v_{0}, 
v_{0}'), f_{n} + f_{n}' ]
$$
It is straightforward to check that the composition law is well 
defined on equivalence classes of framed objects. Indeed, if 
$\varphi: \widetilde M \lra M$ 
is a morphism providing an equivalence of the framed objects 
$[ \widetilde M, \widetilde v_{0}, \widetilde f_{n}] \sim 
[ M, v_{0}, f_{n}]$ then 
$ \varphi \oplus {\rm id}: \widetilde M \oplus M' \lra M \oplus M'$ 
also provides an equivalence. 
The neutral element is $K(0) \oplus K(n)$ with the obvious frame. 
The inversion is given by 
$$
-[M, v_{0}, f_{n}]:= \quad [M, -v_{0}, f_{n}] \quad = \quad [M, v_{0}, -f_{n}]
$$
See ch. 2 of [G8] for a proof of these two facts. 
The composition $f_0\circ v_0: K(0) \to K(0)$ provides an 
isomorphism ${\cal A}_0({\cal M}) = K$. 
It follows from the very definitions that there is a canonical isomorphism
$$
{\cal A}_1({\cal M}) = {\rm Ext}_{\cal M}^1(K(0), K(1))
$$
The tensor product   induces the commutative and associative multiplication
$$
\mu: {\cal  A}_{k}({\cal M}) \otimes {\cal  A}_{\ell}({\cal M})\to {\cal A}_{k+ \ell}({\cal M})
$$
One verifies that it is well defined on equivalence classes using an argument similar 
to the one used for checking the 
additive structure on ${\cal  A}_{k}({\cal M})$. 
The unit is given by $1 \in K = {\cal A}_0({\cal M})$. 
The counit is the projection of ${\cal  A}_{\bullet}({\cal M})$ 
onto its $0$-th component.

Let us define the comultiplication
$$
\Delta = \bigoplus_{0 \leq p \leq n} \Delta_{p, n-p}: \quad 
{\cal  A}_{n}({\cal M})
\to \bigoplus_{0 \leq p \leq n } 
{\cal  A}_{p}({\cal M})\otimes
{\cal  A}_{n-p}({\cal M})
$$
 Choose a basis $\{b_i\}$, where $1 \leq i \leq m$,  of ${\rm Hom}_{\cal M} (K(p),gr_{-2p}^W  M)$   and  the dual basis $\{b'_i\}$ of 
${\rm Hom}_{\cal M} (gr_{-2p}^W M, K(p))$.  
Then
$$
\Delta_{p, n-p } [M,v_{0},f_n]:=  \quad \sum_{i=1}^m [M,v_{0},b'_i] 
\otimes [M, b_i, f_n](-p)
$$
In particular $\Delta_{0,n} = id \otimes 1$ and $\Delta_{n,0} = 1 \otimes id$.  
Set ${\cal A}_{\bullet}({\cal M}):=\oplus {\cal A}_n({\cal M})$.

\begin{theorem} \label{8-23.3/99}
a) ${\cal  A}_{\bullet}({\cal M})$ has the structure of a 
graded Hopf algebra over $K$ with the commutative multiplication
$\mu$ and the comultiplication $\Delta$.

 b) The Hopf algebra ${\cal  A}_{\bullet}({\cal M})$ is 
canonically isomorphic to the Hopf algebra  
${\cal U}_{\bullet}({\cal M})^{\vee}$.
\end{theorem}

\begin{definition} \label{ui}
The Hopf algebra 
${\cal  A}_{\bullet}({\cal M})$ is called the {\em fundamental Hopf algebra}
 of a mixed Tate category ${\cal M}$. 
\end{definition}

{\bf Proof}. a) Let us show that the coproduct is well defined on equivalence 
classes of framed objects. It is sufficient to prove this for equivalences 
given by 
injective and surjective morphisms in ${\cal M}$. Indeed, if 
$\varphi: M \lra M'$ respects the frames in $M$ and $M'$ then the projection 
$M \lra M/{\rm Ker} (\varphi)$ and injection $M/{\rm Ker} (\varphi) 
\hookrightarrow  M'$ 
also respect the frames. Let us suppose that 
$$
\varphi : [M, v_0, f_n] \hookrightarrow [M', v'_0, f'_n]
$$
is an equivalence. Choose a basis $\{b_i\}$ of 
${\rm Hom}_{\cal M} (K(p),gr_{-2p}^W  M')$ 
such that the first $s$ basis 
vectors of $\{b_i\}$ is a basis in 
${\rm Hom}_{\cal M} (K(p),gr_{-2p}^W  M)$. Denote by $b_i'$ the dual basis. 
Then 
$[M', v_0, b_{i}']=0$ for $i>s$. Indeed, we may assume that 
$W_{-2p-2}M'=0$, and also ${\rm gr}^W_0M'= K(0)$. 
Then there is a natural injective morphism 
$K(p) \oplus M  \hookrightarrow M'$ respecting the frames, 
and a projection 
$
K(p) \oplus M \to K(p) \oplus M/W_{-2}M = 
K(p) \oplus K(0)
$.  The statement is proved. 
The arguments in the case of the projection are similar 
(and can be obtained by dualization).

It is straightforward to show 
that $\Delta$ is a homomorphism of algebras, and that 
$\Delta$ is coassociative. 
The part a) of the theorem is proved.

b)  We follow the proof of theorem 3.3 in [G7], making some necessary 
corrections. A canonical isomorphism $\varphi: {\cal  A}_{\bullet}({\cal M}) \to 
{\cal U}_{\bullet}({\cal M})^{\vee}$
is constructed as follows.
Let $F \in End(\widetilde \omega)_n$ and 
$[M,v_{0}, f_n]\in {\cal  A}_{n}({\cal M})$. 
Denote by $F_M$ the endomorphism of $\omega(M)$ provided by $F$. Then 
$$
<\varphi([M,v_{0}, f_n]) , F>:= \quad <f_n, F_M(v_{0})>
$$
It is obviously well defined on equivalence classes of framed objects. 
$\varphi$ is evidently a morphism of graded Hopf algebras. 
Let us show that $\varphi$ is surjective. 
Recall that for a commutative, non-negatively graded 
Hopf algebra ${\cal A}_{\bullet}$ the space of primitives 
$
{\rm CoLie} ({\cal A}_{\bullet}):= 
{\cal A}_{>0}/{\cal A}_{>0}^2
$ 
has a Lie coalgebra structure, and the dual to the 
universal enveloping algebra of the dual Lie algebra is 
canonically isomorphic to  ${\cal A}_{\bullet}$. 
Since $\varphi$ is a map of Hopf algebras it provides a morphism of the 
Lie coalgebras 
$$
\overline \varphi: {\rm CoLie} ({\cal A}_{\bullet}({\cal M})) \lra 
{\rm CoLie} ({\cal U}_{\bullet}({\cal M})^{\vee})
$$
Let us show that this map is surjective. 
Let $f$ be a functional on ${\rm CoLie} ({\cal U}_{\bullet}({\cal M})^{\vee})$ 
which is zero 
on the graded components of degree $\not = n$. Consider $f$ as a functional 
on ${\cal U}_{\bullet}({\cal M})^{\vee}$, and denote by ${\rm Ker}(f)$ its kernel. 
Then ${\cal U}_{\bullet}({\cal M})^{\vee}/{\rm Ker}(f) = :K_{(-n)}$ is a one-dimensional 
space sitting in degree $-n$. The graded $K$-vector space 
$K_{(0)} \oplus {\cal U}_{\bullet}({\cal M})^{\vee}/{\rm Ker}(f) $ has a 
${\cal U}_{\bullet}({\cal M})^{\vee}$-module structure: an element 
$X \in {\cal U}_{>0}({\cal M})^{\vee}$ sends 
$1 \in K_{(0)} \lms f(X) \in K_{(-n)}$ and annihilates $K_{(-n)}$. Since by the 
definition of $f$ the square of the augmentation ideal acts by zero, 
we get a well defined action of ${\cal U}_{\bullet}({\cal M})^{\vee}$. 
Therefore the map $\overline \varphi$ is surjective. 
Dualizing $\overline \varphi$ we get an injective map of 
Lie algebras, and hence an injective map of the corresponding 
universal enveloping algebras. Dualizing the latter map
 we prove the surjectivity of $\varphi$.

Now let us show that $\varphi$ is injective. Using the Tannaka theory we may assume that 
${\cal M}$ is the category of finite dimensional representations of 
the Hopf algebra ${\cal U}_{\bullet}({\cal M})$. Suppose that 
$\overline \varphi\Bigl([M, v_0, f_n] \Bigr) =0$. 
Consider the cyclic submodule 
${\cal U}_{\bullet}({\cal M}) \cdot K_{(0)}$. It has no non-zero 
components in degree $-n$, since otherwise $\varphi[M, v_0, f_n] \not = 0$. 
Thus there are maps 
$$
K_{(0)} \oplus K_{(-n)} \longleftarrow {\cal U}_{\bullet}({\cal M}) \cdot K_{(0)} \oplus K_{(-n)} \lra M / W_{<-n}M
$$
respecting the frames. Thus $[M, v_0, f_n]=0$. 
The part b), and hence the theorem, are proved.

Let $M^*$ be the object dual to $M$. 
Under the  equivalence between the  tensor category ${\cal M}$ and the category of 
graded finite dimensional comodules 
over the Hopf algebra 
${\cal  A}_{\bullet}({\cal M})$ an object $M$ of ${\cal M}$ corresponds to 
the graded comodule 
$\omega (M)$ with 
${\cal  A}_{\bullet}({\cal M})$-coaction $\omega (M )
\otimes \omega (M^{\ast}) \longrightarrow {\cal  A}_{\bullet}({\cal M})$ given by the formula 
\begin{equation} \label{8-24.5/99}
x_{m} \otimes y_{n} \longrightarrow \quad \mbox{the class of $M$ 
framed by $x_{m}, y_{n}$}
\end{equation}
 We call the right hand side the {\it matrix coefficient} of $M$ corresponding to $x_{m}, y_{n}$.

The restricted coproduct
$\Delta'$  provides the quotient 
$
{\cal A}_{\bullet}({\cal M})/({\cal A}_{>0}({\cal M}) )^2
$ 
with the structure of a graded Lie coalgebra  with cobracket 
$\delta$. It is canonically isomorphic to ${\cal L}_{\bullet}({\cal M})$.

\begin{lemma} \label{2.3.01.2}
Each equivalence class of $n$-framed objects contains a unique 
minimal representative, which appears as a subquotient in any 
$n$-framed object from the given equivalence class. 
\end{lemma}

{\bf Proof}. Suppose we have morphisms $M_1 \stackrel{f}{\longleftarrow} N 
\stackrel{g}{\longrightarrow} M_2$ between the $n$-framed objects respecting the frames. For any $n$-framed object $X$ we can assume that ${\rm gr}^W_{2m}X =0$ unless $ -n \leq m \leq 0$ as well as 
${\rm gr}^W_{-2n}X = \Q(n)$; ${\rm gr}^W_{0}X = \Q(0)$. Taking the subquotient ${\rm Im}(f)$ 
of $M_1$ we may suppose that $f$ is surjective. Then ${\rm gr}^W_{0}{\rm Ker}(f) = 
{\rm gr}^W_{-2n}{\rm Ker}(f) = 0$. Therefore $g({\rm Ker}(f))$ has the same property, and hence $M_2$ is equivalent to $N_2':= M_2/{\rm Ker}(f)$. The lemma follows from these remarks.

{\bf 3. Examples of mixed Tate categories}.  
Here is a general method of getting 
mixed Tate categories (cf. [BD]). Let 
${\cal C}$ be any Tannakian $K$-category and $K(1)$ be a 
rank one object of ${\cal C}$ such 
that the objects $K(i):= K(1)^{\otimes i}$, $i \in \Z$, 
are mutually nonisomorphic. 
An object $M$ of ${\cal C}$ is called a {\it mixed Tate object} if it admits 
a finite 
increasing filtration $W_{\bullet}$ indexed by $2\Z$, such that 
${\rm gr}^W_{2p}M$ is a 
direct sum of copies of $K(p)$. Denote by ${\cal T}{\cal C}$ the 
full subcategory of 
mixed Tate objects in ${\cal C}$. Then 
${\cal T}{\cal C}$ is a Tannakian subcategory in ${\cal C}$ and 
$({\cal T}{\cal C}, K(1))$ 
is a mixed Tate category. 
Further, if $({\cal C}_i, K(1)_i)$ are as above and $\varphi: {\cal C}_1 \lra {\cal C}_2$ is a 
tensor functor equipped with an isomorphism $\varphi(K(1)_1) = K(1)_2$, then 
$\varphi({\cal T}{\cal C}_1) \subset  {\cal T}{\cal C}_2$ and the restriction of $\varphi$
to ${\cal T}{\cal C}_1$ is a pure functor.

1. {\it The category of $\Q$-rational Hodge-Tate structures}. 
Applying the above construction 
to the category $MHS/\Q$ of mixed Hodge structures over $\Q$ we 
get the mixed Tate category 
${\cal H}_T$ of $ \Bbb Q$-rational Hodge-Tate structures. Namely, 
a $ \Bbb Q$-rational Hodge-Tate structure is a mixed Hodge structure 
over $\Bbb Q$ 
with $h^{p,q} =0$ if $p \not = q$. Equivalently, a Hodge-Tate structure 
is a $ \Bbb Q$-rational mixed Hodge structure with 
weight $-2k$ quotients isomorphic to 
a direct sum of copies of $\Q(k)$. 
Set ${\cal H}_{\bullet}:= {\cal A}_{\bullet}({\cal H}_{T})$. Then 
\begin{equation} \label{8-24/99}
{\cal H}_1 \stackrel{\sim}{=} Ext^1_{{\cal H}_{T}}(\Q(0),\Q(1)) 
 \stackrel{\sim}{=} 
 \frac{\C}{2 \pi i \Q} \stackrel{\sim}{=}  \C^*\otimes {\Q}
\end{equation}

{\bf Example}. The extension  corresponding 
via the isomorphism  (\ref{8-24/99}) to a given $z \in \C^*$ is provided by the 
$\Q$-rational Hodge-Tate structure 
$H(z)$, also denoted $\log^{\cal H} (z)$,  defined by the period matrix 
$$
 \left (\matrix{1&0 \cr \log( z)& 2 \pi i  \cr }\right ) 
$$ 
Namely, denote by  $H_{\C}$ the two dimensional $\C$-vector space with  basis 
$e_0, e_{-1}$. 
The $\Q$-vector space $H(z)_{\Q}$ is the subspace generated by 
the columns of the matrix, i.e. the vectors 
$e_0 + \log z \cdot e_{-1}$ and $2 \pi i e_{-1}$. 
The weight filtration on it is given by 
$$
W_0H(z)_{\Q} = H(z)_{\Q}, 
\quad W_{-1}H(z)_{\Q} = W_{-2}H(z)_{\Q} = 2 \pi i \cdot e_{-1}, \quad 
W_{-3}H(z)_{\Q} = 0
$$
The Hodge filtration 
is defined by 
$
F^{1}H_{\C} = 0, \quad F^0H_{\C} = <e_0>, 
\quad F^{-1}H_{\C} = H_{\C}
$.

It is easy to see that
$$
Ext^1_{{\cal H}_{T}}(\Q(0),\Q(n)) \stackrel{\sim}{=} \frac{\C}{(2 \pi i)^n \Q} \stackrel{\sim}{=}  \C^*(n-1)\otimes {\Q}
\quad \mbox{for $n > 0$}
$$
It was proved by Beilinson  that the higher Ext groups are zero. 
So the fundamental Lie algebra of the category of 
Hodge-Tate structures is isomorphic to a  free graded pro-Lie algebra over $\Q$ generated by the 
$\Q$-vector spaces $(\C^*(n-1)\otimes {\Q})^{\vee}$ sitting in degree $-n$ 
where $n \geq 1$.

{\it 2. The abelian category of mixed Tate
  motives over a number field}. 
Let  $F$ be a  number field. Then, as explained in [L1], ch. 5 of [G2], and [DG], 
there exists an 
  abelian category  ${\cal M}_T(F)$ 
of mixed Tate
  motives over  $F$ with  all the needed
  properties. In particular it is a Tate category and  
\begin{equation} \label{8-23/99}
Ext^1_{{\cal M}_T(F)}(\Q(0),\Q(n)) = K_{2n-1}(F) \otimes \Q, 
\end{equation}
and the higher Ext groups vanish. 
Therefore the fundamental Lie algebra ${L}_{\bullet}(F)$  of the 
category ${\cal M}_T(F)$, is isomorphic to 
a free graded pro-Lie algebra over $\Q$ generated by the duals to the 
finite dimensional $\Q$-vector spaces 
\begin{equation}\label{5-28.00}
K_{2n-1}(F) \otimes \Q \stackrel{\sim}{=} 
[\Z^{{\rm Hom}(F,\C)}\otimes\Q(n-1)]^+, 
\quad n >1
  \end{equation} 
(where $+$ means the invariants under complex conjugation) 
sitting in degree $-n$, 
 and by the infinite dimensional $\Q$-vector space $(F^* \otimes \Q)^{\vee}$ 
for $n=-1$. 
The isomorphism (\ref{5-28.00}) is given by the regulator map, and is due 
to Borel [Bo]. 

The category  ${\cal M}_T(F)$  is equipped with an array of realization functors  
(see [DG]).

3. {\it The mixed Tate category of lisse $\Q_l$-sheaves on a scheme $X$}. 
Let $X$ be a connected coherent scheme over $\Z_{(l)}$ such that ${\cal O}^*(X)$ does not 
contain  all roots of unity of order equal to a power of $l$. Denote by 
${\cal F}_{\Q_l}(X)$ the Tannakian category of 
lisse $\Q_l$-sheaves on $X$.   
There is  the Tate sheaf $\Q_l(1):= \Q_l(1)_X$. 
Since $\mu_{l^{\infty}} \not \subset {\cal O}^*(X)$ the Tate sheaves $\Q_l(m)_X$ 
 are mutually nonisomorphic. So 
the general 
construction above gives a mixed Tate category ${\cal T}{\cal F}_{\Q_l}(X)$ 
of lisse $\Q_l$-sheaves on  $X$. 
It was considered by Beilinson and Deligne in [BD],  
see an exposition 
in  s. 3.7 of [G1]. Let us stress that in the category ${\cal F}_{\Q_l}(X)$
one may not have 
${\rm Ext}^1(\Q_l(0), \Q_l(n)) =0$ for $n<0$. However in the category 
${\cal T}{\cal F}_{\Q_l}(X)$ we force this to be true. 

In particular, let $F$ be a field which does not
contain all $l^{\infty}$ roots of unity. For example, 
 $F$ can be a number field. Then we get a category 
${\cal T}{\cal F}_{\Q_l}(F)$ of $l$-adic 
mixed Tate Galois modules.  The underlying vector space of the
representation provides another fiber functor on this category. 
It can be related to the canonical fiber functor sa follows. 
Consider the
subcategory ${\cal T}{\cal F}_{\Q_l}(F, p)$ of the Galois modules
unramified at a prime $p$. Then making the choices i) and ii) as in Chapter
1.4, we get an isomorphism of the two fiber functors. 
It follows that for any finite dimensional mixed Tate Galois
representation $V$ unramified at $p$ the Lie algebra ${\cal G}_V$ of the image of ${\rm Gal}(\overline
\Q/F(\zeta_{l^{\infty}}))$ is isomorphic to the image of the
fundamental Lie algebra of this mixed Tate category acting on
$\omega(V)$. The weight filtration on $V$ induces a weight filtration
on the Lie algebra ${\cal G}_V$, and on the associated graded the
above isomorphism  is canonical.  
The isomorphism of fiber functors allows to consider the 
equivalence classes of the framed Galois representations as the
functions on the Galois group ${\rm Gal}(\overline
\Q/F(\zeta_{l^{\infty}}))$. This way we identify the ${\rm I}^{(l)}$
and ${\rm I}_{F_p}^{(l)}$-versions of the $l$-adic iterated integrals
discussed in Chapter 1.4. This identification preserves the
coproduct. So  the description of the
coproduct for ${\rm I}_{F_p}^{(l)}$ follows from the corresponding 
result for their ${\rm I}^{(l)}$-counterpart.  

{\it 4. The hypothetical abelian category of mixed Tate
  motives over an arbitrary field $F$}. It is supposed to be a mixed Tate category 
with the Ext groups given by Beilinson's formula
\begin{equation} \label{8-23/99g}
Ext^i_{{\cal M}_T(F)}(\Q(0),\Q(n)) = {\rm gr}^{\gamma}_nK_{2n-i}(F) \otimes \Q, 
\end{equation}
In the number field case it is known that the right hand side is trivial for $i>1$,
so we arrive at (\ref{8-23/99}). The corresponding motivic Lie algebra 
is not free in general. The formula (\ref{8-23/99g}) leads to the canonical 
isomorphisms
\begin{equation} \label{8-23/99gs}
{\cal A}_1({\cal M}_T(F)) = Ext^1_{{\cal M}_T(F)}(\Q(0),\Q(1)) = K_{1}(F) 
\otimes \Q = F^*\otimes \Q  
\end{equation}
We will often employ a notation ${\cal A}_{\bullet}(F)$ for ${\cal A}_{\bullet}({\cal M}_T(F))$, and call it the {\it motivic Hopf algebra of a field $F$}.   
The grading of the Hopf algebra ${\cal A}_{\bullet}(F)$ 
provides an action of the multiplicative group ${\Bbb G}_m$ on 
the proalgebraic group scheme given by the spectrum 
of ${\cal A}_{\bullet}(F)$ viewed as a  commutative algebra. 
The corresponding semidirect product 
of ${\Bbb G}_m$ and ${\rm Spec}({\cal A}_{\bullet}(F))$  is called the 
{\it motivic Tate Galois group of $F$}.  

{\bf 4. The fundamental Hopf algebra of the abelian category of 
mixed Tate motives over the ring of $S$--integers in 
a number field}. 
Let $M$ be an object  of a mixed Tate category ${\cal M}$.
We say that a framed object $(M, x_n, y_m)$ is a {\it matrix element of 
amplitude $n-m$} of $M$ if $x_n$ and $y_m$ are the framings at the weights 
$-2n$ and $-2m$ respectively. 
%$x_n \in Hom(\Q(n), {\rm gr}^W_{-2n}M)$ and 
%$y_n \in Hom({\rm gr}^W_{-2m}M, \Q(m))$. 
It is, by definition, an element of ${\cal A}_{n-m}({\cal M})$. 

Let ${\cal O}$ be the ring of integers in a number field $F$, 
 $S$ a set of prime ideals in ${\cal O}$, and 
${\cal O}_{F, S}$ the  localization
 of ${\cal O}$ at the primes belonging to  $S$. 
In [DG] an abelian 
 category ${\cal M}_T({\cal O}_{F, S})$ of mixed Tate motives 
over ${\cal O}_{F, S}$ has been defined as a subcategory of 
${\cal M}_T(F)$. Namely, an object $M$ of ${\cal M}_T(F)$ 
belongs to the subcategory ${\cal M}_T({\cal O}_{F, S})$ if and only if 
all amplitude one matrix elements of $M$ provide the extension classes 
with the invariants in ${\cal O}_{F, S}^*\otimes \Q 
\hookrightarrow F^*\otimes \Q$. 
It was proved in [DG] that it is a mixed Tate category with 
\begin{equation} \label{8-24.7/99}
{\rm Ext}^1_{{\cal M}_T({\cal O}_{F, S})}\Bigl(\Q(0), \Q(n)\Bigr)  = 
gr^{\gamma}_nK_{2n-1}({\cal O}_{F, S})\otimes \Q
\end{equation}  
and the higher Ext--groups vanish.  We want to determine the fundamental 
Hopf algebra 
of ${\cal M}_T({\cal O}_{F, S})$.

\begin{lemma} \label{8-24.2/99}
Let ${\cal A}_{\bullet}({\cal O}_{F, S})$ be the maximal 
Hopf subalgebra 
of ${\cal A}_{\bullet}(F):= {\cal A}_{\bullet}({\cal M}_T(F))$ such that 
\begin{equation} \label{8-24.7/99fr}
{\cal A}_{_{1}}({\cal O}_{F, S})  = {\cal O}^*_S \otimes_{\Z} \Q \subset 
{\cal A}_{1}(F)  = F^*\otimes_{\Z} \Q
\end{equation}  
Then ${\cal A}_{\bullet}({\cal O}_{F, S})$ is canonically isomorphic to the 
fundamental Hopf algebra of the category  ${\cal M}_T({\cal O}_{F, S})$. 
\end{lemma}

{\bf Proof}. Observe that 
the Hopf algebra ${\cal A}_{\bullet}({\cal O}_{F, S})$ can  
be defined inductively: 
\begin{equation} \label{5-28.00.1}
{\cal A}_{n}({\cal O}_{F, S}) = \quad \left\{x \in {\cal A}_{n}(F)\quad | \quad \Delta'(x) \in \quad \oplus_{k=1}^{n-1} {\cal A}_{k}({\cal O}_{F, S}) \otimes 
{\cal A}_{n-k}({\cal O}_{F, S}) \right\}
\end{equation}

Consider the subcategory ${\cal M}'_T({\cal O}_{F, S})$ of 
 ${\cal M}_T(F)$ determined 
by the following condition: an object $M$ lies in ${\cal M}'_T({\cal O}_{F, S})$ 
if and only if all its matrix elements are in ${\cal A}_{\bullet}({\cal O}_{F, S})$. 
Clearly ${\cal M}'_T({\cal O}_{F, S})$ 
is a subcategory of ${\cal M}_T({\cal O}_{F, S})$.  
On the other hand it follows from (\ref{5-28.00.1}) that 
any matrix coefficient of an object from ${\cal M}_T({\cal O}_{F, S})$ 
lies in ${\cal A}_{\bullet}({\cal O}_{F, S})$. So 
these two subcategories coincide. The lemma is proved.

 We say that an object $M$ of the category ${\cal M}_{T}(F)$ is defined over 
${\cal O}_{F, S}$ if it belongs to the subcategory ${\cal M}_{T}({\cal O}_{F, S})$. 

\begin{definition}
Let $(v_0, f_n)$ be an $n$-framing on a mixed Tate motive $M$ over $F$. We say that 
the equivalence class of the $n$-framed mixed Tate motive $(M, v_0, f_n)$ 
is defined over ${\cal O}_{F, S}$ if its image in ${\cal A}_{n}(F)$ belongs to 
 the subspace ${\cal A}_{n}({\cal O}_{F, S})$. 
\end{definition}

\begin{proposition} \label{8-25/99}
Let $(M, v_0, f_n)$ be an $n$-framed mixed Tate motive whose equivalence class is defined over 
 ${\cal O}_{F, S}$. Then its minimal representative $\overline M$ is 
defined over ${\cal O}_{F, S}$.
\end{proposition}
Of course $M$ itself is not necessarily defined over 
 ${\cal O}_{F, S}$. 

{\bf Proof}. Let us assume the opposite. Then a certain matrix 
element of $\overline M$ does not belong to ${\cal A}_{\bullet}({\cal O}_{F, S})$. 
Since the minimal representative is a subquotient of every 
object within the  equivalence class of $M$, {\it every} 
motive  from the  equivalence class has such a matrix coefficient. 
Thus no representative in this  equivalence class 
is defined over ${\cal O}_{F, S}$.  
On the other hand by Lemma \ref{8-24.2/99} the Hopf 
algebra ${\cal A}_{\bullet}({\cal O}_{F, S})$ 
is the fundamental Hopf algebra of the category ${\cal M}_T({\cal O}_{F, S})$. 
Thus there must be an $n$-framed object of the latter category 
which is equivalent to the one $(M, v_0, f_n)$. This contradiction proves the 
proposition.

So to produce a mixed Tate motive over ${\cal O}_{F, S}$ it 
is enough to produce 
an element of ${\cal A}_{n}({\cal O}_{F, S})$: then 
we take an $n$-framed mixed 
Tate motive representing this element and its 
minimal representative.

{\bf Remark}. Let us recall  the following rigidity result:
$$
K_{n}({\cal O}_{F, S})\otimes \Q =  
 K_{n}(F)\otimes \Q \quad \mbox{for $n>1$}
$$ 
Combined with (\ref{8-24.7/99}), it implies that for $n>1$ one has 
\begin{equation} \label{8-24.7/99cv}
{\rm Ext}^1_{{\cal M}_T({\cal O}_{F, S})}(\Q(0), \Q(n))  = 
{\rm Ext}^1_{{\cal M}_T(F)}(\Q(0), \Q(n)) 
\end{equation} 
It shows that 
the only difference between the Ext's in the categories of mixed 
Tate motives over $F$ and ${\cal O}_{F, S}$ is the group $Ext^1(\Q(0), \Q(1))$.
 This group is infinite dimensional for the motives over $F$ and finite 
dimensional for the motives over the $S$-integers. As a result the category 
of mixed Tate motives over the 
$S$-integers is ``much smaller'' 
than the one over $F$. For instance  all graded components of 
${\cal A}_{\bullet}(F)$ are infinite dimensional $\Q$-vector spaces, while 
for ${\cal A}_{\bullet}({\cal O}_{F, S})$ they are finite dimensional. 

In the next subsection we show how to apply this remark to obtain estimates from above on 
the $\Q$-vector spaces spanned by periods of mixed Tate motives over ${\cal O}_{F, S}$.

{\bf 5. Periods of mixed Tate motives over ${\cal O}_{F, S}$}. 
Given an embedding $\sigma: F \hra \C$, let us 
define an increasing family of $\Q$-vector subspaces in $\C$, 
depending on the choice of $\sigma$:
$$
\Q = {\cal P}^{\sigma}_{\leq 0}({\cal O}_{F, S}) \hra 
{\cal P}^{\sigma}_{\leq 1}({\cal O}_{F, S}) \hra 
{\cal P}^{\sigma}_{\leq 2}({\cal O}_{F, S}) \hra ... 
$$

Let $(H, v_0, f_n)$ be a $\Q$-Hodge-Tate structure framed by $\Q(0)$ and $\Q(n)$. Let us choose 
a splitting of the weight filtration on $H_{\Q}$. Then there is a period
$
p(H, v_0, f_n) \in \C
$ 
defined as follows. Using the splitting we lift the frame vector $v_0 \in {\rm gr}^W_0H_{\Q}$ 
to a vector  $v_0' \in W_0H_{\Q}$ which projects to $v_0$. Projecting 
$$
W_0H_{\C} \lra W_0H_{\C}/W_{-2n-1}H_{\C} \lra {\rm gr}^{-n}_FH_{\C}
$$ 
we get a vector $v_0''\in {\rm gr}^{-n}_FH_{\C}$ out of $v_0'$. Applying 
the frame functional ${\rm gr}^{-n}_FH_{\C} \stackrel{\sim}{=} {\rm gr}_{-2n}^WH_{\C} \to \C$ to $v_0''$ we get 
the number $p(H, v_0, f_n) \in \C$. 

Let us denote by $r^{\cal H}_{\sigma}$ the 
Hodge realization functor on the category ${\cal M}_T(F)$ corresponding 
to an embedding $\sigma:F \hra \C$.

\begin{definition} \label{6.14.04.1}
The $\Q$-vector space ${\cal P}^{\sigma}_{\leq n}({\cal O}_{F, S})$ is 
spanned over $\Q$ by the periods of $n$-framed Hodge-Tate structures
 $r^{\cal H}_{\sigma}(M, v_0, f_n)$, where $M \in 
{\cal M}_T({\cal O}_{F, S})$, for all possible splittings of the 
weight filtration on $r^{\cal H}_{\sigma}(M)$ and all $n$-framings on $M$.
\end{definition}
We call ${\cal P}^{\sigma}_{\leq n}({\cal O}_{F, S})$ the space 
of {\it weight $\leq n$ periods of mixed Tate motives  over ${\cal O}_{F, S}$}.  
Since the period map on the splitted framed Hodge-Tate structures apparently 
commutes with  the product (i.e. the period of 
the tensor product of the splitted framed Hodge-Tate structures is given by the product of 
the corresponding periods), the space 
$$
{\cal P}^{\sigma}({\cal O}_{F, S}) := \cup_{n\geq 0}{\cal P}^{\sigma}_{\leq n}({\cal O}_{F, S})
$$
is a filtered algebra over $\Q$. Consider its associate graded:
$$
{\cal P}^{\sigma}_{n}({\cal O}_{F, S}):= 
\frac{{\cal P}^{\sigma}_{\leq n}({\cal O}_{F, S})}{{\cal P}^{\sigma}_{\leq n-1}({\cal O}_{F, S})}; 
\qquad {\cal P}^{\sigma}_{\bullet}({\cal O}_{F, S}) := 
\oplus_{n=0}^{\infty}{\cal P}^{\sigma}_{n}({\cal O}_{F, S})
$$

\begin{theorem} \label{6.14.04.5}
There is a surjective homomorphism of commutative algebras
\begin{equation} \label{6.14.04.4}
p_{\sigma}: {\cal A}_{\bullet}({\cal O}_{F, S}) \lra 
{\cal P}^{\sigma}_{\bullet}({\cal O}_{F, S}) 
\end{equation}
\end{theorem}

{\bf Proof}. The crucial claim is that the map $p_{\sigma}$ is well-defined. 
Let us pick a representative $(M, v_0, f_n)$ 
of an element of ${\cal A}_n({\cal O}_{F, S})$. Choosing a splitting 
of the weight filtration of its Hodge reaslization we get its period. 
Changing a splitting we get another period, which differ from the first one by 
an element of ${\cal P}^{\sigma}_{\leq n-1}({\cal O}_{F, S})$. So its projection to 
${\cal P}^{\sigma}_{n}({\cal O}_{F, S})$ is independent of the splitting. 
If $(M', v'_0, f'_n)$ is a subquotient of $(M, v_0, f_n)$ then 
one easily sees that $p_{\sigma}(M', v'_0, f'_n) = p_{\sigma}(M, v_0, f_n)$. 
Finally, the map   $p_{\sigma}$ commutes with the product. The theorem is proved. 

\begin{corollary} \label{6.14.04.6}
One has ${\rm dim}_{\Q}{\cal P}^{\sigma}_{n}({\cal O}_{F, S}) \leq {\rm dim}_{\Q}{\cal A}_{n}({\cal O}_{F, S})$. 
\end{corollary}

Recall  the tensor algebra $T(V_{\bullet})$ generated by the graded $\Q$-vector space 
$V_{\bullet}$. It has a graded Hopf algebra structure, with the commutative 
product given by the shuffle product formula, and the coproduct 
given by the deconcatenation. 
Applying the functor $T$ to the graded space 
$\oplus_{n\geq 1}K_{2n-1}({\cal O}_{F, S})_{\Q}$, where $K_{2n-1}$ 
is in the degree $n$, and using 
 (\ref{8-24.7/99fr}), we get an isomorphism of Hopf algebras
\begin{equation} \label{6.14.04.4df}
{\cal A}_{\bullet}({\cal O}_{F,S}) \stackrel{\sim}{=}
T\Bigl(\oplus_{n\geq 1}K_{2n-1}({\cal O}_{F, S})_{\Q}\Bigr)
\end{equation}
This isomorphism, combined with Corollary \ref{6.14.04.6}, implies 
an estimate from above for ${\rm dim}_{\Q}{\cal P}^{\sigma}_{n}({\cal O}_{F, S})$. 
The following conjecture tells us that this estimate is expected to be exact. 

\begin{conjecture} \label{6.14.04.3}
The map (\ref{6.14.04.4}) is an isomorphism. 
\end{conjecture}

%%% Local Variables: 
%%% mode: latex
%%% TeX-master: "mh"
%%% End: 

\vskip 3mm \noindent
{\bf REFERENCES}
\begin{itemize}
\item[{[BD]}] Beilinson A.A., Deligne P.: 
{\it Motivic polylogarithms and Zagier's conjecture}. 
Manuscript, version of 1992. 
\item[{[BGSV]}] Beilinson A.A., Goncharov A.A., Schechtman V.V., Varchenko A.N.: {\it Aomoto dilogarithms, mixed Hodge structures and motivic cohomology of a pair of triangles in the plane}, the Grothendieck Festschrift, Birkhauser, vol 86, 1990, p. 135-171.
\item[{[Bo]}] Borel A.: {\it Cohomologie de ${\rm SL}\sb{n}$ et valeurs de fonctions z\^eta aux points entiers} Ann. Scuola Norm. Sup. Pisa Cl.
Sci. (4) 4 (1977), no. 4, 613--636.
\item[{[BGK]}] Broadhurst D., Gracey J., Kreimer D.: 
{\it Beyond the triangle and uniqueness relations: non-zeta counterterms of large $N$ at positive knots}. Z. Phys. C75(1997) 559--574. hep-th/9607174.
\item[{[BB]}] Belkale P, Brosnan P.: {\it 
Matroids, motives and a conjecture of Kontsevich}, math.AG/0012198. 
Duke Math. J. 116 (2003). 
\item[{[Chen]}] Chen K.T.: {\it Iterated path integrals}. 
Bull. Amer. Math. Soc. 83 (1977), no. 5, 831--879.
\item[{[CK]}] Connes A., Kreimer D.: {\it Hopf Algebras, Renormalization and Noncommutative Geometry}. hep-th/9808042. Comm. Math. Phys. 199 (1998) 
\item[{[D]}] Deligne P.: {\it Le groupe fondamental de la droite projective 
moins trois points}. In: Galois groups over $\Q$. 
Publ. MSRI, no. 16 (1989) 79-298.   
\item[{[D2]}] Deligne P.: {\it Cat\'egories tannakiennes}. 
 The Grothendieck Festschrift, Vol. II,
111--195, Progr. Math., 87, Birkhauser Boston, Boston, MA, 1990.
\item[{[DG]}] Deligne P., Goncharov A.B.:  {\it Groupes fondamentaux 
motiviques de Tate mixte}.   math.NT/0302267. To appear in Ann. Sci. ENS. 
\item[{[Dr]}] Drinfeld V.G.: {\it On quasi-triangular quasi-Hopf algebras and some group related to 
 Gal$(\overline{\Q}/\Q)$}. Leningrad Math. Journal, 1991. (In Russian).
\item[{[G0]}] Goncharov A.B.: {\it Multiple $\zeta$-numbers,
    hyperlogarithms and mixed Tate motives}. Preprint MSRI 058-93, June 1993.  
\item[{[G1]}] Goncharov A.B.: {\it The dihedral Lie algebras 
and Galois symmetries of 
$\pi_1^{(l)}({\Bbb P}^1 - \{0, \infty\} \cup \mu_N)$}. 
Duke Math. J.  vol 110, N3, (2001), pp. 397-487. 
math.AG/0009121. 
%See also a short version at
%www.math.uiuc.edu/K-theory/ Dec. 1998.
\item[{[G2]}] Goncharov A.B.: {\it Multiple $\zeta$-values, Galois groups and geometry of 
modular varieties} Proc. of the Third European Congress of
Mathematicians. Progress in Mathematics, 
vol. 201, p. 361-392. Birkhauser Verlag. (2001)
math.AG/0005069.
\item[{[G3]}] Goncharov A.B.: {\it Multiple polylogarithms and mixed
    Tate motives}. math.AG/0103059. 
\item[{[G4]}] Goncharov A.B.: {\it Polylogarithms and motivic Galois groups}, 
Motives (Seattle, WA, 1991), 43--96, Proc. Sympos. Pure Math., 55, Part 2,
Amer. Math. Soc., Providence, RI, 1994.
\item[{[G5]}] Goncharov A.B.: {\it Volumes of hyperbolic manifolds and mixed Tate motives}. JAMS 12 (1999),  n 2, 569-618. math.AG/9601021.
\item[{[G6]}] Goncharov A.B.: {\it Polylogarithms in arithmetic and 
geometry}. Proceedings of the International Congress of Mathematicians, Vol. 1 
 (Z\"urich, 1994), 374--387, Birkh\"auser, Basel, 1995.
\item[{[G7]}] Goncharov A.B.: {\it Mixed elliptic motives}. 
Galois representations 
in arithmetic algebraic geometry, (Durham 1996), London Math. Soc. Lecture 
Note Ser., vol. 254, Cambridge University Press, 
1998. 147-221. 
\item[{[G8]}] Goncharov A.B.: {\it Periods and mixed Tate motives}.  
math.AG/0202154. 
\item[{[G9]}] Goncharov A.B.: 
{\it Multiple polylogarithms, cyclotomy and modular complexes},     Math. Res. Letters, 
 vol. 5. (1998), pp. 497-516.  
www.math.uiuc.edu/K-theory/ N 297.
\item[{[K]}] Kreimer D.: {\it On the Hopf algebra structure of 
perturbative quantum field theories}. q-alg/9707029. 
Adv. Theor. Math. Phys. 2 (1998). 
\item[{[L1]}] Levine, M: {\it Tate motives and the vanishing 
conjectures for algebraic 
$K$-theory}. Algebraic $K$-theory and algebraic topology (Lake
Louise, AB, 1991), 167--188, NATO Adv. Sci. Inst. Ser. 
C Math. Phys. Sci., 407, 
Kluwer Acad. Publ., Dordrecht, 1993. 
\item[{[Lo1]}] Loday, J.-L.: {\it Arithmetree}. math.CO/0112034,  J. of Algebra 258, 2002
\item[{[Lo2]}] Loday, J.-L.: {\it Dialgebras}.  math.QA/0102053, 
Springer LNM 1763, 2001. 
\item[{[Z]}] Zagier D.:{\it Periods of modular forms, traces of 
Hecke operators, and multiple $\zeta$-values} Research into automorphic 
forms and $L$-functions (Japanese), Surikaisekikenkyusho Kokyuroku No. 843 (1993), 162-170.. 
\end{itemize}

Address: Dept. of Mathematics, Brown University, Providence, RI 02912, USA.

e-mail sasha@math.brown.edu

\end{document}